\documentclass[11pt]{amsart}

\usepackage[colorlinks = true,
            linkcolor = blue,
            urlcolor  = blue,
            citecolor = black,
            anchorcolor = blue]{hyperref}
\usepackage{color}
\usepackage[shortalphabetic,initials,nobysame]{amsrefs}
\usepackage{upgreek}
\usepackage{tikz}
\usepackage[applemac]{inputenc} 
\usetikzlibrary{arrows}
\usepackage{xcolor}
\usepackage[all]{xy}
\usepackage{graphicx}
\usepackage{wrapfig}
\usepackage{yfonts}
\usepackage{graphicx}
\usepackage{amssymb}
\usepackage{epstopdf}
\usepackage{mathrsfs}
\usepackage[mathcal]{eucal}
\usepackage{bm}

\usepackage{slashed}
\renewcommand{\eprint}[1]{\href{https://arxiv.org/abs/#1}{#1}}

\DeclareMathOperator{\SL}{\mathrm{SL}}

\DeclareMathOperator{\diag}{diag}
\DeclareMathOperator{\Hom}{Hom}

\newcommand{\hcL}{\hat{\cL}}

\BibSpec{article}{%
    +{}  {\PrintAuthors}                {author}
    +{,} { \textit}                     {title}
     +{,} {}                            {note}
    +{.} { }                            {part}
    +{:} { \textit}                     {subtitle}
    +{,} { \PrintContributions}         {contribution}
    +{.} { \PrintPartials}              {partial}
    +{,} { }                            {journal}
    +{}  { \textbf}                     {volume}
    +{}  { \PrintDate}                {date}
    +{,} { \issuetext}                  {number}
    +{,} { \eprintpages}                {pages}
    +{,} { }                            {status}
    +{,} { \DOI}                   {doi}
    +{,} { \eprint}        {eprint}
      +{,} {\publisher}                {publisher}
    +{,} { \address}                   {address}
    +{}  { \parenthesize}               {language}
    +{}  { \PrintTranslation}           {translation}
    +{;} { \PrintReprint}               {reprint}
    +{.} {}                             {transition}
    +{}  {\SentenceSpace \PrintReviews} {review}
}

\newtheorem{Thm}{Theorem}[section]

\newtheorem{Prop}[Thm]{Proposition}

\theoremstyle{definition}

\theoremstyle{remark}
\newtheorem{Rem}[Thm]{Remark}

\newtheoremstyle{named}{}{}{\itshape}{}{\bfseries}{.}{.5em}{#1 #3}
\theoremstyle{named}

\def\Q{\mathbb{Q}}
\def\N{\mathbb{N}}
\def\C{\mathbb{C}}
\def\Z{\mathbb{Z}}

\def\P{\mathbb{P}}

\def\M{\mathbb{M}}

\def\fsh{{\mathcal U}_{\hbar}(\widehat{\frak{sl}}_{2}) }

\def\g{\mathfrak{g}}
\def\Frenkel:2013uda{\mathfrak{h}}

\def\cD{\mathcal{D}}

\def\cF{\mathcal{F}}

\def\cL{\mathcal{L}}

\def\cO{\mathcal{O}}

\def\cV{\mathcal{V}}
\def\cW{\mathcal{W}}

\def\ze{\zeta}

\def\h{\theta}

\def\bA{\textbf{A}}

\def\bo{\textbf{o}}

\def\=>{\Longrightarrow}

\def\to{\longrightarrow}

\def\o+{\oplus}
\def\bo+{\bigoplus}

\def\<{\langle}
\def\>{\rangle}
\def\({\left(}
\def\){\right)}

\def\^{\wedge}
\def\+{\dagger}

\def\inv{^{-1}}

\def\dd[#1,#2]{\frac{d#1}{d#2}}
\def\del[#1,#2]{\frac{\partial #1}{\partial #2}}
\def\over[#1]{\overline{#1}}
\def\vec[#1]{\overrightarrow{#1}}

\makeatletter
\def\mr@ignsp#1 {\ifx\:#1\@empty\else #1\expandafter\mr@ignsp\fi}%
\newcommand{\multiref}[1]{\begingroup
\xdef\mr@no@sparg{\expandafter\mr@ignsp#1 \: }%
\def\mr@comma{}%
\@for\mr@refs:=\mr@no@sparg\do{\mr@comma\def\mr@comma{,}\ref{\mr@refs}}%
\endgroup}
\makeatother

\newcommand{\hypref}[2]{\ifx\href\asklFrenkel:2013udaas #2\else\href{#1}{#2}\fi}

\tikzset{->-/.style={decoration={
  markings,
  mark=at position .5 with {\arrow{latex}}},postaction={decorate}}}
\tikzset{
    >=latex
    }

\newcommand{\ard}{{\mathbf{d}}}

\def\M{\textbf{M}}
\def\N{\textbf{N}}
\def\h{\hbar}
\def\fp{ {\textbf{p}}  }

\newcommand{\wt}{\widetilde}

\newcommand{\nc}{\newcommand}
\nc{\on}{\operatorname}
\nc{\la}{\lambda}
\nc{\wh}{\widehat}
\nc{\ghat}{\wh\g}
\nc{\mb}{\mathbf}

\begin{document}
\title[Geometric realizations of the Bethe ansatz equations]{Geometric realizations of the Bethe ansatz equations}

\author[A.M. Zeitlin]{Anton M. Zeitlin}
\address{\hspace{-0.53cm} School of Mathematics,\newline
Georgia Institute of Technology,\newline
 686 Cherry Street, \newline
Atlanta, GA 30332-0160, USA
\newline
Email: \href{mailto:zeitlin@gatech.edu}{zeitlin@gatech.edu},\newline
 \href{https://zeitlin.math.gatech.edu}{https://zeitlin.math.gatech.edu}}

\date{\today}

\numberwithin{equation}{section}

\begin{abstract}
These lecture notes are devoted to the recent progress in the geometric aspects of quantum integrable systems based on quantum groups solved using the Bethe ansatz technique. One part is devoted to their enumerative geometry realization through the quantum K-theory of Nakajima quiver varieties. The other part describes a recently studied $q$-deformation of the correspondence between oper connections and Gaudin models.
The notes are based on a minicourse at C.I.M.E. Summer School ``Enumerative geometry, quantisation and moduli spaces," September 04-08, 2023.    
\end{abstract}

\maketitle

\tableofcontents

\section{Introduction: geometric wonders of integrable systems}

The theory of integrable systems, both classical and quantum, played an important role in various parts of mathematics and physics for the last 50 years.

What is a classical integrable system? Consider a symplectic $2n$-manifold and $n$ mutually commuting functions under the corresponding Poisson bracket; it is a Hamiltonian system, where one of these commuting quantities is chosen as a Hamiltonian, while $n-1$ functions remain preserved along the Hamiltonian flow as a simple consequence of Poisson commutativity. 
The dynamics of the corresponding system thus belong within the level set of these functions. The critical result about the dynamics of the integrable system is the  {\it Liouville-Arnold theorem} \cite{arnold78}, which states that this level set is a collection of tori and there exists a set of coordinates, known as action-angle variables in which simple linear equations describe the dynamics of the corresponding integrable system. The complicated part is finding these coordinates, which could be a nontrivial task. 

Initially, integrable systems were viewed as rare cases of classical mechanical systems, such as harmonic oscillators and tops. However, in the 1960s, the revolution came through the so-called soliton theory. The crucial example was the Korteweg-de Vries (KdV) partial derivative equation, well-known originally in the theory of shallow waves and later in plasma physics. The discovery of Green, Gardner, Kruskal, and Miura \cite{GGKM} related this  equations' solutions to the Sturm-Liouville operator's spectral data. Then, the work of Lax gave this differential equation the so-called Lax pair formulation: $\dot{L}=[L,A]$, where $L$ is a Sturm-Liouville operator and $A$ is a differential operator as well. That allowed Zakharov and Faddeev \cite{ZahFAd} to show that it is an infinite-dimensional integrable system and relate the spectral data of $L$-operator to action-angle variables, while the corresponding infinite family of mutually commuting Hamiltonians to the monodromy matrix, forming a KdV {\it hierarchy} of integrable partial differential equations.

These findings led to an avalanche of results because a lot of known nonlinear equations admitted Lax pair formulation, such as sine-Gordon, nonlinear Schr\"odinger, and many others (see, e.g.,  \cite{Novikov:1984}). 
On the one hand, that gave a method of finding solutions through the spectral data of $L$-operator, now known as a {\it classical inverse scattering method} \cite{faddeev1987hamiltonian}. On the other hand, the algebraic structures behind the Lax pair and related integrability brought an incredible amount of Lie-theoretic methods, which helped to discover other integrable hierarchies. Some of that technique culminated in the work of Drinfeld and Sokolov \cite{Drinfeld:1985}, where an analog of the KdV model was constructed for any affine Lie algebra $\widehat{\mathfrak{g}}$, while the original KdV corresponds to $\mathfrak{g}=\mathfrak{sl}_2$. At the same time, a lot of new finite-dimensional integrable models (multiparticle systems) were discovered through Lie-theoretic methods, e.g., Toda lattice, Rujsenaars-Schneider and Calogero-Moser systems, etc (see, e.g.,\cite{babelon}).

In the late 70s and 80s, the quantization of these classical integrable models became of interest through the prism of the quantum field theory discoveries. In the quantum world, the mutually commuting conserved quantities under the Poisson bracket are replaced by a mutually commuting family of operators,  quantum Hamiltonians, acting in the Hilbert space. The problem of finding the action-angle variables is replaced by the problem of finding the {\it spectrum} of the model, i.e., simultaneous diagonalization of these operators.  
Although some hesitancy was initially thrown in by some physicists regarding whether integrability would be preserved through quantization, a lot of infinite-dimensional integrable models were successfully quantized \cite{Korepin_1993}: the fundamental idea in this process was to put them on a lattice. At the same time, the lattice models of statistical physics, such as Heisenberg spin chain \cite{Heisenberg:1928aa} were already well-known and studied, producing algebraic structures such as R-matrices and Yang-Baxter equation \cite{Baxter:1982zz}. The increased interest in the understanding of the mathematics behind these structures lead to the discovery of quantum groups by Drinfeld and Jimbo \cite{Drinfeld:1986in}, \cite{Jimbo1985q}, \cite{Chari:1994pz}. The algebraic method for diagonalization of quantum Hamiltonians is known as {\it algebraic Bethe ansatz} \cite{Takhtajan:1979iv}, \cite{Faddeev:aa}, \cite{Korepin_1993}, \cite{Reshetikhin:2010si}. 
Using this method, the quantum Hamiltonians generating and abelian algebra called {\it Bethe algebra} has its spectrum given by the symmetric functions of Bethe roots: roots of an algebraic system of equations called {\it Bethe equations}. Section 2 describes algebraic Bethe ansatz applied to the simplest nontrivial integrable model, called XXZ Heisenberg spin chain, which is based on quantum group $\fsh$.

While all of that progress was happening in the theory of quantum integrable models, through the late 80s to early 90s, in the works of Dubrovin, Givental, Kontsevich, Witten, and many others, classical integrable models penetrated geometry, particularly enumerative geometry, aided by the development of String Theory \cite{Maninfrob}, \cite{Guest}. Among the key results, one may mention Witten's conjecture \cite{witten1990two}, proved by Kontsevich \cite{kontsevich1992intersection},  relating intersection numbers on moduli spaces to the so-called $\tau$-function for KdV hierarchy; on the other hand result of Givental and Kim \cite{givental1995}, \cite{1996alg.geom7001K} identified the equivariant quantum cohomology ring of the flag variety with the algebra of functions on the Lagrangian subvariety generated by Toda lattice Hamiltonians.

Development of the infinite-dimensional Lie algebras representation theory in the 80s and 90s, particularly affine algebras, led to important discoveries. Within the study of conformal field theory (with integrable deformations),  it was observed that matrix elements of the products of the intertwiners (conformal blocks) of affine algebras and their quantum group deformations satisfy differential and difference equations known as (quantum) Knizhnik-Zamolodchikov equations \cite{etingof1998lectures}. The quantum version of the algebraic structures discovered by Drinfeld and Sokolov in the context of KdV-type models produce new objects called $W$-algebras, whose conformal blocks are represented as certain hypergeometric functions are governed by certain differential /difference equations as well. The correspondence between these conformal blocks is now known as quantum Langlands correspondence \cite{Aganagic:2017smx}.  
The early works of Feigin, Frenkel, and collaborators \cite{Feigin:1994in}, \cite{Frenkel:2003qx}, \cite{Frenkel:2004qy} lead to a crucial relation between the spectrum of quantum integrable models of Gaudin type, which is based not on the quantum group, but just any simple Lie algebra $\mathfrak{g}$, and flat $^LG$-connections with singularieties on a projective line called opers, where $^LG$ is a simple group of adjoint type with Langlands dual algebra to $\mathfrak{g}$.  This correspondence is an example of geometric Langlands correspondence  \cite{frenkel2007langlands} and is a {\it critical limit} of quantum Langlands correspondence, where the eigenvalue problem for Gaudin model comes as the limit of KZ equations and the differential operator corresponding to oper connection is the limit of the operator from the differential equation satisfied by conformal blocks of W-algebras.

At the same time, geometric representation theory \cite{CG} gained steam. In the works of Nakajima \cite{nakajima1998, Nakajima:1999hilb, Nakajima:2001qg}, Schiffmann \cite{olivier1998quantum}, Varognolo and Vasserot  \cite{Vasserot:wo}, \cite{varagnolo2000quiver} the realization of (affine) quantum group actions and representations were realized via cohomology and K-theory of a particular class of algebraic varieties called symplectic resolutions.
A very special subclass of such varieties is called Nakajima quiver varieties (see, e.g.,\cite{Ginzburg:}, \cite{Kirquiv}), built via oriented graph called quiver in some way generalizing the notion of the Dynkin diagram in these geometric realizations.

The observations of Nekrasov and Shatashvili in 2009 \cite{Nekrasov:2009ui}, \cite{Nekrasov:2009uh}, who studied 2d and 3d supersymmetric gauge theories, led to the emergence of quantum integrable models in the geometric context, in particular in the context of quiver varieties. Their conjectures and hints led to the work of Okounkov and collaborators \cite{Braverman:2010ei}, \cite{2012arXiv1211.1287M}, \cite{Okounkov:2015aa}, \cite{Okounkov:2016sya} studying the enumerative geometry of symplectic resolutions. The partition functions of gauge theories from \cite{Nekrasov:2009ui}, \cite{Nekrasov:2009uh}  were interpreted as the so-called {\it vertex functions},  the counting functions of {\it quasimaps} to Nakajima varieties. Most importantly, vertex functions appear as solutions of difference equations called quantum Knizhnik-Zamolodchikov equations, which govern intertwining operators for quantum affine algebras.
The isomorphism between quantum cohomology/quantum K-theory rings of Nakajima varieties with Bethe algebras of the corresponding integrable model was  
one of the central conjectures in original papers \cite{Nekrasov:2009ui}, \cite{Nekrasov:2009uh}. It was proven in \cite{Pushkar:2016qvw} for ${\rm T^*Gr}_{k,n}$ and in \cite{Koroteev:2017aa} for Nakajima varieties of type $A$. We will discuss this proof, as well as the related geometric structures, in Section 3. This provides one of the geometric realizations of Bethe equations, as the relations for the quantum K-theory ring are generated by the symmetric functions of Bethe equations: eigenvalues of the elements of Bethe algebra.   

Some elements in the Bethe algebra have a spectrum described by elementary symmetric functions of Bethe roots; they are called {\it Baxter Q-operators}, which generate the entire Bethe algebra. Originally introduced ad-hoc, the proper representation theoretic meaning of those operators took some time:  Bazhanov, Lukyanov, and Zamolodchikov, who studied the Q-operators in the context of quantum KdV model in the 90s, observed that the needed representations belong to the extension of the Grothendieck ring of finite-dimensional representations of the affine quantum group $\fsh$. The complete understanding of that extension was only recently understood and formulated by Frenkel and Hernandez \cite{Frenkel:2013uda}, \cite{Frenkel:2016} for any affine algebra $\widehat{\mathfrak{g}}$. 
These $Q$-operators satisfy the functional relations called the {\it $QQ$-system} which in original work on quantum KdV emerged as an equation for the spectral determinant of an eigenvalue problem for ODE with singular potential (affine oper), which is now known as {\it ODE/IM correspondence} \cite{Dorey:2007zx}, \cite{MRV1}, \cite{MRV2}. For polynomial solutions {\it $QQ$-relation} were known (under certain nondegeneracy conditions) to be equivalent to Bethe equations. 
It turns out that geometrically, the QQ-system is the proper intermediate object between the deformation of oper connections on the projective 
line--multiplicative group-like objects and spin chain models based on quantum groups. Thus, one obtains the deformation \cite{KSZ}, \cite{Frenkel:2020}, \cite{KoroteevZeitlinCrelle} of the example of the geometric Langlands correspondence discussed above and provides another geometrization of Bethe equations. That will be the subject of Section 4, where we also relate it to the enumerative structures in Section 3 and describe some applications.\\

\vspace*{3mm}

\noindent {\bf Acknowledgements}
A.M.Z. is grateful to various collaborations he had over the past several years while working on various geometric aspects of integrable systems, in particular E. Frenkel, P. Koroteev, D.S. Sage, A.V. Smirnov. 
A.M.Z. would like to thank the organizers of the Summer school  ``Enumerative geometry, quantisation and moduli spaces" and CIME foundation for the for the opportunity to teach this minicourse and wonderful hospitality.  A.M.Z. was also partially supported by Simons
Collaboration Grant 578501 and NSF grant
DMS-2203823.

\section{Algebraic Bethe ansatz and Baxter $Q$-operators} 
\subsection{Quantum Integrable systems of spin chain type}
Let $\mathfrak{g}$ is a complex simple Lie algebra. In this section we describe  the  mathematical perspective on how algebraic Bethe ansatz works for integrable models of specific type, namely the so-called spin chains. Let $\mathfrak{g}$ be a simple Lie algebra and 
$\hat{\mathfrak{g}}_{k=0}=\mathfrak{g}[t^{\pm 1}]$ be the corresponding loop algebra: affine algebra with vanishing central charge $k=0$. The finite-dimensional modules $\{V_i\}$ of $\mathfrak{g}$ give rise to the so-called evaluation modules $\{V_i(a_i)\}$,  where $a_i\in\mathbb{C}^\times$ stand for the value of the loop parameter $t$. These modules generate a tensor category, namely every finite-dimensional representation of $\hat{\mathfrak{g}}$ can be written as a tensor product of evaluation modules.  
Passing from $\mathfrak{g}[t^{\pm 1}]$ to the corresponding quantum affine algebra 
$U_{\hbar}(\hat{\mathfrak{g}})$ or the Yangian $Y_{\hbar}({\mathfrak{g}})$ 
one obtains a deformation of such tensor category, known as braided tensor category \cite{Chari:1994pz}. This object features a new  intertwining operator known as R-matrix:
\begin{equation}
R_{V_i(a_i),V_j(a_j)}:V_i(a_i)\otimes V_j(a_j)\to V_j(a_j)\otimes V_i(a_i)
\end{equation}
satisfying famous Yang-Baxter equation.  We note, that in the deformed case the analogue of evaluation map exists only in the type $A$ and the $\{V_i(a_i)\}$ stand here for appropriately ``twisted" by $a_i$ finite-dimensional representations of $U_{\hbar}(\hat{\mathfrak{g}})$.

To describe the integrable model, we choose a specific object in such braided tensor category 
$$\mathscr{H}=V_{i_1}(a_{i_1})\otimes \dots \otimes V_{i_n}(a_{i_n}),$$ 
which we refer to as {\it physical space}, the vectors in this space are called {\it states}. 
For a given module $W(u)$ called {\it auxiliary module} with parameter $u$  known as {\it spectral parameter}, we define the {\it transfer matrix} 
\begin{equation}
T_{W(u)}={\rm Tr}_{W(u)}\Big[P{R}_{\mathscr{H}, W(u)}(1\otimes Z)\Big]; \quad T_{W(u)}:\mathscr{H}\to \mathscr{H}.
\end{equation}
Here the {\it twist Z} is given by $Z=\prod^r_{i=1}\zeta_i^{\check{\alpha}_i}\in e^{\mathfrak{h}}$, where $\mathfrak{h}$ is the Cartan subalgebra in $\mathfrak{g}$, $\{\check{\alpha}_i\}_{i=1,\dots, r}$ are the simple coroots of $\mathfrak{g}$, and  $P$ is a permutation operator. The object standing under the trace, namely $M^Z_{W(u)}=P{R}_{\mathscr{H}, W(u)}(1\otimes Z)$ is called a {\textit quantum monodromy matrix}: it is an operator in  $\mathscr{H}\otimes W(u)$. Therefore, the transfer matrix $T_{W(u)}$ is an operator acting on the physical space $\mathscr{H}$.  
The Yang-Baxter equation implies that transfer matrices, corresponding to various choices of $W(u)$ form a commutative algebra, known as {\it Bethe algebra}. The commutativity of Bethe algebra implies {\it integrability} and the expansion coefficients of the transfer matrix yield (nonlocal) {\it Hamiltonians} of the XXX or XXZ spin chain depending on whether we deal with the Yangian or the quantum affine algebra. From now on, we will primarily focus on quantum affine algebra and the XXZ model, although most of the construction below applies to the Yangian and the XXX models as well. 

The classic example of the XXZ Heisenberg magnet corresponds to the quantum algebra $\fsh$ in which the physical space $\mathcal{H}$ is constructed from $V_i(a_i)={\rm V}(a_i)$ -- the standard two-dimensional evaluation modules of $\fsh$.

The solution of the integrable model implies finding the eigenvalues and eigenvectors of simultaneously diagonalized Hamiltonians, i.e. elements of the Bethe algebra. 
One way to accomplish the task is to follow the old-fashioned procedure from the 1980s known as {\it algebraic Bethe ansatz} \cite{Korepin_1993}, \cite{Faddeev:aa},\cite{Reshetikhin:2010si}.  
It implies that the eigenvalues of the transfer-matrices (upon rescaling) are symmetric functions of the roots of
the system of algebraic equations, known as Bethe ansatz equations. Although this approach is straightforward and effective, we will discuss other modern techniques, which provide insights into representation-theoretic aspects of the problem.

\subsection{The quantum algebra $\fsh$}
Let us remind the structure of the quantum affine algebra $\fsh$; we use the following normalization: $\hbar={\mathrm q}^2$, where $\mathrm{q}$ is a standard deformation parameter \cite{Chari:1994pz}.

The quantum affine algebra $U_{\hbar}(\hat{\mathfrak{sl}}_2)$ is the Hopf algebra over 
$\mathbb{C}(\h^{1/2})$ generated by Chevalley generators
$e_{\pm \alpha_1}$, $e_{\pm \alpha_0}$, 
$k_{\alpha_1}={\h}^{h_{\alpha_1}/2}$, $k_{\alpha_0}={\h}^{h_{\alpha_0}/2}$ with the following relations:
\begin{eqnarray}
k_{\gamma}k_{\gamma}^{-1}&\!\!=\!\!&
k_{\gamma}^{-1}k_{\gamma}^{}=1~,\qquad\quad\;
[k_{\gamma}^{\pm 1},k_{\gamma'}^{\pm 1}]=0~,\qquad
\nonumber\\
k_{\gamma}^{}e_{\pm\alpha}^{}k_{\gamma}^{-1}&\!\!=\!\!&
{\h}^{\pm(\gamma,\alpha)/2}e_{\pm\alpha}~,
\qquad k_{\gamma}^{}e_{\pm\alpha_0}k_{\gamma}^{-1}=
{\h}^{\pm(\gamma,\alpha_0)/2}e_{\pm(\alpha_0)},\nonumber
\end{eqnarray}
\begin{eqnarray}
[e_{\alpha_1},e_{-\alpha_0}]&\!\!=\!\!&0~,\qquad\qquad\qquad\quad
[e_{-\alpha_1},e_{\alpha_0}]=0~,
\\[7pt]
[e_{\alpha},e_{-\alpha}]&\!\!=\!\!&[h_{\alpha}]_{{\h}}~,\qquad\qquad\;\;\,
[e_{\alpha_0},e_{-\alpha_0}]=[h_{\alpha_0}]_{{\h}}~,\nonumber
\end{eqnarray}
\begin{eqnarray}\nonumber
[e_{\pm\alpha_1}^{},[e_{\pm\alpha_1}^{},[e_{\pm\alpha_1}^{},
e_{\pm\alpha_0}^{}]_{{\h}}]_{{\h}}]_{{\h}}&\!\!=\!\!&0~,\quad\;\;
\\\nonumber
[[[e_{\pm\alpha_1},e_{\pm(\alpha_0)}]_{{\h}},e_{\pm\alpha_0}]_{{\h}},
e_{\pm\alpha_0}]_{\h}&\!\!=\!\!&0~,\quad\;\;
\end{eqnarray}
where $\gamma~=~\alpha_1,\alpha_0$,  and we assume for simplicity that $(\alpha_1, \alpha_1)=(\alpha_0, \alpha_0)=-(\alpha_0, \alpha_1)=2$ and
$$[h_\beta]_{{\h}}\!:=\!(k_\beta\!-\!k_\beta^{-\!1})/({\h}^{1/2}-\!{\h}^{-1/2}).$$
The bracket $[\cdot,\cdot]_{\h}$ are the ${\h}$-commutator:
$
[e_{\beta}^{},e_{\beta'}^{}]_{{\h}}\!\!=\!\!e_{\beta}^{}e_{\beta'}^{}-
{\h}^{(\beta,\beta')/2}e_{\beta'}^{}e_{\beta}^{}~.
$

The Hopf algebra structure on $\mathcal{U}_{\h}(\widehat{\mathfrak{sl}}_2)$ is given by the following coproduct $\Delta_{\sqrt{\h}}$ and antipode $S$:
\begin{eqnarray}
\begin{array}{rcccl}
\Delta(k_\gamma^{\pm 1})&\!\!=\!\!&k_\gamma^{\pm 1}\otimes k_\gamma^{\pm 1}~,
\qquad\qquad\qquad S(k_\gamma^{\pm 1})&\!\!=\!\!&k_\gamma^{\mp 1}~,
\\[7pt]
\Delta(e_{\beta}^{})&\!\!=\!\!&e_{\beta}^{}\otimes 1
+ k_{\beta}^{-1}\otimes e_{\beta}^{}~,
\qquad\quad S(e_{\beta}^{})&\!\!=\!\!&-k_{\beta}e_{\beta}^{}~,
\\[7pt]
\Delta(e_{-\beta}^{})\!\!&=\!\!&
e_{-\beta}^{}\otimes k_{\beta}+1 \otimes e_{-\beta}^{}~,
\qquad S(e_{-\beta}^{})&\!\!=\!\!&-e_{-\beta}^{}k_{\beta}^{-1}~,
\end{array}
\end{eqnarray}
where $\beta=\alpha_1,\,\alpha_0$.\\

We are interested in the case when central charge is zero: in this case $k_{\alpha_1}=k^{-1}_{\alpha_0}$.

At the same time, for geometric applications Drinfeld realization is more natural. In this case 
$\fsh$ is generated by elements $\{E_n, H_n, F_n, K\}_{n\in \mathbb{Z}}\setminus \{0\}$ over $\C(\hbar^{1/2})$, satisfying the following relations (see, e.g.,\cite{Chari:1994pz},\cite{Khoroshkindrinf}):
\begin{eqnarray}
\begin{array}{l} \label{sl2rel}
	K K^{-1}=K^{-1} K=1,\\
	{[}H_{m},H_n{]}={[}H_m,K^{\pm 1}{]}=0,\\
	K E_{m} K^{-1}=\hbar E_{m}, K F_{m} K^{-1}=\hbar^{-1 } F_{m},\\
	{[}E_m, F_l{]}=\dfrac{\psi^{+}_{m+l}-\psi^{-}_{m+l}}{\h^{1/2}-\h^{-1/2}}, \\
	{[} H_{k}, E_{l} {]}=\dfrac{{[}2k{]}_{\h}}{k} E_{k+l}, \quad
	{[} H_{k}, F_{l} {]}=-\dfrac{{[}2k{]}_{\h}}{k} F_{k+l},
\end{array}
\end{eqnarray}
where  
$$
\begin{array}{l}
\sum\limits_{m=0}^{\infty} \psi^{+}_m z^{-m}= K \exp\Big( (\h^{1/2}-\h^{-1/2})\sum\limits_{k=1}^{\infty} H_{k} z^{-k} \Big)\\
\sum\limits_{m=0}^{\infty} \psi^{-}_{-m} z^{m}= K^{-1} \exp\Big( -(\h^{1/2}-\h^{-1/2})\sum\limits_{k=1}^{\infty} H_{-k} z^{k} \Big).
\end{array}
$$

The following formulas allow to restore Chevalley generators in terms of Drinfeld generators with zero central charge:
\begin{eqnarray} \label{drinf}
&&e_{\alpha_0}=F_0K^{-1},\quad e_{\alpha_1}=E_{-1},\quad
e_{-\alpha_0}= K E_0,\quad
e_{-\alpha_1}=F_{1} ,\\
&&k_{\alpha_1}=k_{\alpha_0}^{-1}=K.\nonumber
\end{eqnarray}

The universal R-matrix \cite{Chari:1994pz}, \cite{KhoroshkinR} is an element in the completed tensor product $\mathcal{U}_{\hbar}(\hat{\mathfrak{b}}_{+})\hat{\otimes }\mathcal{U}_{\hbar}(\hat{\mathfrak{b}}_{-})$. 
Here $\mathcal{U}_{\hbar}(\hat{\mathfrak{b}}_{+})$, $\mathcal{U}_{\hbar}(\hat{\mathfrak{b}}_{-})$ are the upper and lower Borel subalgebras of $\mathcal{U}_{{\h}}(\widehat{\mathfrak{sl}}_2)$, generated by $k^{\pm 1}_{\alpha_0}, k^{\pm 1}_{\alpha_1}, 
e_{ \alpha_0}, e_{ \alpha_1}$ and $k^{\pm 1}_{\alpha_0}, k^{\pm 1}_{\alpha_1}, 
e_{ -\alpha_0}, e_{ -\alpha_1}$ correspondingly. 
Considering its action in evaluation modules of $\fsh$ leads to R-matrix intertwining operators we discussed in the previous section. 
The R-matrix satisfies the following relations with respect to the coproduct $\Delta$ and opposite coproduct $\tilde{\Delta}=\sigma\Delta$ and $\sigma(a \otimes b)=b\otimes a$:
\begin{eqnarray}
&&\tilde{\Delta}(a)\!\!=\!\!R\Delta(a)R^{-1} \qquad\quad\;\;
\forall\,\,a \in \fsh~,
\nonumber\\
&&(\Delta\otimes{\rm id})R\!\!=\!\!R^{13}R^{23}~,\qquad
({\rm id}\otimes\Delta R=R^{13}R^{12}.
\end{eqnarray}
The relations above can be understood in the following way:
$R^{12}=\sum a_{i}\otimes b_{i}\otimes{\rm id}$,
$R^{13}=\sum a_{i}\otimes{\rm id}\otimes b_{i}$,
$R^{23}=\sum {\rm id}\otimes a_{i}\otimes b_{i}$ if $R$ has the form
$R=\sum a_{i}\otimes b_{i}$. The celebrated Yang-Baxter equation

\subsection{Bethe ansatz for XXZ spin chain}

Evaluation homomorphism $\pi(x)$ for $\fsh$ is as follows:
\begin{eqnarray}
 \pi(x): \left\{
    \begin{array}{rl}
      &h_{\alpha_1} \rightarrow  \mathcal{H},\quad h_{\alpha_0}\rightarrow - \mathcal{H},\\
      &e_{-\alpha_1}  \rightarrow \mathcal{F},\quad e_{\alpha_1}\rightarrow \mathcal{E}\\
      &e_{-\alpha_0} \rightarrow x\mathcal{E},\quad e_{\alpha_0}\rightarrow x^{-1}\mathcal{F}
    \end{array} \right.
\end{eqnarray}

where the generators $\mathcal{E}, \mathcal{F}, \mathcal{H}$ of $\mathcal{U}_{{\h}}(\mathfrak{sl}_2)$ satisfy standard commutation relations:
$$
[\mathcal{E}, \mathcal{F}]=\frac{\hbar^{\frac{\mathcal{H}}{2}}-\hbar^{-\frac{\mathcal{H}}{2}}}{\hbar^{\frac{1}{2}}-\hbar^{-\frac{1}{2}}}.
$$
Under the evaluation homomorphism one can produce the repesentation of 
$\fsh$ out of representation of $\mathcal{U}_{{\h}}(\mathfrak{sl}_2)$. 
The representations of special interest are Verma modules $V_n^+=\{{\rm span}(v_k=\mathcal{F}^kv_0);~\mathcal{E}v_0=0, \mathcal{H}v_0=nv_0\}$,  
where the action of generators is given by the following formulas:
\begin{eqnarray}
\label{pin}
&&\mathcal{H}v_k=(n-2k)v_k,\nonumber\\
&&\mathcal{F}v_k=v_{k+1}, \\
&&\mathcal{E}v_k=[k]_{\h}[n-k+1]_{\h}v_{k-1}.\nonumber
\end{eqnarray}
We also are be interested in the finite-dimensional modules $V_n\equiv V_n^+/V^+_{n-2}$. We will refer to the corresponding modules of $\fsh$ as $\pi^{+}_{n}(x)$ and $\pi_{n}(x)$ correspondingly. 
We also denote two-dimensional representation ${\rm V}(x):=\pi_1(x^{-1})$, which as a vector space is just a ${\rm span}(v_0, v_1)$.

We will be interested in the XXZ model which has the Hilbert space 
\begin{equation}
\mathscr{H}={\rm V}(a_1)\otimes \cdots \otimes {\rm V}(a_n)
\end{equation}

The normalized universal R-matrix as an element of the completed tensor product 
$$\mathcal{U}_{\hbar}(\hat{\mathfrak{b}}_{+})\hat{\otimes} \mathcal{U}_{\hbar}(\hat{\mathfrak{b}}_{-})$$ with $\mathcal{U}_{\hbar}(\hat{\mathfrak{b}}_{+})$ being represented via evaluation homomorphism and $\mathcal{U}_{\hbar}(\hat{\mathfrak{b}}_{-})$ considered in 2-dimensional representation $\pi_1$, is given by the following expression, see e.g. \cite{MR3199543}:
\begin{eqnarray}
&&(\pi_n(x)\otimes\pi_1(1))R=\phi_n(x)\mathcal{L}_n(x)\nonumber\\
&&\mathcal{L}_n(x)=
\left( {\begin{array}{cc}
{\hbar}^{\frac{\mathcal{H}}{4}}-{\h}^{-\frac{1}{2}} x^{-1}{\hbar}^{-\frac{\mathcal{H}}{4}} & (\hbar^{\frac{1}{2}}-{\hbar}^{-\frac{1}{2}})\mathcal{F}{\hbar}^{-\frac{\mathcal{H}}{4}}\\
x^{-1}(\hbar^{\frac{1}{2}}-\hbar^{-\frac{1}{2}})\mathcal{E}{\hbar}^{\frac{\mathcal{H}}{4}}&{\hbar}^{-\frac{\mathcal{H}}{4}}-{\hbar}^{-\frac{1}{2}} x{\hbar}^{\frac{\mathcal{H}}{4}}
\end{array}}
\right),
\end{eqnarray}
where $\phi_n(x)=\exp(\sum^{\infty}_{l=1}\frac{{\h}^{\frac{n+1}{2}}+{\h}^{-\frac{n+1}{2}}}{1+{\h}^{l}}\frac{x^{-l}}{l}$). The operator $\mathcal{L}(x)$ is known as the {\it $\mathcal{L}$-operator}  for the XXZ spin chain and the reason is as follows.

Now we  consider the tensor product 
$$\pi_1(u/a_1)\otimes \dots\otimes \pi_1(u/a_n)\otimes \pi_1(1),$$ where each of the 2-dimensional evaluation modules is refered as $site$ of the lattice and the last site is considered as auxilliary.

Thus we can define normalized quantum monodromy matrix as follows:
\begin{equation}\label{mon}
\mathcal{M}^Z(u)\equiv\mathcal{L}_1(u/a_1)\dots\mathcal{L}_n(u/a_n)\left( {\begin{array}{cc}
\zeta& 0\\
0&\zeta^{-1}
\end{array}}
\right),
\end{equation}
where the i-th $\mathcal{L}$-operator acts in the tensor product of the $i$-th cite and the auxilliary module. It is indeed has the form of discretized monodromy matrix for the linear differential problem produced by some  differential $\mathcal{L}$-operator. Note, that up to to the scalar multiple, it coincides with the quantum monodromy matrix from previous section, and thus the following integrability relation holds:
\begin{eqnarray}\label{int}
\Big[tr_{\pi_1(u_1)}\mathcal{M}^Z(u_1), tr_{\pi_1(u_2)}\mathcal{M}^Z(u_2)\Big]=0
\end{eqnarray} 
We can represent $\mathcal{M}^Z(u)$ as a matrix acting in $\pi_1(1)$ as follows: 
\begin{eqnarray}
\mathcal{M}^Z(u)=\left( {\begin{array}{cc}
A(u)& B(u)\\
C(u)&D(u)
\end{array}}
\right),
\end{eqnarray}
so that $tr_{\pi_1(1)}(\mathcal{M}^Z(u))=A(u)+D(u)\in End(\mathscr{H})[u,u^{-1}, \zeta, \zeta^{-1}]$.
Then the Yang-Baxter equation produces commutation relations between $A, B,C,D$ for different values of parameter $u$. Let us denote the highest weight vector in the product of our $N$ cites as
\begin{eqnarray}
\Omega_{+}\equiv \nu_0\otimes...\otimes \nu_0
\end{eqnarray}
Then $\{C(u)\}$ operators annihilate $\Omega_+$ and we could use $\{B(u)\}$-operators as creation operators and we obtain the following statement for  (see, e.g., \cite{Faddeev:aa}, \cite{Korepin_1993}, \cite{Reshetikhin:2010si}).
\begin{Thm}\label{theig}
Vectors $\{B(v_1)\dots B(v_k)\Omega_{+}\}$ are the eigenvectors of $tr_{\pi_1(1)}(\mathcal{M}^Z(u))$ with eigenvalues
\begin{eqnarray}
\label{transeig}
&&\Lambda(u|\{v_i\},\zeta)=\nonumber\\
&&\alpha(u)\prod_{i=1}^k\frac{v_i\hbar^{1/2}-u\hbar^{-1/2}}{v_i-u}+\delta(u)\prod_{i=1}^k\frac{v_i\hbar^{-1/2}-u\hbar^{1/2}}{v_i-u}
\end{eqnarray}
so that
\begin{eqnarray}
\alpha(u)=\zeta\hbar^{\frac{n}{4}}\prod_{i=1}^n\Big(1-\frac{a_i}{u\hbar}\Big ), \quad \delta(u)=\zeta^{-1}\hbar^{-\frac{n}{4}}\Big (1-\frac{a_i}{u}\Big)
\end{eqnarray}
are the eigenvalues of $A(u)$ and $D(u)$ on $\Omega_+$ and parameters $v_i$ are subject to Bethe equations:
\begin{eqnarray}\label{beq}
\prod_{j=1}^n\frac{\hbar^{1/2}v_i-\hbar^{-1/2}a_j}{v_i-a_j} =-
\zeta^{-2}\prod_{j=1, i\neq j}^k \frac{v_i\hbar^{1/2}-v_j\hbar^{-1/2}}{v_i\hbar^{-1/2}-v_j\hbar^{1/2}}.
\end{eqnarray}
\noindent ii) Vectors $\{B(v_1)\dots B(v_k)\Omega_{+}\}$ indexed by the solutions of Bethe equations, span the weight subspace $\mathscr{H}_k\subset\mathscr{H}$ corresponding to the eigenvectors of $\mathcal{H}$ with eigenvalue $n-2k$. 
\end{Thm}

An interesting question is the {\it completeness} of Bethe ansatz, namely whether this construction produces all eigenvectors (see, e.g.  Section 4.3 of \cite{Reshetikhin:2010si}.)

Denoting $x=u^{-1}$ and introducing functions $g_i(x)=(1-a_ix)$ and 
$\mathcal{Q}(x)$ is the operator with eigenvalues $\prod^k_{i=1}(1-v_ix)$, we have the following formula for the transfer matrix in terms of $Q$-operator:
$$
T_1(x):=tr_{\pi_1(1)}(\mathcal{M}^Z(u))=\zeta\hbar^{\frac{h_{\alpha_1}}{4}}\Bigg[ \prod^n_{i=1}g_i(x\hbar^{-1})\Bigg]\frac{\mathcal{Q}(\hbar x)}{\mathcal{Q}(x)} +\zeta^{-1}
\hbar^{-\frac{h_{\alpha_1}}{4}}\Bigg[ \prod^n_{i=1}g_i(x)\Bigg]\frac{\mathcal{Q}(\hbar^{-1} x)}{\mathcal{Q}(x)}
$$ 
In the next subsection we will discuss the representation theoretic construction of the operator $\mathcal{Q}(x)$, which is known as the Baxter $Q$-operator. 
 
\subsection{Baxter Q-operator and prefundamental representations} 
 
Let us return back to the fact that universal $R$-matrix is an element of 
$\mathcal{U}_{\hbar}(\hat{\mathfrak{b}}_{+})\hat{\otimes }\mathcal{U}_{\hbar}(\hat{\mathfrak{b}}_-)$.  That indicates that in principle one can extend the category of our finite-dimensional modules and consider representations of the Borel subalgebra $\mathcal{U}_{\hbar}(\hat{\mathfrak{b}}_-)$, thus  using the $R$-matrix as an intertwiner between those representations.  There is a family of such modules which generate the entire extended Grothendieck ring, known as prefundamental representations. They were introduced in the case 
of $\fsh$ by Bazhanov, Lukyanov and Zamolodchikov \cite{Bazhanov:1998dq}. The conceptual understanding of the prefundamental representations was given by Frenkel and Hernandez in \cite{HJ}, \cite{Frenkel:2013uda}, \cite{Frenkel:2016}. The construction of such modules in the case of $Y_{\hbar}(\mathfrak{g})$ is not yet known, although some attempts were given in \cite{BFLMS}. 

We will now describe this construction in the case of $\fsh$ following \cite{Bazhanov:1998dq}.

First, we introduce the deformed oscillator algebra $\mathcal{B}$ with the generators $H, \mathcal{E}_{\pm}$:
\begin{equation}
{\h}^{\frac{1}{2}} \mathcal{E}_+\mathcal{E}_- - {\h}^{-\frac{1}{2}}\mathcal{E}_-\mathcal{E}_+=\frac{1}{{\h}^{\frac{1}{2}}-{\h}^{-\frac{1}{2}}},\quad  [H, \mathcal{E}_{\pm}]=\pm 2\mathcal{E}_{\pm}.
\end{equation}
The algebra  $\mathcal{B}$ has the following Fock space representations:
\begin{equation}
\mathbf{F}_{\pm}=
\{{\rm span}\{\mathcal{E}^k_{\mp}|0\rangle_{\pm}\};~ \mathcal{E}_{\pm}|0\rangle_{\pm}=0,~ H|0\rangle_{\pm}=0\}
\end{equation}

The representations of the lower Borel subalgebra $\mathcal{U}_{\hbar}(\hat{\mathfrak{b}}_-)$ of $\mathcal{U}_{{\h}}(\widehat{\mathfrak{sl}}_2)$, which we are interested in, 
can be described by the following homomorphisms $\rho_{\pm}(x): \mathcal{U}_{\hbar}(\hat{\mathfrak{b}}_-)\to\mathcal{B}$:
\begin{equation}
 \rho_{\pm}(x): \left\{
    \begin{array}{rl}
      &h_{\alpha_1}\to \pm H,\quad  h_{\alpha_0}  \to \mp H\\
      &e_{-\alpha_1}  \to \mathcal{E}_{\mp},\quad e_{-\alpha_0} \to x\mathcal{E}_{\pm},
    \end{array} \right.
\end{equation}
The above homomorphisms thus automatically produce representations of $\mathcal{U}_{\hbar}(\hat{\mathfrak{b}}_-)$ on $\mathbf{F}_{\pm}$. In the following we will refer to the corresponding representations as $\rho_{\pm}(x)$ as well.
 
We are interested in the decomposition of the tensor product
\begin{equation}
\rho_{-}(x{\h}^{-\frac{n+1}{2}})\otimes \rho_{+}(x{\h}^{\frac{n+1}{2}})\label{prodqosc}
\end{equation}
for $n\in \mathbb{Z}$. 
 
To describe the answer, let us also introduce 1-dimensional representations of $\mathcal{U}_{\hbar}(\hat{\mathfrak{b}}_-)$ with eigenvalue of $h_{\alpha_1}$ equal to $s$ as $\omega_{s}$. It is clear that $\omega_s\otimes \omega_{s'}=\omega_{s+s'}$.
Then we have the following Proposition \cite{blz}.
\begin{Prop}
Decomposition of the product \ref{prodqosc} in the Grothendieck ring of $\mathcal{U}_{\hbar}(\hat{\mathfrak{b}}_-)$ is given by the following expresion:
\begin{eqnarray}
\rho_{-}(x\h^{-\frac{n+1}{2}})\cdot \rho_{+}(x{\h}^{\frac{n+1}{2}})=\omega_{-n}(1-\omega_{-2})\inv\pi_n^{+}(x),\label{grot}
\end{eqnarray}
where $(1-\omega_{-2})^{-1}$ is understood as the geometric series expansion.
\end{Prop}
 
Now we describe the weighted traces of R-matrices in those representations. First of all, from now on, we assume that whenever we write the trace of the universal R-matrix, which belongs to $\mathcal{U}_{\hbar}(\hat{\mathfrak{b}}_{+})\hat{\otimes} \mathcal{U}_{\hbar}(\hat{\mathfrak{b}}_{-})$ over certain representation of its $\mathcal{U}_{\hbar}(\hat{\mathfrak{b}}_{-})$ part, its $\mathcal{U}_{\hbar}(\hat{\mathfrak{b}}_{+})$-part is considered in some finite-dimensional representation $V$, so that the following trace
\begin{eqnarray}
\tilde{Q}_{\pm}(x)=tr_{\rho_{\pm}(x)}\Big[ (I\otimes\rho_{\pm}(x))R(I\otimes \zeta^{h_{\alpha_1}})\Big]
\end{eqnarray}
is well-defined as an element of ${\rm End}(V)[[x, \zeta^{-1}]]$,  and we will treat it as such in the following (see \cite{Frenkel:2013uda}).  One can also define it analytically for $Re\{\ln(\zeta)\}>0$ and then analytically continue the result.  
 
We would like to normalize this operator, so that at $x=0$ it is equal 1.  Note, that  from the structure of the universal $R$-matrix we have $\lim_{x\rightarrow 0} (I\otimes\rho_{\pm}(x))R={\h}^{\frac{h_{\alpha_1}\otimes h_{\alpha_1}}{4}}$. 
Thus introducing 
\begin{eqnarray}
&&\mathcal{Z}_{\pm}(h_{\alpha_1})=\nonumber\\
&&tr_{\rho_{\pm}(x)}\Big[ (I\otimes\rho_{\pm}(x)){\h}^{\frac{h_{\alpha}\otimes h_{\alpha}}{4}}(I\otimes \zeta^{h_{\alpha}})\Big]=\sum^{\infty}_{s=0}{\h}^{-s\frac{h_{\alpha_1}}{2}}\zeta^{-2s}=\frac{1}{1-\h^{-\frac{h_{\alpha_1}}{2}}\zeta^{-2}},
\end{eqnarray}
again, understood as an element of ${\rm End}(V)[[ \zeta^{-1}]]$, we obtain that $Q_{\pm}(x)\equiv\mathcal{Z}_{\pm}^{-1}(h_{\alpha})\tilde{Q}_{\pm}(x)$ have the property $Q_{\pm}(0)=1$. 

Then we have the following proposition \cite{blz}, which is a consequence of the relation (\ref{grot}).

\begin{Prop} The product of two $Q$-operators, acting in certain finite-dimensional representation $V$, gives a trace of the R-matrix in $\pi_n^+(x)$ representation:
\begin{equation}
\h^{\frac{n{h_{\alpha_1}}}{4}}\zeta^nQ_{+}(\h^{\frac{n+1}{2}}x)Q_{-}(\h^{-\frac{n+1}{2}}x)=
tr_{\pi^{+}_n(x)}\Big[ R(I\otimes {\zeta}^{h_{\alpha_1}})\Big]\Big (1-\hbar^{-\frac{h_{\alpha_1}}{2}}\zeta^{-2}\Big ),
\end{equation}
where the trace in the RHS is understood as an element in ${\rm End}(V)\zeta^{n}[\zeta^{-1},x]]$. 
 \end{Prop}
 
To prove that, it is enough to see that
\begin{equation}
tr_{\omega_s}\Big[ R(I\otimes {\zeta}^{h_{\alpha_1}})\Big]={\h}^{\frac{sh_{\alpha_1}}{4}}\zeta^s.
\end{equation}
Then formula (\ref{grot}) corrected by normalization coefficients gives the result of the Theorem.

Let us introduce the following notation for the following traces, again, understood as power series from ${\rm End}(V)\zeta^{n}[[\zeta^{-1},x]]$ for certain finite-dimensional representation $V$ of $\mathcal{U}_{\hbar}(\hat{\mathfrak{b}}_+)$:
\begin{eqnarray}
&&T_n^{+}(x)\equiv tr_{\pi^{+}_n(x)}\Big[ R(I\otimes {\zeta}^{h_{\alpha_1}})\Big],
\nonumber\\
&&T_n(x)\equiv tr_{\pi_n(x)}\Big[ R(I\otimes {\zeta}^{h_{\alpha_1}})\Big]=T_n^{+}(x)-T_{-n-2}^{+}(x).
\end{eqnarray}
There are two direct consequences of the Theorem we proved. One is known as the {\it QQ-relation} between $Q_{\pm}(x)$-operators, illustrating that they are not independent. The second, known as the {\it Baxter TQ-relation}, expresses the dependence of $T_1$ on  $Q_{\pm}$- operators (see e.g. \cite{blz}). 
\begin{Thm}\label{wr}
Opertaors $Q_{\pm}(u)$, $T_1(u)$, acting in a finite-dimensional representation $V$ obey the following relations:
\begin{eqnarray}
&&\zeta{\h}^{\frac{h_{\alpha_1}}{4}}Q_{+}(\h^{\frac{1}{2}} x)Q_{-}(\h^{-\frac{1}{2}}x)-\zeta^{-1}{\h}^{-\frac{h_{\alpha_1}}{4}}Q_{+}(\h^{-\frac{1}{2}}x)Q_{-}(\h^{\frac{1}{2}} x)=\zeta{\h}^{\frac{h_{\alpha_1}}{4}}-\zeta^{-1}{\h}^{-\frac{h_{\alpha_1}}{4}},\\
&&\nonumber\\
&&T_1(x)Q_{\pm}(x)={\h}^{\pm\frac{h_{\alpha_1}}{4}}\zeta^{\pm 1}Q_{\pm}({\h}x)+{\h}^{\mp\frac{h_{\alpha_1}}{4}}\zeta^{\mp 1}Q_{\pm}(\h^{-1}x).
\end{eqnarray}

\end{Thm}

Thus we can identify the nomalized operator $Q_+(x)$ with the operator $\mathcal{Q}_+(x)$, which we introduced in the previous subsection in the ad hoc way. Namely, introducing the $q$-Pochhammer symbol:
\begin{eqnarray}\label{poch}
(a;q)_{\infty}=\prod^{\infty}_{k=0}(1-aq^k),
\end{eqnarray}
the following is true.

\begin{Prop}\label{norm}
The eigenvalues of the operator
\begin{eqnarray}
\mathcal{Q}_+(x)=\prod^n_{i=1}\frac{(a_ix;{\h} ^2)_{\infty}}{ (a_ix{\h};{\h}^2)_{\infty}}Q_+(x),
\end{eqnarray}
 are equal to 
$\prod_{i=1}^s\Big(1-{x}{v_i}\Big)$ and they are attained on the vectors
$\{P(v_1,\dots v_k)\}$, provided the Bethe ansatz equations (\ref{beq})
are satisfied.
\end{Prop}

One can reproduce
$Q_-(x)$ from algebraic Bethe ansatz as well. In order to do that, one has to use lowest weight vector
$\Omega_-=\nu_1\otimes\dots\otimes \nu_n$. It is annihilated by $B(u)$-operators and the space of states is spanned by the operators $C(w_1)\dots C(w_k)$ acting on $\Omega_-$.
The following Theorem is an analogue of Theorems \ref{theig} and \ref{norm}.
\begin{Prop}
Vectors $\{P_{-}(w_1,\dots w_k)=C(w_1)\dots C(w_k)\Omega_-\}$ are the eigenvectors for $tr(\mathcal{T}(u))$ with eigenvectors
 $\Lambda(u|\{w_i\},\zeta^{-1})$ (see \ref{transeig}) , so that $\{w_i\}^k_{i=1}$ are solutions of  Bethe equations (\ref{beq}) with $\zeta\rightarrow \zeta^{-1}$.
The operator
\begin{eqnarray}
\mathcal{Q}_-(x)=\prod^n_{i=1}\frac{(a_ix;{\h} ^2)_{\infty}}{ (a_ix{\h};{\h}^2)_{\infty}}Q_-(x)
\end{eqnarray}
has eigenvalues
$\prod_{i=1}^n(1-x w_i)$ on the vectors
$P_-(w_1,\dots w_k)$.

\noindent ii) Vectors $P_-(s_1,\dots s_k)$ span the weight subspace $\tilde{\mathscr{H}}_k\subset\mathscr{H}$ corresponding to the eigenvectors of ${\mathcal{H}}$ with eigenvalue $2k-n$. 
\end{Prop}

The following Proposition gives a normalized version of the quantum Wronskian relation (see Theorem \ref{wr}).\\

\begin{Prop}
The Wronskian relation between $\mathcal{Q}_{\pm}$-operators reads as follows:
\begin{eqnarray}
&&\zeta\h^{\frac{h_{\alpha}}{4}}\mathcal{Q}_{+}({\h}^{\frac{1}{2}} x)
\mathcal{Q}_{-}({\h}^{-\frac{1}{2}}x)-
\zeta^{-1}{\h}^{-\frac{h_{\alpha}}{4}}
\mathcal{Q}_{+}({\h}^{-\frac{1}{2}}x)
\mathcal{Q}_{-}({\h}^{\frac{1}{2}} x)=\nonumber\\
&&(\zeta{\h}^{\frac{h_{\alpha}}{4}}-{\h}^{-\frac{h_{\alpha}}{4}}\zeta^{-1})\prod^n_{i=1}(1-a_i{\h}^{-\frac{1}{2}}x).
\end{eqnarray}
\end{Prop}
 Reducing it to eigenvalues, one can view it as a functional relation, known as a $QQ$-system, which we will discuss in the next subsection.

\subsection{Modern approaches to Bethe ansatz}
\subsubsection{Quantum Knizhnik-Zamolodchikov equations}
The intertwining operators for the quantum affine algebra $U_{\hbar}(\hat{\mathfrak{g}})$ and thus the matrix elements of their products, known as {\it conformal blocks} satisfy certain difference equations known as quantum Kniznik-Zamolodchikov (qKZ) equations (also known as Frenkel-Reshetikhin equations) \cite{FR1998}.  Explicitly, qKZ equations can be written as the following difference equations: 
\begin{eqnarray}\label{qkz}
\Psi(a_{i_1}, \dots, q a_{i_k}, \dots , a_{i_n}, \{z_i\})=H^{(q)}_{i_k}\Psi(a_{i_1}, \dots,, a_{i_n}, \{z_i\}).
\end{eqnarray}
Here the solution $\Psi$ takes values in the physical space $\mathcal{H}$ and operators $H^{(q)}_i$ are expressed in terms of products of R-matrices. 
Moreover, it is known that $H^{(1)}_i=\lim_{q\rightarrow 1}H^{(1)}_i$ coincide with certain transfer matrices of the corresponding XXZ model (see e.g. \cite{FHRnew}). 

There is also a commuting system of equations in  $\{z_i\}$-variables for $\Psi$, known as {\it dynamical equations} 
 see e.g. \cite{Tarasov_2002,Tarasov_2005}.
 
 The solution to the qKZ equation is given by an integral expression, so that the integrand has the following asymptotic behavior in the limit $q\rightarrow 1$ (or $\eta=\log(q)\rightarrow 0$): 
\begin{eqnarray} 
 e^{\frac{Y(\{a_i\}, \{z_i\}, \{x_i\})}{\eta}}\Big[\phi_0(\{a_i\}, \{z_i\}, \{x_i\})+O(\eta)\Big],
 \end{eqnarray}
where $\{x_i\}$ are the variables of integration.  
In the limit $q\rightarrow 1$ the stationary phase approximation gives $\Psi=e^{\frac{S}{\eta}}(\Psi_0+O(\eta))$, where $S=Y|_{\sigma_i},$  
where $\sigma_i$ are the solutions of the equations $\partial_{x_i}Y=0$  which coincide with the Bethe equations, and $\Psi_0$ is the eigenvector for transfer matrices $H^{(1)}_i$:
$$H^{(1)}_i\Psi_0=e^{p_i}\Psi_0,~{\rm where} ~p_i={a_i}\partial_{a_i}S.$$ 

 In Section 3 we will discuss the geometric realization \cite{Okounkov:2015aa}, \cite{Okounkov:2016sya} of the qKZ equations and dynamical equations.

\subsubsection{$QQ$-systems and Baxter operators}

When we earlier discussed the transfer matrices $T_{W(u)}$ we considered $W(u)$ to be a finite-dimensional module of $\mathcal{U}_{\hbar}(\hat{\mathfrak{g}})$. We also noticed that the {\it universal R-matrix}, which produces particular braiding operators $R_{V_i(a_i), V_j(a_j)}$ belongs to the completion of the tensor product $\mathcal{U}_{\hbar}(\hat{\mathfrak{b}}_{+})\otimes \mathcal{U}_{\hbar}(\hat{\mathfrak{b}}_{-})$, where $\mathcal{U}_{\hbar}(\hat{\mathfrak{b}}_{\pm})$ are the Borel subalgebras of $\mathcal{U}_{\hbar}(\hat{\mathfrak{g}})$. Therefore, there is no obstruction in taking auxiliary representations $W(u)$ to be representations of $\mathcal{U}_{\hbar}(\hat{\mathfrak{b}}_{+})$.  We discussed in the previous section using the example of $\mathcal{U}_{\hbar}(\widehat{\mathfrak{sl}}_2)$ the prefundamental representations of $\mathcal{U}_{\hbar}(\hat{\mathfrak{b}}_{+})$ which are infinite-dimensional. If one extends the braided tensor category of finite-dimensional modules by such representations, the Grothendieck ring of the resulting category is generated by those modules. The corresponding transfer matrices, namely Baxter $Q$-operators, turn out to be well-defined, and upon rescaling by a certain coefficient function, the eigenvalues of the transfer matrices are polynomials of the spectral parameter, generating elementary symmetric functions of the solutions of Bethe equations. Although originally introduced by Baxter via their eigenvalues, their representation-theoretic meaning was realized much later in the papers of Frenkel and Hernandez \cite{Frenkel:2013uda}, \cite{Frenkel:2016}, following earlier ideas of Bazhanov, Lukyanov, and Zamolodchikov \cite{Bazhanov:1998dq} and later, by Hernandez and Jimbo \cite{HJ} (see most recent ideas in \cite{frenkel2023extended}). We note that similar construction and the analog of the $QQ$-system should also exist for Yangians, with some progress being made in \cite{BFLMS}, \cite{ekhammar2021extended}.

For a simple Lie algebra $\mathfrak{g}$, there are two series of prefundamental representations $\{V^i_{\pm}(u)\}^r_{i=1}$ and the associated Baxter operators $\{Q^i_{\pm}(u)\}^r_{i=1}$.  They obey the following relations \cite{Frenkel:2016}:
\begin{eqnarray}\label{QQsec2}
{{\widetilde\xi_i}}Q^i_{-}({u})Q^i_{+}({\hbar}  {u})-{\xi_i}Q^i_{-}({\hbar}  {u})Q^i_{+}({u}) &=&\Lambda_i({u})\prod_{j\neq i}\Bigg[\prod^{-a_{ij}}_{k=1} Q^j_{+}({\hbar} ^{b_{ij}^k}{u})\Bigg]\nonumber\\ 
i&=&1,\dots,r, \quad b_{ij}^k\in \mathbb{Z}
\end{eqnarray}
Here polynomials $\Lambda_i({u})$ are related to {\it Drinfeld polynomials} \cite{Chari:1994pz}, characterizing the representation $\mathcal{H}$ of $U_{{\hbar}}(\hat{\mathfrak{g}})$ and 
${\xi_i}$, ${{\widetilde\xi_i}}$ are some monomials of $\{z_i\}$.  
This system of equations, known as the {\it $QQ$-system}, considered as equations on  $\{Q^i_{+}(u)\}_{i=1,\dots, r}$ and subject to some nondegeneracy conditions, is equivalent to the system of Bethe ansatz equations, where the {\it Bethe roots} are the roots of $Q_+$-polynomials. This was noted and studied in various instances, see, e.g.,   \cite{MVflag}, \cite{2006math......4048M}, \cite{Mukhin_2005}, \cite{2012arXiv1210.2315L}, \cite{MR2189873}. 

We will discuss the various contexts in which $QQ$-systems emerge in Section 4 of these notes. From a geometric perspective, we show that equations (\ref{QQsec2}) characterize a difference analog of $G$-connections ($G$ here is a simple simply connected Lie group) on a projective line, called  $(G, \hbar)$-opers. 
This is a deformation of a celebrated example of geometric Langlands correspondence \cite{Frenkel_LanglandsLoop}, where instead of the XXZ model, one has its limit, 
called Gaudin model associated to Langlands dual Lie algebra $^L\mathfrak{g}$, and the actual $G$-oper connections on the projective line, studied by B. Feigin, E. Frenkel and collaborators over the years \cite{Feigin:1994in}, \cite{Frenkel:2003qx}, \cite{FFTL:2010}, \cite{Rybnikov:2010}.  
In this non-deformed case, the differential limit of equations (\ref{QQsec2}) plays a role as well \cite{mukhvarmiura}, \cite{Brinson:2021ww}.

\section{Quantum K-theory of Nakajima quiver varieties}
\subsection{Nakajima quiver varieties}\label{qvar}
{\it Nakajima quiver varieties} are important examples of the so-called {\it symplectic resolutions} (see ,e.g.,\cite{Ginzburg:},\cite{Kaledinreview}, \cite{Kamnitzer:2022aa}). Their localized equivariant cohomology and K-theory have the structure  of representations of the corresponding Yangian $Y_{\hbar}(\mathfrak{g}^{Q})$/quantum algebra $U_{\hbar}(\mathfrak{g})$, which corresponds to the quiver, which in the simply laced $ADE$ case is in one-to-one correspondence with the Dynkin diagram \cite{nakajima1998, Nakajima:1999hilb, Nakajima:2001qg},  \cite{olivier1998quantum}, \cite{Vasserot:wo}, \cite{varagnolo2000quiver}.  For simplicity throughout these lectures we will consider only simply-laced case, i.e. the number of quiver vertices $I=\{1, \dots, {\rm rank}(\mathfrak{g})\}$.

The quiver varities naturally emerge as the Higgs vacua of 2d and 3d supersymmetric gauge theories. Here we give necessary definitions with a dictionary for the physical terminology.

A {\it quiver} consists of vertices $I$ and oriented edges connecting them. 
A framed quiver is a quiver, where the set of vertices is doubled, and each of the vertices in the added set has an edge going from it to the vertex, whose copy it is. 

To a given framed quiver one assigns the following data:

\begin{itemize}

\item  A set of the following vector spaces $V_i,W_i$, where $V_i$  correspond to original vertices, and $W_i$, known to physicists as {\it flavor space} correspond to their copies, together with a set of morphisms between these vertices, corresponding to the edges of the quiver. One assigns the dimension labels to the corresponding vertices: $\mathbf{v}_i=\text{dim}\ V_i$, $\mathbf{w}_i=\text{dim}\ W_i$. 
We also denote by $Q_{ij}$ the incidence matrix of the quiver, i.e. the number of edges between vertices $i$ and $j$.

\item For a given framed quiver we construct the affine space $M:=\text{Rep}(\mathbf{v},\mathbf{w})$, where $\mathbf{v}$ and $\mathbf{w}$ are dimension vectors of the quiver in the following way. It is a direct sum of: \\

\begin{enumerate} 
\item  $\sum_{i\in I}\text{Hom}({W}_i,{V}_i)$;\\
\item  $\sum_{i,j\in I}{ Q}_{ij}\otimes \text{Hom}({V}_i,{V}_j).$\\
\end{enumerate}
\end{itemize}
known to physicists as the {\it matter} data of {\it fundamental fields} labeled by flavor framing space $W_i$ and {\it  bifundamental fields} correspondingly.

In other  words, for a given framed quiver,  
$$M=\sum_{i\in I}\text{Hom}({W}_i,{V}_i)\oplus\sum_{i,j\in I}{ Q}_{ij}\otimes \text{Hom}({V}_i,{V}_j).$$ 
This space  has the natural group action: 
$$G=GL(V_1)\times\dots\times GL(V_{{\rm rank}(\mathfrak{g})}).$$ 
To have enough supersymmetries in the gauge theory to describe vacua one should  consider a doubling of such space, i.e. $T^*M$. In this case the action of $G$ is hamiltonian with the moment map $\mu$. 
Let $L_{\mathbf{v},\mathbf{w}}=\mu^{-1}(0)$ be the zero locus of the moment map or the {\it vacua} condition.

The Nakajima variety $X$ corresponding to the quiver is an algebraic symplectic reduction
$$
X=N_{\mathbf{v},\mathbf{w}}=L_{\mathbf{v},\mathbf{w}}/\!\! /_{\theta}G=L^{ss}_{\mathbf{v},\mathbf{w}}/G\,.
$$
Here superscript $ss$ stands for semi-stable locus which depends on a choice of stability parameter $\theta\in {\mathbb{Z}}^I$ (see, e.g., \cite{Ginzburg:} for a detailed definition). The group
$$\prod_{i.j} GL(Q_{ij})\times
\prod_i GL({W}_i)\times \mathbb{C}^{\times}_\hbar$$
acts as automorphisms of $X$, coming form its action on $\text{Rep}(\mathbf{v},\mathbf{w})$. Here  $\mathbb{C}^{\times}_{\hbar}$ scales cotangent directions with weight $\hbar$ and therefore symplectic form with weight $\hbar^{-1}$. Let us denote by $\mathsf T=\bA\times \mathbb{C}^{\times}_\hbar$ as the maximal torus of this group.

As a Nakajima quiver variety $\N_{\bf v, \bf w}$ is a symplectic resolution,
i.e. it comes with natural projective morphism to affine variety, Theorem 5.2.2 in \cite{Ginzburg:}:
\begin{equation} \label{affineN}
\N_{{\bf v},{\bf w}}\rightarrow \N^{0}_{{\bf v},{\bf w}}:={\rm{Spec}}\Big(\C[\mu^{-1}(0)]^{G}\Big).
\end{equation}

\noindent 
{\bf Examples.}  Our main example is the simplest nontrivial situation with one vertex with a framing with assigned vector spaces $V$, $W$. In this case 
$\N(\mathbf{v},\mathbf{w})=\N_{k,n}={\rm T^*Gr}_{k,n}$, where we denoted $\dim(V)=k$, $\dim(W)=n$ In this case
$$
M={\rm Hom}(V,W), \quad \N_{k,n}={\rm Hom}(V, W)\oplus {\rm Hom}(W, V)/\!\! /_{\theta}GL(V)
$$
In this case momentum map $\mu:T^*M\to \mathfrak{gl}(V)^*$ is given by $\mu(A,B)=BA$. In this case the {stable points in $T^*M$ are given by $\{ (A,B): rank(A)=k \}$.

One can generalize this example to the cotangent the space of partial flags, i.e. we are considering the following quiver of type $A_n$ with single framing vertex:

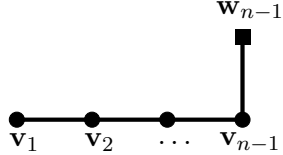
\vspace{0.1in}
\begin{center}
\begin{tikzpicture}
\draw [ultra thick] (0,0) -- (3,0);
\draw [ultra thick] (3,1) -- (3,0);
\draw [fill] (0,0) circle [radius=0.1];
\draw [fill] (1,0) circle [radius=0.1];
\draw [fill] (2,0) circle [radius=0.1];
\draw [fill] (3,0) circle [radius=0.1];
\node (1) at (0.1,-0.3) {$\mathbf{v}_{1}$};
\node (2) at (1.1,-0.3) {$\mathbf{v}_2$};
\node (3) at (2.1,-0.3) {$\ldots$};
\node (4) at (3.1,-0.3) {$\mathbf{v}_{n-1}$};
\fill [ultra thick] (3-0.1,1) rectangle (3.1,1.2);
\node (5) at (3.1,1.45) {$\mathbf{w}_{n-1}$};
\end{tikzpicture}
\end{center}
\vspace{0.1in}

The stability condition is chosen so that maps $V_{n-1}\to W_{n-1}$ and $V_i\to V_{i+1}$ are injective. For the variety to be non-empty the sequence $\mathbf{v}_{1},\ldots ,\mathbf{v}_{n-1}, \mathbf{w}_{n-1}$ must be non-decreasing.

For a Nakajima quiver variety $\N_{\bf{v}, \bf{w}}$ one can define a set of tautological bundles 
$$
\mathcal{V}_i=L_{{\bf v}, {\bf w}}^{ss}\times V_i/G,\quad  \mathcal{W}_i=L_{{\bf v}, {\bf w}}^{ss}\times W_i/G.
$$
From this construction it follows that all bundles $\mathcal{W}_i$ are topologically trivial.  It is known that tensorial polynomials of these bundles and their duals generate $K_{\mathsf{T}}(X)$: that is known as a Kirwan surjectivity conjecture proven in \cite{McGerty:2016kir}.  For more details regarding quiver varieties, one can consult with e.g. \cite{Ginzburg:},  introduction to \cite{2012arXiv1211.1287M} or Section 4 of \cite{Okounkov:2015aa}.

From now on we will mainly focus on the simplest example of $\N_{{\bf v},{\bf w}}$: the cotangent bundles to Grassmannians $\N_{k,n}={\rm T^* Gr}_{k,n}$.

\subsection{${\rm T^* Gr}_{k,n}$ and its equivariant K-theory}

Let us choose a basis in the $n$-dimensional vector space with coordinates $(x_1,\dots,x_n)$. 
The corresponding $n$-dimensional torus $\bA=(\C^{\times})^n$ acting in this vector space by scaling the coordinates in the chosen basis:
$$
(x_1,\dots,x_n)\rightarrow (x_1 a_1,\dots,x_n a_n).
$$
Thus we obtain the action of $\bA$ on $k$-subspaces, i.e., on $\textrm{Gr}_{k,n}$:  a $k$-subspace representing a point in  $\textrm{Gr}_{k,n}$ is fixed under this action if and only if it is spanned by $k$ basis vectors.  The fixed points $\textrm{Gr}_{k,n}^{\bA}$ are in one-to-one  correspondence with the $k$-subsets of the set $\{1,2,\dots,n\}$, so that each such subset corresponds to a choice of coordinate $k$-vectors from the coordinate $n$-vectors. As we discussed in previous section, one can further extend the torus action on the symplectic manifold $\N_{k,n}$ to $\mathsf T=\bA\times \C^{\times}_{\hbar}$, so that the extra one-dimensional torus $\C^{\times}_{\hbar}$ acts on $\N_{k,n}$ by scaling the cotangent directions with $\hbar$. This new action contracts fibers, i.e. contracts $\N_{k,n}$ to $\textrm{Gr}_{k,n}$, the fixed point set
$\N_{k,n}^{\mathsf T}$ is the same as $\textrm{Gr}_{k,n}^{\bA}$. 
In other words, the fixed set $\N_{k,n}^{\mathsf T}$ is a set of
$n!/k!/(n-k)!$ points labeled by $k$-subsets 
$\{i_1,\dots,i_k\}\subset \{1,\dots, n\}.$  

The equivariant $K$-theory $K_{\mathsf T}(\N_{k,n})$ is generated by tensorial polynomials in  tautological bundle $\mathcal{V}$: this is as a particular case of Kirwan surjectivity we discussed in the previous section. Localization theorem in equivariant K-theory implies that the classes of fixed points is a basis of the localized K-theory. In other words, one can think about
the localized K-theory $K_{\mathsf T}(\N_{k,n})_{loc}$ as a 
$n!/k!/(n-k)!$-dimensional vector space spanned by the classes of fixed points from $\N_{k,n}^{\mathsf T}$,  forgetting all other structures. One can  label them by $k$-subsets
$\fp=\{x_1,\cdots,x_{k}\}\subset \{a_1,\cdots, a_n\}$.

Let $K_{GL(V)\times \mathsf T}(pt) = \Z[s_1^{\pm 1},s^{\pm 1}_2,\cdots,s_k^{\pm 1},a_1^{\pm 1},\dots, a_{n}^{\pm},\hbar^{\pm}]^{\frak{S}_k}$ be the ring of symmetric Laurent polynomials in $k$ variables with coefficients in $K_{\mathsf T}(pt)$. For every such polynomial $\tau$, we denote by the same symbol $\tau \in K_{\mathsf T}(\N_{k,n})$ the corresponding Schur functor of $\mathcal{V}$. For example, 
	$$\tau(s_1,\cdots,s_k)=(s_1+\cdots + s_k)^3 - \sum_{1\leq i_1<i_2\leq k} s_{i_1}^{-1} s_{i_2}^{-1} $$
	corresponds to $\tau(V) = V^{\otimes 3} - \Lambda^{2} V^{*}$. 

Let us set the following notation for the disjoint union of $\N_{k,n}$ for all $k$:
$\N(n)=\coprod\limits_{k=0}^{n}\,\N_{k,n},$  
so that the fixed point set $\N(n)^{\mathsf T}$ consists of total $2^n$ points.

The equivariant $K$-theory $K_{\mathsf T}(\N(n))$ is a module over the ring of equivariant constants: $R=K_{\mathsf T}(\cdot)=\Z[a_1^{\pm },\cdots,a_n^{\pm 1},\hbar^{\pm 1}]$. The localized $K$-theory
\begin{eqnarray}
K_{\mathsf T}(\N(n))_{loc}=K_{\mathsf T}(\N(n)) \bigotimes\limits_{R} {\mathcal{A}} = \bigoplus_{k=0}^{n} K_{\mathsf T}(\N_{k,n}) \bigotimes\limits_{R} {\mathcal{A}}, ~{\rm where}~  {\mathcal{A}}=\Q(a_1,\cdots,a_n,\hbar),\nonumber
\end{eqnarray}
is an ${\mathcal{A}}$-vector space  of dimension $2^n$ spanned by the $K$-theory classes of fixed points $\mathcal{O}_{\fp}$. If $\tau$ is a $K$-theory class then the operation 
of tensor multiplication by $\tau$, namely $\gamma \to \tau \otimes \gamma$ is an $\mathcal{A}$-linear operator
acting in the vector space above.  These operators are diagonal in the basis of fixed points:
\begin{eqnarray}
\label{clev}
\tau \otimes \mathcal{O}_{\fp} =\tau(a_{i_1},\cdots, a_{i_k}) \, \mathcal{O}_{\fp} \ \ \ \textrm{for}  \ \ \ \fp=\{i_{1},\cdots, i_{k}\}\subset \{1,\cdots, n\}.
\end{eqnarray}

In other words, the eigenvalues of the operators of multiplication by tautological bundles on ${K_{\mathsf T}(\N_{k,n})}_{loc}$ are given by the values of the corresponding Laurent polynomials $\tau(s_1,\cdots,s_k)$ evaluated at the solutions of the following equations:
\begin{eqnarray}
\label{clbeth}
\prod\limits_{j=1}^{n}(s_i-a_j) =0,  \ \ \ i=1, \dots, k
\end{eqnarray}
with $s_i\neq s_j$.  The solutions of (\ref{clbeth}) with $s_i\neq s_j$  are in one-to-one correspondence with the
$k$-subsets $\{a_{i_1},\cdots,a_{i_k}\}\subset \{a_{1},\cdots, a_{n}\}$ and, therefore, with the set of the fixed points $\N_{k,n}^{\mathsf T}$. 
We are looking to deform this statement to the case of the  quantum $K$-theory.
We note that the system of equations (\ref{clbeth}) turns out to be the l limit $\zeta\rightarrow 0$ of the Bethe ansatz equations~(\ref{beq}).

\subsection{Geometric quantum group action}

We start with the variety $\N(1)$: it is a union of two components (points) corresponding to the Grassmannians of zero- and one-dimensional hyperplanes in one-dimensional vector space.  Thus $K_{\mathsf T}(\N(1))$ is a two-dimensional vector space over the field $\mathcal{A}$. 
For  the equivariant parameter $a$ of $\bA\cong\C^{\times}$ we denote 
$
\mathbb{F}(a):=K_{\mathsf T}(\N(1))
$
the corresponding two-dimensional vector space. 
  
For general Nakajima variety $X=\N(\mathbf{v}, \mathbf{w})$, if $\sigma:\mathbb{C}^{\times}\to \bA$ is such that $\mathbf{w}=a'\mathbf{w}'+a''\mathbf{w}''$, where $a',a''\in \mathbb{C}^{\times}$, we obtain (see Section 2.4 of \cite{2012arXiv1211.1287M}) the property which is known as a tensor product structure : 

\begin{equation}
 X^{\sigma}=\bigsqcup_{\mathbf{v}'+\mathbf{v}''=\mathbf{v}}\N(\mathbf{v}',\mathbf{w}')\times \N(\mathbf{v}'',\mathbf{w}'')\,. \notag
 \end{equation}
In our particular case, we can iterate this procedure to obtain:
$$
\N(n)^{\bA}=\N(1)^n=\N(1)\times \cdots \times \N(1),
$$
which leads to the following isomorphism of $\mathcal{A}$-vector spaces:
\begin{eqnarray} \label{tp}
K_{\mathsf T}(\N(n))_{loc}=K_{\mathsf T}(\N(n)^{\bA})\bigotimes\limits_{R} \mathcal{A}=\mathbb{F}(a_1)\otimes \cdots \otimes \mathbb{F}(a_n).
\end{eqnarray} 

In section 9 of \cite{Okounkov:2015aa} the following maps in localized K-theory, known as stable envelopes: 
$$\text{Stab}^s_{\pm}(a): K_{\mathsf T}(\N_{\bf v, w}^{\sigma})\to K_{\mathsf T}(\N_{\bf v, w}),$$ 
where $a=a'/a''$ were introduced for general Nakajima variety $\N_{\bf v, w}$ (see \cite{2012arXiv1211.1287M} for the cohomological version).  They depend on the \textit{slope} $s$,  a generic element of ${\rm Pic}(X)\otimes~\mathbb{R}$. After localization these maps become isomorphisms of vector spaces and thus one can define the following linear operators:
 \begin{equation}
R^s(a)={(\text{Stab}^s_{-}(a))}^{-1}\text{Stab}^s_{+}(a)\,.
\end{equation}
known as geometric R-matrices, which give rise to the trigonometric $R$-matrices in the $ADE$ case for a certain value of the slope $s$ . 
For a multidimensional $\bA=\prod_i\mathbb{C}^{\times}_{a_i}$-torus action,
 the construction of stable envelopes can be generalized by choosing a chamber 
 $\mathfrak{C}$ in the real Lie algebra of $\bA$, producing the maps
  $\{{\rm Stab}^s_{\mathfrak{C}}\}$ and $R$-matrices $\{R_{\mathfrak{C}, \mathfrak{C}'}^s(\{a_i/a_j\})\}$.  The explicit formulas for these maps in the case of $\N(n)$ and more generally, for flag varieties, can be found in Section 5 of \cite{Rimanyi:2014ef}. We also note, that in \cite{hernandez2022stable} an entirely algebraic procedure for the construction of stable envelopes is given, while at \cite{dedushenko2023interfaces} an alternative physics construction is proposed.

By RTT-procedure \cite{DF}, \cite{frt} the set of $R$-matrices
endows the vector spaces $K_{\mathsf T}(X)_{loc}$ for a Nakajima variety $\N_{\bf v, \bf w}$  with a structure of modules over certain Hopf  algebra, which we denote by $\mathcal{U}_{\hbar}(\frak{g}^Q)$, see Section 3 of \cite{Okounkov:2016sya} for the definitions. The Hopf structure gives the tensor product (\ref{tp}) a structure of module over $\mathcal{U}(\frak{g}^Q)$.

To see how that explicitly works in the case of $\N(n)$, for a given permutation $\pi \in \mathfrak{S}_{n}$ one can fix a chamber in the real Lie algebra of $\bA$ given in coordinates by $a_{\pi(1)}>\cdots>a_{\pi(n)}$. The slope $s$ is a generic element of ${\rm Pic}(\N(n))\otimes~\mathbb{R}\cong~\mathbb{R}$.
For every such choice of $\pi$ and $s$ we have canonical maps:
$$
\textrm{Stab}_{\pi}^{s}:\ \ K_{\mathsf T}(\N(n)^{\bA}) \longrightarrow K_{\mathsf T}(\N(n)),
$$
Then for a given transposition $\sigma=(k,m)$ one can define geometric R-matrices: 
$$
\mathcal{R}^{s}_{k,m}:=(\textrm{Stab}_{(k,m)}^{s})^{-1}
\circ \textrm{Stab}^{s}_{\textrm{id}},
$$
which by construction satisfy the quantum Yang-Baxter equations:
\begin{equation} \label{qybe}
\mathcal{R}^{s}_{i,j} \mathcal{R}^{s}_{i,k} \mathcal{R}^{s}_{j,k}=\mathcal{R}^{s}_{j,k} \mathcal{R}^{s}_{i,k} \mathcal{R}^{s}_{i,j}
\end{equation}
for an arbitrary triple $i,j,k$ and a slope $s$.

As in the previous section, let ${\rm V}(a)$ denote two-dimensional evaluation module
of quantum affine algebra $\mathcal{U}_{\hbar}(\widehat{\mathfrak{sl}}_2)$.  Let $e_1,e_2$ be the standard basis of ${\rm V}(a)$ corresponding to $1$ and $0$ weight vectors which in the language of spin chains these correspond to ``spin up'' and ``spin down'' states.
The corresponding $\fsh$-module given by the tensor product of two-dimensional evaluation modules 
$
{\rm V}(a_{1})\otimes \cdots \otimes {\rm V}(a_{n}).
$
has the \textit{standard  basis} in the weight $k$-subspace. It is labeled by $k$-subsets ${\bf p}=\{i_1,\dots, i_k\} \subset \{1,\dots, n\}$. Such a subset corresponds to the basis vector $v_{\bf p}=e_{k_1}\otimes\cdots \otimes e_{k_n}$
with $e_1$ at the positions $i\in {\bf p}$ and $e_2$ at  $i\in {\bf n}\setminus {\bf p}$. 

The following theorem gives geometric realization of quantum group \cite{Okounkov:2016sya},  \cite{Rimanyi:2014ef}:

\begin{Thm} \label{procan}
i) Let $s_{0}\in\epsilon [{\mathcal{O}}(1)]$ with $\epsilon \in (-1,0)$. Then the geometric $R$-matrix $\mathcal{R}^{s_0}_{l,l+1}$ evaluated in the basis of fixed points of (\ref{tp}) coincides with the trigonometric $\fsh$ $R$-matrix $R_{l,l+1}(a_{l}/a_{l+1})$ evaluated in the standard basis of $v_{\bf p}$. 	\\
ii) $\mathcal{U}({\frak{g}}^Q)\cong \fsh$ for $\normalfont \N_{\bf v, \bf w}=\N_{k,n}$. For every $n$ the map: 
$$
\normalfont
 {\rm V}(a_1) \otimes\cdots
\otimes {\rm V}(a_n) \to K_{\mathsf T}(\N(n)^{\bA})
$$
mapping a standard basis vector $e_{{\bf p}}$ to the class of the corresponding fixed point from $\normalfont K_{\mathsf T}(\N(n)^{\bA})$,
is an isomorphism of $\fsh$-modules, so that in particular, $\mathbb{F}(a)\cong {\rm V}(a)$.
\end{Thm}

The RTT procedure \cite{frt} constructs the corresponding quantum groups from $R$-matrices. The action of $\fsh$ in ${\rm V}(a_1) \otimes\cdots \otimes {\rm V}(a_n)$ in the standard basis coincides with its action on $K_{\mathsf T}(\N(n)^{\bA})$ in the basis of fixed points because the corresponding $R$-matrices coincide in these bases.  Other types of slopes 
correspond to  internal automorphisms of this representation given by operators of classical multiplication by the line bundles $\mathcal{O}(m)$. The action of this automorphism on the generators of $\fsh$ amounts to 
conjugation by $\mathcal{O}(m)$ iwhich shifts the loop index of Drinfeld's generators of $\fsh$ by $m$ units. These automorphisms act as the lattice elements of the quantum affine Weyl group 
of $\fsh$.


\subsection{Quasimaps, nonsingular and relative conditions}

Now we turn to the construction of the equivariant quantum K-theory. There are several version of this object. Our construction uses the quasimap theory instead of a canonical theory of stable maps.

First we define quasimaps to Nakajima variety. We will be very informal here: for a more detailed exposition one should follow original paper \cite{Ciocan-Fontanine:2011tg}, lecture notes  \cite{Okounkov:2015aa}, or \cite{Pushkar:2016qvw} for more detailed summary.

A quasimap $f$ 
\begin{equation}
f:{\mathcal{C}}\dashrightarrow \N_{\bf v, \bf w}\nonumber
\end{equation}

from $\mathcal{C}\cong \mathbb{P}^1$ to Nakajima variety $\N_{\bf v, \bf w}$ is:
\begin{itemize}
\item the collection of vector bundles $\mathscr{V}_i$ on $\mathcal{C}$ of ranks $\mathbf{v}_i$;\\ 
\item section of the bundle
\begin{equation}
f\in H^0(\mathcal{C},\mathscr{M}\oplus \mathscr{M}^{\ast}\otimes \hbar),\nonumber
\end{equation}
satisfying $\mu =0$, where
$$
{\mathscr{M}}=\sum_{i\in I}Hom(\mathscr{W}_i,\mathscr{V}_i)\oplus\sum_{i,j\in I}
{\mathscr Q}_{ij}\otimes Hom(\mathscr{V}_i,\mathscr{V}_j),
$$
\end{itemize}
so that $\mathscr{W}_i$ are trivial bundles of rank $\mathbf{w}_i$ and $\mu$ is the moment map. Here $\hbar$ is a trivial line bundle with weight $\hbar$ introduced to have the action of $\mathsf{T}$ on the space of quasimaps.
The degree $\bf{d}$ of a quasimap is a the vector of degrees of bundles $\mathscr{V}_i$.
and thus can be considered as an element of $H_2(\N_{\bf v, \bf w},\mathbb{Z})$. 

We say that a quasimap is stable if $f(p)$ is also a stable (belongs to the orbits giving rise to Nakajima variety) for all but finitely many points of $\mathcal{C}$. If $f(p)$ is
not stable, we will say that the quasimap $f$ is {\it singular} at $p$.

One can also define the {\it twisted} version of quasimaps. 
Let $\sigma: \mathbb{C}^{\times}\to A$ be a cocharacter of subtorus $A=\text{Ker}\, \hbar\subset T$, preserving symplectic form, which determines the twist of a quasimap to $\N_{\bf v, \bf w}$. As $T$ is acting on $W_i$,  $\sigma$ determines $\mathscr{W}_i$ over the curve $\mathcal{C}$ as bundles associated to $\mathscr{O}(1)$.

For a point on the curve $p\in \mathcal{C}$ we have an evaluation map from the moduli space of stable quasimaps to the quotient stack $\text{ev}_p : \textsf{QM}^{\ard} \to L_{\mathbf{v},\mathbf{w}}/G$ defined by
$\text{ev}_p(f)=f(p)$. Note that the quotient stack contains $\N_{\bf v, \bf w}$ as an open subset corresponding to locus of semistable points: $$\N_{\bf v, \bf w}={\mu}_{ss}^{-1}(0)/G\subset L_{\mathbf{v},\mathbf{w}}/G.$$
A quasimap $f$ is called nonsingular at $p$ if $f(p)\subset X$. 

In short, we conclude that the open subset ${{\textsf{QM}}^\ard}_{\text{nonsing  p}}\subset {{\textsf{QM}}^\ard}$ of stable quasimaps of degree $\bf{d}$, nonsingular at the given point $p$,
is endowed with a natural evaluation map:
\begin{equation}
{{\textsf{QM}}^\ard}_{\text{nonsing }\, p} \stackrel{{\text{ev}}_p}{\longrightarrow} \N_{\bf v, \bf w}
\end{equation}
which sends a quasimap to its value at $p$. The moduli space of relative quasimaps ${{\textsf{QM}}^\ard}_{\text{relative} \, p}$ is a resolution of ${\text{ev}}_p$ (or compactification), meaning we have a commutative diagram:
\begin{center}
\begin{tikzpicture}[node distance =5.1em]

  \node (rel) at (2.5,1.5) {${{\textsf{QM}}^\ard}_{\text{relative}\, p}$};
  \node (nonsing) at (0,0) {${{\textsf{QM}}^\ard}_{\text{nonsing}\,  p}$};
  \node (X) at (5,0) {$\N_{\bf v, \bf w}$};
  \draw [->] (nonsing) edge node[above]{$\text{ev}_p$} (X);
  \draw [->] (rel) edge node[above]{$\widetilde{\text{ev}}_p$} (X);
  \draw [right hook->] (nonsing) edge  (rel);
    \end{tikzpicture}
\end{center}
with a proper evaluation map $\widetilde{\text{ev}}_p$ from ${{\textsf{QM}}^\ard}_{\text{relative}\, p}$ to $\N_{\bf v, \bf w}$. The construction of this resolution and the moduli space of relative quasimaps is explained in \cite{Okounkov:2015aa}. The main idea of this construction is to allow the base curve to change in cases, when the relative point becomes singular. When this happens we replace the relative point by a chain of non-rigid projective lines, such that the endpoint and all the nodes are not singular. We refer to the original curve as {\it parametrized component} of relative quasimap. Similarly, for nodal curves, we do not allow the singularity to be at the node, and if that happens we instead paste in a chain of non-rigid projective lines.
The properness is needed for the pushforward in the K-theory along evaluation map to be well-defined, see below. Let us specialize again to the case of $\N_{\bf v, \bf w}=\N_{k,n}$ for simplicity.

The construction above can be generalized to introduce the following spaces of quasimaps relative at the points $p_1,\dots, p_m$, and nonsingular at the points $r_1, \dots, r_s$:
$$
{\textsf{QM}}^\ard_{{{\rm{nonsing}}, r_1,\cdots,r_s } \atop {{{\rm{relative}},  p_1,\dots,p_m}}}\subset {\textsf{QM}}^\ard_{ {{{\rm{relative}},  p_1,\dots,p_m}}},
$$ 
where the former space is the open subset of the latter. In the first space $f(r_i)$ are stable for all $i=1,\dots, s$.  In the algebro-geometric language both spaces above are Deligne-Mumford stacks of finite type with perfect obstruction theory and the latter space is a stack proper over affine space over $\N^0_{k,n}$.

While the evaluation maps 
$$
{\rm{ev}}_{p_i}: {\textsf{QM}}^\ard_{{{\rm{nonsing}}, r_1,\cdots,r_s } \atop {{{\rm{relative}},  p_1,\dots,p_m}}}\rightarrow \N_{k,n}, \ \ \ 
{\rm{ev}}_{r_j}: {\textsf{QM}}^\ard_{{{\rm{nonsing}}, r_1,\cdots,r_s } \atop {{{\rm{relative}},  p_1,\dots,p_m}}}\rightarrow \N_{k,n}
\ \ \ 
$$
for $i=1,\dots,m, j=1,\dots,s$ are well-defined, ${\textsf{QM}}^\ard_{{{\rm{nonsing}}, r_1,\cdots,r_s } \atop {{{\rm{relative}},  p_1,\dots,p_m}}}$ is not necessarily proper over $\N^0_{k,n}$ and thus ${\rm{ev}}_{r_j}$ does not necessarily provide push-forward maps in equivariant K-theory.

To define pushforwards in the latter case,  we introduce another equivariant variable.  
Let  $\C^{\times}_q$ be a 1-dimensional torus acting  on $\mathbb{P}^1$ in the following way: if $[x:y]$ denotes homogeneous coordinates on $\mathbb{P}^1$ then the corresponding action of $\C^{\times}_q$ is defined by:
$$
[x:y] \rightarrow [x q :y] =[x :y q^{-1}].
$$ 
The set  ${\mathbb{P}^1}^{\C^{\times}_q}$ consists of two fixed points: 
\begin{eqnarray} \label{points}
p_1=[1,0] =0 \in \mathbb{P}^{1},  \ \ p_2=[0,1]=\infty \in \mathbb{P}^{1}.
\end{eqnarray}
Therefore a $\C^{\times}_{q}$-fixed quasimap may only have special marked points on $\mathbb{P}^1$  (nonsingular or relative)  at $p_1$ or $p_2$. 
The ${\textsf{QM}}^\ard$-moduli spaces have a natural action of maximal torus $\mathsf{T}$, lifting its action from $\mathbb{N}_{k,n}$. Thus we obtain full torus $\mathsf{T}_q=\mathsf{T}\times \mathbb{C}^{\times}_q$.

We will need the following statements about pushforward maps proved in Section 7.2 of \cite{Okounkov:2015aa}:
\begin{Thm}
i)	For $p_1,p_2$ as in (\ref{points}) the evaluation maps 
\begin{equation}
\normalfont\mathrm{ev}_{p_1}: ( {\normalfont\textsf{QM}}^\ard_{\mathrm{nonsing}, p_1 } )^{\C^{\times}_q} \rightarrow \N_{k,n}, \quad 
	\mathrm{ev}_{p_1}: ({\normalfont \textsf{QM}}^\ard_{{\mathrm{nonsing}, p_1} \atop {\mathrm{relative}, p_2} } )^{\C^{\times}_q} \rightarrow \N_{k,n}
\end{equation}
	are proper. \\
ii)	The following pushforward maps are well-defined:
$$
\normalfont\mathrm{ev}_{p_1,*}: K_{\mathsf{T}_q}({\normalfont{\textsf{QM}}}^\ard_{\mathrm{nonsing}, p_1 } ) \rightarrow K_{\mathsf{T}_q}(\N_{k,n})_{loc}
$$	
$$
\normalfont\mathrm{ev}_{p_1,*} \otimes \mathrm{ev}_{p_2,*}: K_{\mathsf{T}_q}({\normalfont\textsf{QM}}^\ard_{{\mathrm{nonsing}, p_1} \atop {\mathrm{relative}, p_2} } ) \rightarrow K_{\mathsf{T}_q}(\N_{k,n})^{\otimes 2}_{loc}
$$ 
\end{Thm}

\subsection{Capped tensors and the gluing formula}

Let us consider the following map:
$$
\textrm{ev}_{p_1 *} \otimes \dots\otimes \textrm{ev}_{p_m *}: K_{\mathsf{T}}\Big( {{\textsf{QM}}^\ard}_{\text{relative }\, p_1,\dots,p_m} \Big)
\rightarrow K_{\mathsf{T}}(\N_{k,n} )^{\otimes m}.
$$
which is well-defined since the corresponding evaluation maps are proper. 
Let 
${\widehat{\mathcal{O}}}_{\rm vir}^\ard\in K_{\mathsf T}\Big( {{\textsf{QM}}^\ard}_{\textsf{relative }\, p_1,\dots,p_m} \Big)$ 
denote the twisted (we discuss in more detail later) virtual structure sheaf of the moduli space provided by the perfect obstruction theory of ${{\textsf{QM}}^\ard}_{\text{relative }\, p_1,\dots,p_m}$. 

The  power series with coefficients in $K$-theory of $\N_{k,n}$	
\begin{equation} \label{cten}
\hat{V}_{p_1,\dots,p_m}(z):=\sum\limits_{d=0}^{\infty}\, z^d\,  \mathrm{ev}_{p_1 *} \otimes 
\dots\otimes \mathrm{ev}_{p_m *}\Big({\widehat{\mathcal{O}}}_{\rm vir}^\ard\Big)\in  K_{\mathsf T}(\N_{k,n} )^{\otimes m}[[z]].
\end{equation}
are called {\it capped $K$-theoretic tensors}. 

Let $p \in \mathbb{P}^1$ be a point disjoint from $p_1,\dots,p_m$ and let $\tau$ be a Schur functor with coefficients in $K_{\mathsf T}(pt)$. The value $f(p)$ for such point is not necessary stable, therefore it does not provide an evaluation map to $\N_{k,n}$, but only to the quotient stack:
$$
\textrm{ev}_{p}: {{\textsf{QM}}^\ard}_{\text{relative }\, p_1,\dots,p_m}\rightarrow [\mu^{-1}(0)/GL(k)].
$$  
Let $\mathcal{V}_{stack}=(V \times T^*M)/GL(k)$ be the element in $K$-theory class of this stack associated to $k$-dimensional $GL(k)$-module $V$. We define 
$$
\tau(\left.\mathscr{V}\right|_{p}):=\textrm{ev}_{p}^{*}( \tau(\mathcal{V}_{stack}) ). 
$$
Then the power series 	
\begin{equation}\label{phfd}
\widehat{V}^{(\tau)}_{p_1,\dots,p_m}(z):= \sum\limits_{d=0}^{\infty}\, z^d\,  \mathrm{ev}_{*,p_1} \otimes \dots\otimes \mathrm{ev}_{*,p_m}\Big({{\widehat{\mathcal{O}}}_{\rm vir}^\ard \otimes \tau(\left.\mathscr{V}\right|_{p})}\Big) 
\end{equation}
is called \textit{capped K-theoretic tensor with descendent $\tau$ inserted at $p$}. 
Unless stated otherwise, we use $p=\{0\}$ on the parametrized component 
$\mathcal{C}=\mathbb{P}^1$ of the quasimap. 

Now let us discuss a tensor of rank $m$  produced by the pushforwards (\ref{phfd}) in more details. In the basis of fixed points of $K_{\mathsf T}(\N_{k,n})$ we can describe such tensor explicitly by its components $T_{i_1,...,i_m} (\mathcal{C})$. If $\mathcal{C}$ degenerates to a union of two rational curves $\mathcal{D}=\mathcal{C}_1\cup_p \mathcal{C}_2$ with one common nodal point $p$, the so-called \textit{degeneration formula} in the theory of quasimaps states that the tensors for quasimaps with domain $\mathcal{D}$ and its degeneration $\mathcal{C}_1\cup_p \mathcal{C}_2$ are equal:
$$
T_{i_1,...,i_m} (\mathcal{C})=T_{i_1,...,i_m} (\mathcal{C}_1\cup_p \mathcal{C}_2).
$$ 

The statement of the \textit{gluing formula} gives the factorization of tensors (\ref{phfd}) corresponding to  parametrized component  $\mathcal{C}_1\cup_p \mathcal{C}_2$ into a product of tensors given by quasimaps with parametrized components $\mathcal{C}_1$ and $\mathcal{C}_2$. Assume that after degeneration $p_1,.\cdots,p_s \in \mathcal{C}_1$ while the rest of the points lie within the second component, i.e. $p_{s+1},\dots, p_m \in \mathcal{C}_2$.  This gives rise to two  tensors for first and second components correspondingly: i) on $\mathcal{C}_1$ we have ${T}_{i_1,\dots,i_s,a}(\mathcal{C}_1)$  with $s+1$ indices, so that the index $a$ corresponds to the new point $p$; ii) on $\mathcal{C}_2$ we have a tensor 
${T}_{i_{s+1},\dots,i_{m},b}(\mathcal{C}_2)$, were $b$ is the index corresponding to the relative point $p$. Then gluing theorem states that there exists the \textit{gluing matrix} ${\bold{G}} \in  K^{\otimes 2}_{\mathsf T}(\N_{k,n})[[z]]$ such that
\begin{equation} \label{gluingth}
{T}_{i_1,\dots,i_m}(\mathcal{C}_1\cup_p \mathcal{C}_2)=\sum\limits_{a,b}\, {(\bold{G^{-1})}^{{a,b}}} {T}_{i_1,\dots,i_s,a}(\mathcal{C}_1)\, {T}_{i_{s+1},\dots,i_{m},b}(\mathcal{C}_2).
\end{equation}

It turns out that the gluing matrix is given by the following pushforward: 
$$
{\bold{G}}=\sum\limits_{d=0}^{\infty}\, z^d {\text{ev}}_{p_1,*} \otimes {\text{ev}}_{p_2,*} \left({\widehat{\mathcal{O}}}_{\rm vir}^\ard \right). 
$$

\subsection{Capping operator and difference equations}
The simplest capped tensor one can use is as follows:
$$
\widehat{V}^{(\tau)}(z)=\sum\limits_{\ard=0}^{\infty} z^{\ard} {\rm{ev}}_{p_2, *}\Big(\widehat{{\mathcal{O}}}^{\ard}_{{\rm{vir}}} \otimes\tau (\left.\mathscr{V}_{\{i\}}\right|_{p_1}) \Big) \in  K_{\mathsf{T}_q}(\N_{k,n})[[z]],
$$
where $\widehat{{\mathcal{O}}}^{\ard}_{{\rm{vir}}}$ is a structure sheaf on  $\textsf{QM}^{\ard}_{{\rm{relative}} \, p_2}$. 
We will refer to it as {\it capped vertex} with descendants $\tau$. Since the relative pushforward is well-defined, we can set 
$$
\mathcal{\tau}(z):=\lim_{q\rightarrow 1}\widehat{V}^{(\tau)}(z)\in  K_{\mathsf{T}}(\N_{k,n})[[z]]
$$
We will refer to these series as quantum tautological class for $\tau(V)\in K_{\mathsf{T}}(X)$. 
At the same time, the element
$$
V^{(\tau)}(z)=\sum\limits_{\ard=0}^{\infty} z^{\ard} {\rm{ev}}_{p_2, *}\Big(\widehat{{\mathcal{O}}}_{{\rm{vir}}} \otimes\tau (\left.\mathscr{V}_{\{i\}}\right|_{p_1}) \Big) \in  K_{\mathsf{T}_q}(\N_{k,n})_{loc}[[z]],
$$
where $\widehat{{\mathcal{O}}}_{{\rm{vir}}}$ is a structure sheaf on  $\textsf{QM}^{\ard}_{{\rm{nonsing}} \, p_2}$
are the series with coefficients in localized K-theory and is called {\it  bare vertex} with descendant $\tau$. 

The operator which relates capped and bare vertices, is known as {\it capping operator} (see Section 7.4.4. of \cite{Okounkov:2015aa}) and is defined as the following class in the localized K-theory:
\begin{equation}
\Psi(z)=\sum\limits_{\ard=0}^{\infty}\, z^{\ard} {\rm ev}_{p_1,*}\otimes{\rm ev}_{p_2,*}
\Big( {\it \widehat{\mathcal{O}}_{{\rm{vir}}}} \Big)
\in K^{\otimes 2}_{\mathsf{T}_q}(\N_{k,n})_{loc}[[z]],
\end{equation}
where $\widehat{{\mathcal{O}}}_{{\rm{vir}}}$ is a structure sheaf on $\textsf{QM}^{\ard}_{{{\rm{relative}\, p_1}} \atop {{\rm{nonsing}\, p_2}}}$.
One can define the following bilinear form $(\cdot, \cdot)$ on $K_{\mathsf{T}}(\N_{k,n})$ via the equivariant Euler characteristic 
\begin{equation}
(\mathcal{F},\mathcal{G})= \chi\left(\mathcal{F}\otimes\mathcal{G}\otimes K^{-\frac{1}{2}}\right)\,,
\end{equation}
twisted by a root of the canonical class, which exists for Nakajima quiver varieties. This twisting is introduced due to the similar twisting of the structure sheaf $\widehat{{\mathcal{O}}}_{{\rm{vir}}}$ (see \cite{Okounkov:2015aa}). 
Because of this pairing, $\Psi(z)$ can be considered as an operator acting from the second to the first copy of $K_{\mathsf{T}_q}(\N_{k,n})_{loc}[[z]]$. 

Using this operator intepretation of $\Psi(z)$, the relation between the capped vertex with descendent $\tau$ and bare vertex 
is given by the following formula:
\begin{equation}
\label{verrel}
\widehat{V}^{(\tau)}(z) = \Psi(z) {V}^{(\tau)}(z).
\end{equation}
This capping operator as well as bare and capped vertices could be defined for general Nakajima variety $\N_{\bf v, w}$. Upon suitable normalization:  ${\Phi}(a,z)=\Psi(a,z)\Theta(a)$ 
by a certain K-theory class $\Theta(a)$,  it is known to satisfy difference equations (see Proposition 6.5.30, Section 8.2.17,  Theorem 8.2.20 and Corollary 8.2.24 of \cite{Okounkov:2015aa}):
\begin{eqnarray}\label{qkzgeom}
\Phi(a, q^{\mathcal{L}_i}z)={\bf M}_{\mathcal{L}_i}\Phi( a,z);\quad  \Phi(q^{\sigma}a, z)={\bf S}_{\sigma}\Phi(a, z)
\end{eqnarray}
  where $\mathcal{L}_i$ are the line tautological bundles 
  ${\mathcal{L}_i}=\Lambda^{{\bf v}_i}V_i$, the operator ${\bf S}_{\sigma}$ is given by:
  $$
{\mathbf S}_{\sigma}(a,z)=
\sum\limits_{\ard=0}^{\infty}\, z^{\ard}\,
{\rm ev}_{{p_1},*}\otimes {\rm ev}_{{p_2},*}\left( \widehat{{\mathcal{O}}}_{{\rm{vir}}} \right) {\mathbf G}^{-1}\in K_{\mathsf{T}_q}^{\otimes 2}(X)_{loc}[[z]], 
$$
and $\widehat{{\mathcal{O}}}_{{\rm{vir}}}$ is the virtual structure sheaf on the space of twisted quasimaps: $\textsf{QM}^{~\sigma,\ard}_{{{\rm{relative}\, p_1}} \atop {{\rm{relative}\, p_2}}}$. 
Here we also use the notation $(q^{\mathcal{L}_i}z)^{\ard}=q^{\langle \mathcal{L}_i, \ard \rangle} z^{\ard}$, where $\ard\in H_2(\N_{{\bf v},{\bf w}}, \mathbb{Z})$, $\mathcal{L}_i\in Pic(\N_{{\bf v},{\bf w}}).$
For our main example $\mathbf{N}_{k,n}={\rm T^* Gr}_{k,n}$ there exist just one operator 
$\mathbf{M}$ associated with tautological bundle $\mathcal{O}(1)=\Lambda^kV$. 
Without normalization, in this case the corresponding difference equation involving operator ${\bf M}$ is as follows:
\begin{equation}
	\label{difference}
	\Psi(a,qz)={\bf M}(a,z) \Psi(a, z) \mathcal{O}(1)^{-1},
\end{equation}	
We will return to the properties of the operator ${\bf M}$ and this difefrence equations in the next subsection.

Let us discuss now this system of equations in general. It turns out that it coincides with quantum Knizhnik-Zamolodchikov equations and the related dynamical equations. 
In particular, it turns out that the following Theorem is true (see  Theorem 9.3.1 from  \cite{Okounkov:2015aa} as well as Theorem 3.1 of \cite{Okounkov:2016sya})):

\begin{Thm}
The shift operator $\mathbf S_{\sigma}$ is related to the R-matrix 
$R(a)$ in the following way
$$
\tau_\sigma^{-1}z^{{\mathbf v}'}R(a)=\text{Stab}_{+}^{-1}\tau_{\sigma}^{-1}\mathbf S_{\sigma}(a,z)\text{Stab}_{+},
$$
where $\tau_{\sigma}f(a)=f(q^{\sigma} a)$, identifying the second difference equation  (\ref{qkzgeom}) with the quantum Knizhnik-Zamolodchikov equation.
\end{Thm} 
\subsection{Quantum equivariant K-theory ring}

Now we will discuss how the quantum K-theory ring emerges from equivariant prushforwards. 

Let $\textsf{QM}^{\ard}_{ \,p_1,p_2,p_3}$ be the moduli space of quasimaps from ${\mathbb{P}}^1$ with 3 relative points 
and let $\mathcal{O}^\ard$ be the virtual structure sheaf on this moduli space.  For a given a class $\mathcal{F}\in K_{\mathsf T}(\N_{k,n})$, we construct the following equivariant pushforward:
\begin{equation}
\label{qprod}
\mathcal{F}\circledast:= 
\left(\sum\limits_{\ard=0}^{\infty}\,z^{\ard} {\text{ev}}_{p_1\ast}\times {\text{ev}}_{p_3\ast} \left(\text{ev}^{\ast}_{p_2}({\bold{G}}^{-1}\mathcal{F})\otimes \mathcal{O}_{\rm vir}^\ard \right)\right) {\bold{G}}^{-1}.
\end{equation}
By definition, $\mathcal{F}\circledast$ is a rank two tensor, which, thanks to the scalar product in K-theory, can be identified with the linear operator:
$$
\mathcal{F}\circledast\in End({K_{\mathsf T}(\N_{k,n})})[[z]].
$$
Note, that the moduli space of degree zero quasimaps is isomorphic to $\N_{k,n}$, since a degree zero quasimap maps the entire curve to a single point in $\N_{k,n}$. This implies that the zeroth coefficient
of the power series (\ref{qprod}) has the form: 
\begin{equation} \label{relcla}
\mathcal{F} \circledast  =\mathcal{F} \otimes +\dots ,
\end{equation}
where $\mathcal{F} \otimes$ is the operator of tensor multiplication by $\mathcal{F}$ in the equivariant $K$-theory,
where dots in (\ref{relcla}) stand for the terms vanishing in $z\rightarrow 0$ limit, which could be thought of classical limit in this quantum K-theoretic sense.

We will refer to  $QK_{\mathsf T}(\N_{k,n}):=K_{\mathsf T}(\N_{k,n})[[z]]$ endowed with the multiplication~(\ref{qprod}) the  {\it quantum equivariant $K$-theory ring} of $\N_{k,n}$.

Let $\hat{\bold{1}}(z)\in QK_{\mathsf T}(\N_{k,n})$  be the quantum tautological bundle for $\tau={\mathcal{O}}_{\N_{k,n}}$:
\begin{equation}
\label{mulid}
\hat{\bold{1}}(z)= \sum\limits_{\ard=0}^{\infty} z^d \textrm{ev}_{p_2, *}\Big( \mathcal{O}_{\rm vir}^{\ard}   \Big),
\end{equation}
where $\mathcal{O}_{\rm vir}^{\ard}$ above is the structure sheaf on $\textsf{QM}^{\ard}_{\textrm{relative} \, p_2}$. 
We proved the following in \cite{Pushkar:2016qvw}:

\begin{Thm}
	The quantum $K$-theory ring $\normalfont QK_{\mathsf T}(\N_{k,n})$ is a commutative, associative unital algebra, where the unit element is given in (\ref{mulid}).
\end{Thm}

Let $\widehat{\mathcal{O}(1)}(z)$ be the quantum tautological bundle for $\tau=\mathcal{O}(1)=\det\mathcal{V}$. It turns out that under there exist a limit  $q\rightarrow 1$  the operator ${\bf M}(z)$ coincides with the operator of quantum multiplication by the quantum line bundle:
$$
{\bf M}(z)|_{q\rightarrow 1}=\widehat{\mathcal{O}(1)}(z) \circledast.
$$
We remind that localized $K$-theory $K_{\mathsf T}(\N_{k,n})_{loc}$ is a vector space over the field ${\mathcal{A}}$ with a natural basis given by the $K$-theory classes of fixed points ${\mathcal{O}}_{\fp}$, where $\fp=\{i_1,\cdots,i_k\}\subset \{1,\cdots,n\}$ is a $k$-subset labeling the elements of $\N_{k,n}^{\mathsf T}$. Thus the  operators of classical multiplication by $\mathcal{F}\in K_{\mathsf T}(\N_{k,n})$ act on this vector space as ${\mathcal{A}}$-linear operators, i.e. $\mathcal{F}\in {\rm End}_{{\mathcal{A}}}(K_{\mathsf T}(\N_{k,n})_{loc})$. These operators are diagonal in the basis of fixed points ${\mathcal{O}}_{\fp}$:
$$
\mathcal{F}\otimes {\mathcal{O}}_{\fp}=\left.\mathcal{F}\right|_{\fp} {\mathcal{O}}_{\fp}
$$ 
where $\left.\mathcal{F}\right|_{\fp} \in {{\mathcal{A}}}$ denotes the restriction of K-theory class to the corresponding fixed point $\fp$. 

The operators of quantum multiplication by $K$-theory classes are defined as power series 
$\mathcal{F}\circledast \in {\rm End}_{{{\mathcal{A}}}}(K_{\mathsf T}(\N_{k,n})_{loc})[[z]]$. Therefore one could study the spectral problem, namely finding classes $\psi_{\fp}(z) \in K_{\mathsf T}(\N_{k,n})[[z]]$ and power series $v_{\fp}(z)\in{{\mathcal{A}}}[[z]]$ with fixed degree zero terms:
\begin{eqnarray} \label{norpow}
&&\psi_{\bf p}(z)={\mathcal{O}}_{\bf p}+O(z), \quad v^{\mathcal{F}}_{\bf p}(z)=\left.\mathcal{F}\right|_{\bf p} +O(z)\\
&& {\rm such~ that}\quad \mathcal{F}\circledast \psi_{\bf p}(z)= v^{\mathcal{F}}_{\bf p}(z) \psi_{\bf p}(z).\nonumber 
\end{eqnarray}
The power series $\psi_{\bf p}(z)$ normalized as in (\ref{norpow}) do not depend on $\mathcal{F}$, because the quantum multiplication is commutative and thus $\psi_{\bf p}(z)$ are the  eigenvectors of quantum multiplication and the power series $v^\mathcal{F}_{\bf p}(z)$ are the corresponding eigenvalues.  

The operator of classical multiplication by the line bundle $\mathcal{O}(1)$ in the basis of fixed points is given by the diagonal matrix with the eigenvalue $\lambda^{0}_\fp=a_{i_1}\cdots a_{i_k}$ corresponding to a fixed point $\fp=\{i_1,\cdots,i_k\}$. The quantum line bundle is  a $z$-deformation of this class: $\widehat{\mathcal{O}(1)}(z)=\mathcal{O}(1)+\dots$, where dots stand for terms vanishing in $z\rightarrow 0$ limit.  Hence, the eigenvalues of the quantum line bundle $\widehat{\mathcal{O}(1)}(z)$ have the form $\lambda_{\bf p}(z)=\lambda^{0}_{\bf p}+\lambda^{1}_{\bf p} z +\cdots$.

The $q$-difference equation (\ref{difference}) in the basis of fixed points is very useful to compute the eigenvalues and eigenvectors of quantum multiplication explicitly. By construction, the capping operator is a $K_{\mathsf T}(\N_{k,n})^{\otimes 2}$ -valued power series of the form:
$$
\Psi(0)=I+\Psi_1 z +\Psi_2 z^2 +\dots,
$$
with the first term being the identity matrix.  
Therefore, at $z=0$ the equation  (\ref{difference}) holds trivially because
$\M(0)=\mathcal{O}(1)$. One obtains the higher terms by solving corresponding linear the linear problem.  In \cite{Pushkar:2016qvw}, \cite{Koroteev:2017aa} we performed such analysis and obtained that the eigenvalues of quantum multiplication by tautological bundles can be expressed in terms of bare vertex functions, as stated in the Theorem below.


\begin{Thm}
	The eigenvalues of the operator of quantum multiplication by $\hat{\tau}(z)$ can be expressed as follows:
	\begin{equation}
	\label{eval}\normalfont
	\tau_{\fp}(z)=\lim\limits_{q \rightarrow 1 } \frac{V^{(\tau)}_{\fp}(z)}{V^{(1)}_{\fp}(z)}
	\end{equation}
\end{Thm}

\subsection{Bare vertex with descendants and Bethe ansatz equations}
As we have seen in the previous section, the bare vertex functions give rise to  
the eigenvalues of the quantum tautological bundles. These bare vertex functions  are  computed using equivariant localization in K-theory. It turns out that these functions are given by some standard $q$-hypergeometric series expansions, so that the limit (\ref{eval}) provides an explicit formula for the eigenvalues of quantum tautological bundles. As we shall see, these eigenvalues coincide with the symmetric polynomials of the solutions of Bethe equations. 

The moduli spaces of nonsingular and relative quasimaps have a perfect deformation-obstruction theory \cite{Ciocan-Fontanine:2011tg}. 
In particular, this leads to the construction of a tangent virtual bundle $T^{\textrm{vir}}$, a virtual structure sheaf $\mathcal{O}_{\rm vir}^{\ard}$ and a virtual canonical bundle over these moduli spaces. Let $(\mathscr{V},\mathscr{W})$ be the data defining a point on ${\textsf{QM}^{\ard}}_{\text{nonsing}\, p}$. We define the fiber of the reduced virtual tangent
bundle to ${\textsf{QM}^{\ard}}_{\text{nonsing }\, p}$ at this point to be equal to:
\begin{equation}
T^{\textrm{vir}}_{(\mathscr{V},\mathscr{W})} 
\textsf{QM}^{\ard}_{\text{nonsing  p}}=H^{\bullet}(\mathscr{P}\oplus\hbar\, \mathscr{P}^{*} ) - T_{f(p)} \N_{k,n},
\end{equation}
where $T_{f(p)} \N_{k,n}$ is a normalizing term and  $\mathscr{P}$ is the polarization bundle on the curve $\mathcal{C}$:
$$
\mathscr{P}=\mathscr{W}^{*} \otimes \mathscr{V} - \mathscr{V}^{*} \otimes \mathscr{V}.
$$

We use here the reduced virtual tangent space which differs from standard one by subtracting $T_{f(p)} \N_{k,n}$. Namely, if $\mathcal{V}$, $\mathcal{W}$ are the tautological bundles on $\N_{k,n}$, then the tangent bundle has the form: $T \N_{k,n} = \mathcal{P} + \hbar \mathcal{P}^* \ \  \textrm{for} \ \  \mathcal{P}=\mathcal{W}^{*}\otimes \mathcal{V} - \mathcal{V}^{*} \otimes \mathcal{V}.$ This term does not depend on $d$ and thus produces a simple multiple in the vertex function: it is the multiple normalizing the  vertex such that $V^{(\tau)}_{\fp}(0)=\tau$.

The symmetrized virtual structure sheaf is defined by:
\begin{equation}
\label{virdef}
\mathcal{O}^{\ard}_{\rm vir}=\mathcal{O}^{\ard}_{\textrm{vir}}\otimes {\mathscr{K}^{1/2}_{\textrm{vir}}} q^{\deg(\mathscr{P})/2},
\end{equation}
where $\mathcal{O}^{\ard}_{\textrm{vir}}$ is a standard structure sheaf and $\mathscr{K}_{\textrm{vir}}={\det}^{-1}T^{\textrm{vir}} \textsf{QM}^{\ard}$ is the virtual canonical bundle.

Recall that the degree $\ard$ of quasimaps to $\N_{k,n}$ is given by degree of rank $k$ bundle  $\mathscr{V}$ on ${\mathbb{P}}^1$. Let us consider the set of $\mathsf T$-fixed points, namely  $(\mathscr{V},\mathscr{W}) \in  (\textsf{QM}^{\ard}_{\text{nonsing} \, p_2})^{\mathsf T}$ such that the value of the evaluation map at $p_2$ is $\fp \in \N_{k,n}^{\mathsf T}$. Since the virtual tangent space is a representation of the torus $\mathsf T$, let us denote 	its $\mathsf T$-character by:
\begin{equation} \label{pcontr}
\chi(\ard)=\textrm{char}_{\mathsf T}\Big( T^{\rm vir}_{(\mathscr{V},\mathscr{W})} \textsf{QM}^{\ard}\Big).
\end{equation}
Localization in $K$-theory gives the following formula for the equivariant pushforward:
\begin{eqnarray}\label{vertexexp}
V^{(\tau)}_{\fp}(z)=\sum_{\ard=0}^{\infty}\sum\limits_{(\mathscr{V},\mathscr{W}) \in  (\textsf{QM}^{\ard}_{\text{nonsing} \, p_2})^{\mathsf T}}\, \hat{s}(  \chi(\ard) )\, z^{\ard} q^{\deg(\mathscr{P})/2} \tau (\left.\mathscr{V}\right|_{p_1}),
\end{eqnarray}	
where the sum runs over the $\mathsf T$-fixed quasimaps which take value $\fp$ at the nonsingular point $p_2$ and the function $\hat{s}$ is defined as follows:
$$	
\hat{s}(x)=\frac{1}{x^{1/2}-x^{-1/2}}, \ \ \ \hat{s}(x+y)=\hat{s}(x)\hat{s}(y).
$$
Note, that the tangent weight contribute to vertex via the roof function 
$\hat{s}(x)$ because the symmetrized virtual structure sheaf (\ref{virdef})is defined together with a shift on the square root of canonical bundle 
${\mathscr{K}^{1/2}_{\textrm{vir}}}$.

The contributions to the formula (\ref{vertexexp}) are computed in the following way. The polarization  bundle $\mathscr{P}$ on $\mathbb{P}^1$  corresponding to a  $\mathsf T$-fixed point on $\textsf{QM}^{d}_{\textrm{nonsing} \, p_2}$ splits into a sum of $\mathsf T$-equivariant line bundles $\mathscr{P}=\bigoplus_{i} a_i q^{-d_i} \mathcal{O}(d_i)$ for some characters $a_i$ of the framing torus $\bA$. Explicit computation gives the following  contribution of $\mathsf T$-characters:
\begin{equation}	
\textrm{char}_{\mathsf T}\Big(H^{\bullet}(a_i q^{-d_i} \mathcal{O}(d_i))\Big)=a_iq^{-d_i}(d_i+1)_q, \quad {\rm where} \quad (n)_q=\frac{q^n-1}{q-1}.
\end{equation}	
	
	
Using the notation 
\begin{eqnarray}
\varphi(x)=\prod^{\infty}_{i=0}(1-q^ix),\quad  \{x\}_{d}=\dfrac{(\hbar/x,q)_{d}}{(q/x,q)_{d}} \, (-q^{1/2} \hbar^{-1/2})^d, \ \ \textrm{where}  \ \ (x,q)_{d}=\frac{\varphi(x)}{\varphi(q^dx)}.\nonumber
\end{eqnarray}
If $\fp=\{x_1,\cdots,x_k\} \subset \{a_1,\cdots,a_n\}$ be a $k$-subset defining a torus fixed point $\fp\in \N^{\mathsf T}_{k,n}$. The vertex function component corresponding to that point is given by the following $q$-hypergeometric function:
	$$
	V^{(\tau)}_{\fp}(z) = \sum\limits_{d_1,\cdots, d_k\in {\mathbb{Z}_{\ge 0}}}\, z^{d} q^{n d/2}\, \prod\limits_{i,j=1}^{k}\{x_i/x_j\}^{-1}_{d_i-d_j}   \prod\limits_{i=1}^{k} \prod\limits_{j=1}^{n} \{x_i/a_j\}_{d_i} \tau(x_1 q^{-d_1},\cdots, x_k q^{-d_k}),
	$$
	where $d=\sum^k_{i=1}d_i$.
This function has a  Mellin-Barnes type integral representation: 
\begin{eqnarray}
\label{verint}
	&&V^{(\tau)}_{\fp}(z)= \nonumber\\
	&& \dfrac{1}{2 \pi i \alpha_p} \int\limits_{C_{\fp}} \prod\limits_{i=1}^{k} \dfrac{d s_i}{s_i} \, e^{-\frac{\ln(z_{\sharp}) \ln(s_i) }{\ln(q)}} \,  \prod\limits_{i,j=1}^{k} \dfrac{\varphi\Big( \frac{s_i}{ s_j}\Big)}{\varphi\Big(\frac{q}{\hbar} \frac{s_i}{ s_j}\Big)} \prod\limits_{i=1}^{n}\prod\limits_{j=1}^{k} \dfrac{\varphi\Big(\frac{q}{\hbar} \frac{s_j}{ a_i }\Big)}{\varphi\Big(\frac{s_j}{ a_i }\Big)}
	\tau(s_1,\cdots,s_k).
\end{eqnarray}
Here  the contour of integration $C_{\fp}$, corresponding to a fixed point $\fp=\{x_1,\cdots,x_k\}\subset\{a_1,\cdots,a_n\}$
	is a positively oriented contour enclosing the poles at $s_i=q^{-d_i} x_i$ for $i=1, \dots, k$, $d_i\in \mathbb{Z}_{\ge 0}$. We  also  used a shifted degree counting parameter $z_{\sharp}=(-1)^n  \hbar^{n/2}z$ and $\alpha_p$ is a normalization constant:
	$$
	\alpha_p=\prod\limits_{i,j=1}^{k} \dfrac{\varphi\Big(\frac{x_i}{x_j}\Big)}{\varphi\Big(\frac{q x_i}{\hbar x_j}\Big)}\,
	\prod\limits_{i=1}^{n}\prod\limits_{j=1}^{k} \dfrac{\varphi\Big(\frac{q x_j}{\hbar a_i}\Big)}{\varphi\Big(\frac{x_j}{a_i}\Big)} \prod\limits_{i=1}^{k}  e^{-\frac{\ln(z_{\sharp}) \ln(x_i) }{\ln(q)}}.
	$$

Such integral representation is convenient for the computation of the $q\rightarrow 1$ asymptotical behavior of the vertex function. 
In this limit a single term in the $q$-hypergeometric series dominates. 
Let integrand in (\ref{verint}) behaves in the following way:
 $$
e^{\frac{\phi(\{a_i\}, \{s_i\},z)}{\ln(q)}}(1+O(\ln(q)) 
 $$
The saddle point is defined by the equations:	
$s_i \partial_{s_i} \phi =0$ for $i=1,\dots,k$, or more explicitly:
\begin{equation}
\label{baeq}
\prod\limits_{j\neq i} \,\dfrac{s_i-s_j \hbar}{s_i \hbar-s_j} \prod\limits_{j=1}^{n} \, \dfrac{s_i-a_j}{a_j \hbar - s_i}=z \hbar^{-n/2}, \ \ i=1,\dots,k.
\end{equation}

They coincide with the Bethe equations up to  the substitution $\hbar\rightarrow \hbar^{-1}$, $\zeta^2\rightarrow (-1)^{n+1}z$. Thus we obtain the following Theorem.
\begin{Thm}\cite{Pushkar:2016qvw}\label{symmeig}
	The operators of quantum multiplication by quantum tautological classes
	$$
	{\mathcal{F}} \rightarrow \hat{\tau}(z) \circledast {\mathcal{F}}
	$$	
	are diagonal in the basis $\normalfont\psi_{\fp}(z)$. The corresponding eigenvalues are given
	by values of the symmetric polynomial $\tau(s_1,\cdots,s_k)$  at the solutions of the Bethe equations (\ref{baeq}), corresponding to $\fp$.
The algebra $\normalfont QK_{\mathsf T}(T^*{\rm Gr}_{k,n})$ coincides with Bethe algebra of $XXZ$ model, namely the algebra of symmetric functions of Bethe roots.
\end{Thm}

\subsection{Physics perspective}

In 2009 Nekrasov and Shatashvili \cite{Nekrasov:2009ui}, \cite{Nekrasov:2009uh} studied 3d supersymmetric gauge theories on $\mathcal{C}\times S^1$. These theories have gauge group ${G^c}=\times^{rank (\mathfrak{g})}_{i=1}{U}({\bf v}_i)$, the compact form of group $G$ considered in subsection \ref{qvar}. 

As usual in physics, the theory is defined by a certain action functional  $S({\phi_{\{\alpha\}}, A_{\{i\}}})$, where 

\begin{itemize}
\item $A_{\{i\}}$: ${U}({\bf v}_i)$-connections (gauge fields) on $\mathcal{C}\times S^1$;
\item $\phi_{\{\alpha\}}$: sections of associated vector ${U}({\bf v}_i)$-bundles on $\mathcal{C}\times S^1$, corresponding to the quiver data (matter fields).
\end{itemize}

Physicists compute the corresponding path integrals:
$$
\langle F\rangle=\int [d\phi_{\{\alpha\}}][dA_{\{i\}}] e^{- \beta S(\phi_{\{\alpha\}}, A_{\{i\}})}F({\phi_{\{\alpha\}}, A_{\{i\}}}).
$$

The supersymmetry properties of this integral allow to compute it through localization, relying on contribution from the minima of functional $S$, known  as {\it moduli of Higgs vacua}, which in Nekrasov-Shatashvili case are described by Nakajima quiver variety. There the momentum map relation $\mu=0$ iproduced a low energy configuration condition, while stability conditions are necessary for supersymmetry datum.

Due to supersymmetry properties the path integral above corresponds to the computation of index of a certain Dirac operator acting in certain infinite-dimensional space of field configurations:
$$\langle F\rangle={\rm str} (e^{-\beta \slashed{D}^2} F)={\rm tr}_{{\rm Ker} \slashed{D}_{even}}(F)-{\rm tr}_{{\rm Ker} \slashed{D}_{odd}}(F)={\rm str}_{{\rm index} \slashed{D}}(A).$$
The zero modes  ${\rm Ker} \slashed{D}$ corresponds to {\it Abrikosov vortices}, which we know as quasimaps.  The suitable observables-insertions $F$, producing Bethe algebra elements are the {\it line operators}, namely traces of holonomy of gauge fields. In our language the expression $\langle F\rangle$ corresponds to weighted K-theoretic counts of quasimaps, which are K-theoretic tensors with insertions corresponding to quantum tautological classes.

\subsection{Geometric interpretation of the $Q$-operator and $QQ$-systems}

Let us consider the $K$-theory class of $x$-weighted exterior algebra of 
$\mathcal{V}$:
$$
\Lambda_x = \bigoplus_{m=0}^{\infty}\, x^m \Lambda^{\!m}( \mathcal{V}) =\Lambda^{\! \bullet}_{x} \mathcal{V}.
$$

From previous subsection we obtain that the eigenvalues of quantum multiplication by exterior powers are the elemntary symmetric functions of Bethe roots.  The following Theorem identifies the above exterior powers with the expansion of the Baxter $Q$-operator in terms of spectral parameter.
\begin{Thm} \cite{Pushkar:2016qvw}\label{baxter}
The operator of quantum multiplication by the quantum tautological bundle $\hat{\Lambda}_x(z)$ coincides with the Baxter operator for  the $XXZ$ spin chain under the identification $K^{loc}_{\mathsf T}(\N(n)) = {\mathscr{H}}_{XXZ}$ as $\fsh$-modules: $
\hat{\Lambda}_x(z)\circledast\cdot ={ Q_+}(x).
$
\end{Thm}

At the same time, given that ${ Q_+}(x)$ is given as a trace over $\hbar$-oscillator representation $\rho_+(x)$ as we discussed in Section 2, it turns out that the following formula is true.

\begin{Thm}\cite{Pushkar:2016qvw}\label{Bform}
For arbitrary $n$ and $k$ the operator of quantum multiplication by the quantum $l$-th exterior power of tautological bundle
is given by the following universal formula:
$$
\widehat{\Lambda^{\! l} \mathcal{V}}(z)\circledast \cdot = \Lambda^{\! l} \mathcal{V} + a_1(z) \, F_0 \Lambda^{\! l-1} \mathcal{V} E_{-1}+ a_2(z) \, F_{0}^2 \Lambda^{\! l-2} \mathcal{V} E_{-1}^2 + \cdots
+ a_l(z) F_{0}^l  E_{-1}^l
$$
with
$$
a_m(z)= \frac{(\hbar-1)^m\ \hbar^{\frac{m(m+1)}{2}} K^{m}}{(m)_{\hbar}!\prod\limits_{i=1}^
{m}(1-(-1)^nz^{-1}\hbar^{i} K)}, \ \ \ \textrm{for}  \ \ \ (m)_{\hbar}=\dfrac{1-\hbar^m}{1-\hbar}, \ \ (m)_{\hbar}!=(1)_{\h}\cdots (m)_\hbar.
\nonumber
$$
Here $\widehat{\Lambda^{\! i} \mathcal{V}}(z)$  ($\Lambda^{\! i} \mathcal{V} $)  stands for the operators of quantum (classical) multiplication
by the quantum (classical) exterior powers and $E_r, F_r, H_r, K$ are the Drinfeld's generators of $\fsh$.
\end{Thm}

We note, that in the classical limit $z\rightarrow 0$ all $a_l(z)$ vanish and we obtain $\widehat{\Lambda^{\! l} \mathcal{V}}(0)= \Lambda^{\! l} \mathcal{V}$.

From the appropriate formulas for the operators of quantum multiplication by the quantum exterior powers $\widehat{\Lambda^{\! i}\mathcal{V}}(z)$ we obtain the {\it universal formulas} for any tautological classes:  that means the operators of quantum multiplication by $\hat{\tau}(z)$ are the elements in  $\fsh[[z]]$ which do not depend on $k,n$ defining the representation of $\fsh$ and thus the Nakajima variety.

Finally, let us discuss the geometric interpretation of the $QQ$-system in this case. So far we have seen interpretation of $Q_+(x)$. Now what about $Q_-(x)$? 
Let us consider the following short exact sequence of bundles:
\begin{eqnarray} 
0\to \mathcal{V}\to \mathcal{W}\to \hbar \otimes \mathcal{V}^{\vee }\to 0
\end{eqnarray}

The latter bundle can be interpreted as the tautological bundle of $ {\rm T}^*{\rm Gr}_{n-k,n}$ which is the same GIT quotient, but with the inverted stability parameter. The $Q_-(x)$-operator can be  identified with the generating function of the quantum exterior powers of $\mathcal{V}^{\vee}$, with K\"ahler parameter being inverted $z\rightarrow z^{-1}$.

\subsection{Further remarks}

\subsubsection{Generalizations to other Nakajima varieties} 
In this Section, we described the quantum K-theory of $ T^*{\rm Gr}_{k,n}$ as the Bethe algebra. This can be generalized to flag varieties and A-type Nakajima varieties in general \cite{Koroteev:2017aa},\cite{Koroteev:2023aa}. 
One can rewrite the capping operator $\Psi(z)$, which serves as a solution of qKZ equations using vertex functions with insertions \cite{Aganagic:2017be}. These insertions are described through K-theoretic stable envelopes. This, in particular, gives a geometric description for Bethe eigenfunctions. This principle can also be extended to, e.g., hypertoric varieties \cite{Smirnov:2020aa}. 

The universal formulas for quantum multiplication based on the combinatorial formulas as in Theorem \ref{Bform} are not known in the general case (even for Nakajima quiver varieties of type A) and could be an interesting project.

\subsubsection{3D Mirror symmetry}

{\it Symplectic duality} \cite{Kamnitzer:2022aa} is a duality operation on symplectic resolutions known in physics as 3D mirror symmetry \cite{Intriligator:1996ex}, \cite{Hanany:1996ie}, \cite{Aharony:1997bx}. While no complete definition exists, here is the enumerative geometry prescription for this definition. As we have seen, the vertex functions are the  $z$-solutions, namely analytic functions in K\"ahler variables in a punctured neighborhood of $\{z_i=0\}$, and have an infinite amount of poles in a punctured neighborhood of $\{a_i=0\}$. As we have seen, the vertex functions with insertions corresponding to stable envelopes are the solutions of difference equations in both $\{a\}$- and $\{z\}$-variables. The monodromy matrix relating $z$-solutions to $a$-solutions (analytic in equivariant variables) is known as an elliptic stable envelope introduced by Aganagic and Okounkov \cite{Aganagic:2016aa}. The $a$-solutions are conjectured to be the vertex functions for the 3d mirror dual variety, for which the  $\{a\}$ and $\{z\}$-parameters are interchanged as well as q-difference equations (\ref{qkzgeom}). The 3d mirror symmetry has been proven on the level of vertex functions in the example of ${\rm T}^*{\rm Gr}_{k,n}$ \cite{Dinkins:2020ab} and full flag varieties \cite{Rimanyi:2019ab}, \cite{Dinkins:2020ab} and hypertoric varieties \cite{Smirnov:2020aa}. On the level of quantum K-theory algebra, this 3d mirror symmetry has been studied in \cite{Koroteev:2023aa} for Nakajima varieties of type $A$ and instanton spaces using the relation to multiparticle systems.

\subsubsection{Quantum multiparticle systems and vertex functions}
Since the works of Cherednik \cite{Cherednik:1991mg} and Matsuo \cite{Matsuo1992}, the deep connection between solutions of the (deformed) Knizhnik-Zamolodchikov equations and quantum multiparticle systems has been investigated in various circumstances. A specific limit of that construction is known as quantum/classical duality \cite{Gorsky:2013xba}, \cite{Zabrodin:2017td}, \cite{Zabrodin:2017vt}, \cite{Zabrodin:}, \cite{KSZ} which we will discuss in Section 4.

In the first papers on quantum cohomology and various versions of quantum K-theory, these multiparticle systems played a vital role \cite{2001math8105G}. It was shown recently \cite{Koroteev:2018a} that the quasimap vertex function of partial flag varieties is the eigenfunction of the quantum Hamiltonians for the trigonometric Ruijsenaars-Schneider model. That, in particular, led to a nontrivial identity for stable envelopes and provided a geometric interpretation of Cherednik and Matsuo correspondence. At the same time, it was instrumental in the proof of 3D mirror symmetry on the level of vertex functions \cite{Dinkins:2020aa}, \cite{Dinkins:2020ab}.

\section{$(G,\hbar)$-opers and $QQ$-systems}

\subsection{$QQ$-systems and patterns in representation theory}

The $QQ$-systems we discussed in the Bethe ansatz context in Section 2 as deformed quantum Wronskian system have natural realization in the representation theoretic context.
\\

Here is a couple of situations how they emerge:

\begin{itemize}

 \item Let $\{V_{\omega_i}\}_{i=1,\dots, r}$ -- fundamental representations of simple Lie algebra $\mathfrak{g}$.

There are natural homomorphisms $m_i$:

$$
m_i:\quad \Lambda^2V_{\omega_i} \quad \to\quad \otimes_{j\neq i} V^{\otimes^{-a_{ji}}}_{\omega_j},
$$

where $v_i\wedge f_i v_i\to$  is highest weight of the above wedge product, and  $\{e_i, f_i, \check{\alpha}_i\}_{i=1,\dots, r}$ are the Chevalley generators of $\mathfrak{g}$.

This is how $QQ$-system appears in the celebrated ODE/IM correspondence \cite{Bazhanov:1996dr}, \cite{Dorey:2007zx}, \cite{MRV1}, \cite{MRV2}, \cite{Frenkel:2016}, where this formula plays a role in the spectral determinant decomposition for certain ODE problem corresponding to the so-called affine oper. 

The statement of the correspondence can be roughly formulated as follows. The vacuum eigenvalues of the Baxter operators quantum KdV model associated with affine Lie algebra $\widehat{\mathfrak{g}}$ appear as spectral determinants of certain singular differential operators associated with the so-called affine opers associated with ${^L}\widehat{\mathfrak{g}}$. In a particular case of standard quantum KdV, these operators are just singular Sturm-Liouville operators. As it was shown in \cite{MRV1,MRV2} they turn out to be the solution of the $QQ$-system with different analyticity conditions on entire $Q$-functions, which are just entire functions in this case. 

\vspace{3mm}

\item Given an $n\times n$ matrix $M$ the relations between minors is given by the following formula:

$$ \det(M^1_1)\det (M^k_k)-\det(M_1^k)\det(M^1_k)=\det (M) ~\det (M^{1,k}_{1,k}) $$

Here $M^i_j$ is a matrix with eliminated $i$-th row and $j$-th column. 
It has many names, e.g., Lewis Carrol identity or Pl\"ucker identity depending on the context.
When applied to the $\hbar$-deformed version of Wronskian matrix that leads to the $QQ$-system for $\mathfrak{g}=\mathfrak{sl}(n)$ \cite{MVflag}, \cite{2006math......4048M}, \cite{2012arXiv1210.2315L}, \cite{MR2189873}.

\vspace{3mm}

\item In this section we discuss a geometrization of $QQ$-system, which is related to the study of $\h$-defomed version of $G$-connections on the projective line ($G$-simple simply connected Lie group)  called 
$(G, \hbar)$-opers \cite{KSZ}, \cite{Frenkel:2020}. This relation is the generalization of the well-studied example of the geometric Langlands correspondence between $G$-local systems on projective line, represented by $G$-oper connections with regular singularties and $D$-modules on the moduli space of $^LG$-bundles with parabolic structures, explicitly described by $^L \mathfrak{g}$-Gaudin model. 

This correspondence encompasses the $\mathfrak{g}=\mathfrak{sl}(n)$ approach which we discussed above using $\hbar$-deformed Wronskians, which is related to alternative definition of $(G,\hbar)$-opers specific to $G=SL(n)$ case. At the same time, the so-called {\it generalized minors}: regular functions on $G$ satisfying relations similar to Lewis Carroll identities leading to alternative definition of $(G, \hbar)$-opers.

\end{itemize}

\subsection{$Z$-twisted Miura $(SL(2), \hbar)$-opers and the  $QQ$-system}
Fix $\hbar\in\C^*$.  For a given a vector bundle $E$ over $\P^1$, let $E^{\hbar}$
denote the pullback of $E$ under the map $$M_{\hbar}: \quad u\mapsto \hbar u.$$  We will always
assume that $E$ is trivializable.  Consider a
map of vector bundles $A:E\to E^{\hbar}$.  
Upon picking a trivialization, the map $A$ is determined by a matrix
$A(u)$ giving the linear map $E_u\to E_{\hbar u}$ in the given
bases.  A change in trivialization by $g(u)$ changes the
matrix via 
\begin{equation}\label{gauge tr}A(u)\mapsto g(\hbar u)A(u)g^{-1}(u);
\end{equation}
thus, $\hbar$-gauge change is twisted conjugation. Let $M_{\hbar}:E\to E^\hbar$ be
the operator that takes a section $s(u)$ to $s(\hbar u)$.  We associate the
map $A$ to the difference equation
$M_\hbar(s)=As$.  

A meromorphic {\it$(GL(r+1),\hbar)$-connection} over $\P^1$
  is a pair $(E,A)$, where $E$ is a (trivializable) vector bundle of
  rank $N$ over $\P^1$ and $A$ is a meromorphic section of the sheaf
  $\Hom_{\cO_{\P^1}}(E,E^\hbar)$ for which $A(u)$ is invertible.  The pair
  $(E,A)$ is called an $(SL(r+1),\hbar)$-connection if there exists a
  trivialization for which $A(u)$ has determinant $1$.\\

A {\it $(GL(2),\hbar)$-oper} on $\P^1$ is a triple
  $(E,A,\cL)$, where: \\

\begin{itemize}  
  
\item $E$ is a vector bundle on $\mathbb{P}^1$ of rank 2.\\
  
\item   $A$ is a $(GL(2),\hbar)$-connection \\

\item  $\cL$ is a line subbundle such that the induced map $\bar{A}:\cL\to (E/\cL)^\hbar$ is an isomorphism.  

\end{itemize}

The triple is called an \emph{$(\SL(2),\hbar)$-oper} if  $(E,A)$ is an $(\SL(2),\hbar)$-connection.

The condition that $\bar{A}$ is an isomorphism can be made explicit in
terms of sections.  Indeed, it is equivalent to $$s(\hbar u)\wedge A(u)
s(u)\neq 0$$ for $s(u)$ any section generating $\cL$ over either of
the standard affine coordinate charts.

However, one can violate this condition, by introducing the notion of a
$\hbar$-oper with regular singularities.  Let $\Lambda(u)$ be the polynomial.
A \emph{$(\SL(2),\hbar )$-oper with regular singularities} is a meromorphic
$(\SL(2),\hbar)$-oper $(E,A,\cL)$ for which $\bar{A}$ is an isomorphism  everywhere on except at the roots of $\Lambda(u)$.

This condition can be restated in terms of a section $s(u)$
generating $\cL$ over $\P^1$, namely 
\begin{eqnarray}
s(\hbar u)\wedge A(u) s(u)=\Lambda(u)
\end{eqnarray}

Next, we define twisted $\hbar$-opers; turns out that these are $\hbar$-analogues of the 
opers with a double pole singularity. Let $Z=\diag(\ze,\ze^{-1})$ be a diagonal
matrix with $\ze\ne\pm 1$. A $(\SL(2),\hbar)$-oper $(E,A,\cL)$ with regular singularities  is called a \emph{$Z$-twisted $\hbar$-oper} if $A$ is gauge-equivalent  to $Z$.

Finally, we will need the notion of a \emph{Miura $(SL(2), \hbar)$-oper}.  
As in the classical case, this is a quadruple $(E,A,\cL,\hcL)$ where
$(E,A,\cL)$ is a $SL(2,\hbar)$-oper and $\hcL$ is a line bundle preserved by $A$.

This way, for Miura $Z$-twisted $(SL(2),\hbar)$-oper, one can reformulate the condition as follows:
\begin{eqnarray}
s(\hbar u)\wedge Z s(u)=\Lambda(u), \quad s(u)=\begin{pmatrix} Q_-(u)\\Q_+(u)
  \end{pmatrix},
\end{eqnarray}
or, simply a familiar expression:
\begin{eqnarray}
\zeta Q_+(\hbar u)Q_-(u)-\zeta^{-1} Q_-(\hbar u)Q_+(u)=\Lambda(u)
\end{eqnarray}
Under the non-degeneracy conditions this conditions are equivalent to Bethe equations:
$$
\frac{\Lambda(w_i)}{\Lambda(\hbar^{-1}w_i)}=-\zeta^2\frac{Q_+(\hbar w_i)}{Q_+(\hbar^{-1} w_i)}, \quad Q_+(u)=\prod_{j}(u-w_j)
$$
Now, let us put the gauge connection of Miura oper in the canonical form in the basis of sections of $\mathcal{L}, \tilde{\mathcal{L}}$. 
First, we consider the gauge change by 
\begin{equation}
g(u)=\begin{pmatrix}
   Q_+(u) & -Q_-(u)\\
   0 &  Q_+^{-1}(u)
 \end{pmatrix},
\end{equation} which takes the section
$s(u)$ into $g(u)s(u)=\left(\begin{smallmatrix}0\\1
  \end{smallmatrix}\right)$.
  
In this gauge, the $\hbar$-connection matrix has the form
\begin{equation}\label{acon}
A(u)=U(\hbar u)ZU(u)^{-1}=
\begin{pmatrix} 
\zeta \frac{Q_+(\hbar u)}{Q_+(u)}& \Lambda(u)\\
   0 &  \zeta^{-1} \frac{Q_+(u)}{Q_+(\hbar u)}
 \end{pmatrix},
\end{equation}  

In terms of standard Chevalley generators of $\mathfrak{sl}(2)$ generators in 2-dimensional representations, the corresponding connection has the form:
\begin{eqnarray}
A(u)=g(u)^{\check{\alpha}} e^{\frac{\Lambda(u)}{g(u)}e}, \quad g(u)=\frac{Q_+(\hbar u)}{Q_+(u)},
\end{eqnarray}
where $e, f, \check{\alpha}$ is a standard Chevalley triple of $\mathfrak{sl}(2)$. 

At the end of this subsection, we mention that for $\zeta\neq \pm 1$, i.e., $Z$-regular semisimple, there are two Miura opers for a given $Z$-twisted Miura oper, corresponding to exchange $Q_+(u)\leftrightarrow Q_-(u)$, which is produced by the lift of the generator of the Weyl group element
in the gauge when $Z$ is diagonal.   
The transformation between two canonical forms of the canonical form of $\hbar$-connection (\ref{acon}) via lower-triangular matrix. We will discuss this in general in the next subsection.

\subsection{Miura $(G, \hbar)$-opers}

Let $G$ be the simple simply connected Lie group with Lie algebra $\mathfrak{g}$. Given a principal $G$-bundle $\cF_G$ over $\P^1$ (in Zariski
topology), let $\cF_G^\hbar$ denote its pullback under the map $M_\hbar$. 
A meromorphic $(G,\hbar)$-{\em
  connection} on a principal $G$-bundle $\cF_G$ on $\P^1$ is a section
$A$ of $\Hom_{\cO_{U}}(\cF_G,\cF_G^\hbar)$, where $U$ is a Zariski open
dense subset of $\P^1$. We can always choose $U$ so that the
restriction $\cF_G|_U$ of $\cF_G$ to $U$ is isomorphic to the trivial
$G$-bundle. Choosing such an isomorphism, i.e. a trivialization of
$\cF_G|_U$, we also obtain a trivialization of
$\cF_G|_{M_\hbar^{-1}(U)}$. Using these trivializations, the restriction
of $A$ to the Zariski open dense subset $U \cap M_\hbar^{-1}(U)$ can be
written as section of the trivial $G$-bundle on $U \cap M_\hbar^{-1}(U)$,
and hence as an element $A(z)$ of $G(z)$. Here we use the notation: if $K$ is a complex algebraic group, we set $K(z)=K(\C(z))$. Changing the trivialization of $\cF_G|_U$ via $g(u) \in G(u)$ changes
$A(u)$ by the $\hbar$-gauge transformation (\ref{gauge tr}). 
This shows that the set of equivalence classes of pairs $(\cF_G,A)$ as
above is in bijection with the quotient of $G(u)$ by the $\hbar$-gauge
transformations \eqref{gauge tr}.

Following the ideas from \cite{Frenkel1998,1998CMaPh.192..631S} we defined a {\it $(G,\hbar)$-oper} in \cite{Frenkel:2020} as a
$(G,\hbar)$-connection of a $G$-bundle on $\P^1$ equipped with a reduction
to the Borel subgroup $B_-$ that is not preserved by the
$(G,\hbar )$-connection but instead satisfies a special ``transversality
condition'' which is defined in terms of the Bruhat cell associated to
the Coxeter element $c$.

Let us state it in precise terms.  A meromorphic $(G,\hbar)$-{\em oper} on
  $\mathbb{P}^1$ is a triple $(\cF_G,A,\cF_{B_-})$:
  \begin{itemize}
\item   $G$-bundle $\cF_G$ on
  $\mathbb{P}^1$,
\item    $A$ is a
  meromorphic $(G,\hbar)$-connection on  $\cF_G$, 
   
 \item  $\mathcal{F}_{B_-}$ is the reduction of $\cF_G$  to $B_-$ 
\end{itemize}  
 satisfying the following condition: there exists a Zariski
  open dense subset $U \subset \P^1$ together with a trivialization
  $\imath_{B_-}$ of $\mathcal{F}_{B_-}$, such that the restriction of
  the connection $A: \cF_G \to \cF_G^\hbar$ to $U \cap M_\hbar^{-1}(U)$,
  written as an element of $G(z)$ using the trivializations of
  ${\mathcal F}_G$ and $\cF_G^\hbar$ on $U \cap M_\hbar^{-1}(U)$ induced by   $\imath_{B_-}$ takes values in the Bruhat cell $$B_-(\C[U \cap
  M_\hbar^{-1}(U)]) c B_-(\C[U \cap M_\hbar^{-1}(U)]).$$

Since $G$ is assumed to be simply-connected, any $\hbar$-oper connection
$A$ can be written using a particular trivialization $\imath_{B_-}$,
in the form
\begin{eqnarray}    \label{qop1}
A(z)=n'(z)\prod_i (\phi_i(z)^{\check{\alpha}_i} \, s_i )n(z)
\end{eqnarray}
where $\phi_i(z) \in\C(z)$ and  $n(z), n'(z)\in N_-(z)$ are such that
their zeros and poles are outside the subset $U \cap M_\hbar^{-1}(U)$ of
$\P^1$, and $s_i$ are the lifts of Weyl group $w_i$ tcorresponding to simple reflections to $G$.

Now let us define Miura $(G,\hbar)$-opers. \\

A {\em Miura $(G,\hbar)$-oper} on $\mathbb{P}^1$ is a quadruple:
  $(\cF_G,A,\cF_{B_-},\cF_{B_+})$, where:\\
   
\begin{itemize}  
\item   $(\cF_G,A,\cF_{B_-})$ is a
  meromorphic $(G,\hbar)$-oper on $\P^1$ \\  
  \item  $\cF_{B_+}$ is a reduction of
  the $G$-bundle $\cF_G$ to $B_+$ that is preserved by the
  $\hbar$-connection $A$. \\
\end{itemize}

Forgetting $\cF_{B_+}$, we associate a $(G,\hbar)$-oper to a given Miura
$(G,\hbar)$-oper. We will refer to it as the $(G,\hbar)$-oper underlying this
Miura $(G,\hbar)$-oper.

Now suppose we are given a principal $G$-bundle $\cF_G$ on any smooth
complex manifold $X$ equipped with reductions $\cF_{B_-}$ and
$\cF_{B_+}$ to $B_-$ and $B_+$, respectively. Then we assign to any
point $x \in X$ an element of the Weyl group $W_G$. Namely, the fiber
$\cF_{G,x}$ of $\cF_G$ at $x$ is a $G$-torsor with reductions
$\cF_{B_-,x}$ and $\cF_{B_+,x}$ to $B_-$ and $B_+$,
respectively. Choose any trivialization of $\cF_{G,x}$, i.e. an
isomorphism of $G$-torsors $\cF_{G,x} \simeq G$. Under this
isomorphism, $\cF_{B_-,x}$ gets identified with $aB_- \subset G$ and
$\cF_{B_+,x}$ with $bB_+$. Then $a^{-1}b$ is a well-defined element of
the double quotient $B_-\backslash G/B_+$, which is in bijection with
$W_G$. Hence we obtain a well-defined element of $W_G$.

We will say that $\cF_{B_-}$ and $\cF_{B_+}$ have a {\em generic
  relative position} at $x \in X$ if the element of $W_G$ assigned to
them at $x$ is equal to $1$ (this means that the corresponding element
$a^{-1}b$ belongs to the open dense Bruhat cell $B_- \cdot B_+ \subset
G$). The following theorem was proven in \cite{Frenkel:2020}.

\begin{Thm}    \label{gen rel pos}
  For any Miura $(G,\hbar)$-oper on $\mathbb{P}^1$, there exists an open
  dense subset $V \subset \P^1$ such that the reductions $\cF_{B_-}$
  and $\cF_{B_+}$ are in generic relative position for all $x \in V$.
\end{Thm}

 Let $U$ be a Zariski open dense subset on $\P^1$ as in the defininition of 
 $(G,\hbar)$-oper. Choosing a trivialization $\imath_{B_-}$ of $\cF_G$ on $U
  \cap M_\hbar^{-1}(U)$, we can write the $\hbar$-connection $A$ in the form
  \eqref{qop1}. On the other hand, using the $B_+$-reduction
  $\cF_{B_+}$, we can choose another trivialization of $\cF_G$ on $U
  \cap M_\hbar^{-1}(U)$ such that the $q$-connection $A$ acquires the form
  $\wt{A}(u) \in B_+(u)$. Hence there exists $g(u) \in G(u)$ such that
\begin{equation}    \label{connecting}
g(\hbar u) A(u)
g(u)^{-1} = \wt{A}(u) \in B_+(u).
\end{equation}
From the Bruhat decomposition $
G(u) = \bigsqcup_{w \in W_G} B_+(u) w N_-(u)$, we see that  
the statement of the proposition is equivalent to the statement that
$$
g(u) \in B_+(u) N_-(u),
$$
corresponding to $w=1$.

Then we have the following Theorem.

\begin{Thm}  \label{gen rel pos1}
i) For any Miura $(G,\hbar)$-oper on $\mathbb{P}^1$, there exists a
trivialization of the underlying $G$-bundle $\cF_G$ on an open
dense subset of $\P^1$ for which the oper $\hbar$-connection has the form
\begin{equation}    \label{genmiura}
A(z)\in N_-(z)\prod_i((\phi_i(z)^{\check{\alpha}_i}s_i
)N_-(z) \; \cap \; B_+(z).
\end{equation}
ii) Any element from the above intersection can be expressed in the following way:
\begin{equation}    \label{gicheck}
\prod_i g_i(z)^{\check{\alpha}_i}e^{\frac{t_i \phi_i(z)}{g_i(z)}e_i}, \qquad
g_i(u) \in \mathbb{C}(u),
\end{equation}
where $\{e_i, f_i, \check{\alpha}_i\}^r_{i=1}$ are the Chevalley generators, and $t_i$ depend on the lifting of $s_i$.
\end{Thm}
The first statement is a corollary of the previous Theorem. However, then  second requires non-trivial Lie-theoretic proof. In the following we set $t_i=1$.

\subsection{$Z$-twisted Miura opers with regular singularities}

 A $(G,\hbar)$-{\em oper with regular singularities determined by $\{
    \Lambda_i(u) \}_{i=1,\ldots,r}$} is a $(G,\hbar)$-oper on $\P^1$ whose
  $\hbar$-connection \eqref{qop1} may be written in the form
\begin{equation}    \label{Lambda}
A(u)= n'(u)\prod_i(\Lambda_i(u)^{\check{\alpha}_i} \, s_i)n(u), \qquad
n(u), n'(u)\in N_-(u).
\end{equation}

  {\em A Miura $(G,\hbar)$-oper with regular singularities determined by
polynomials $\{ \Lambda_i(u) \}_{i=1,\ldots,r}$} is a Miura
  $(G,\hbar)$-oper such that the underlying $(G,\hbar)$-oper has
regular singularities determined by $\{ \Lambda_i(u)
\}_{i=1,\ldots,r}$ and thus there exists trivialization where it is written in the form:
\begin{equation}    \label{form of A}
A(u)=\prod_i
g_i(u)^{\check{\alpha}_i} \; e^{\frac{\Lambda_i(u)}{g_i(u)}e_i}, \qquad
g_i(u) \in \C(u)^\times.
\end{equation}

Next, following the $SL(2)$ case, we consider a class of (Miura) $\hbar$-opers that are gauge equivalent to a constant element of $G$ as $(G,\hbar)$-connections.  Let $Z$ be an element of the maximal torus $H$. Since $G$ is
simply-connected, we can write:
\begin{equation}    \label{Z}
Z = \prod_{i=1}^r \zeta_i^{\check\alpha_i}, \qquad \zeta_i \in
\C^\times.
\end{equation}

 Then a {\em $Z$-twisted $(G,\hbar)$-oper} on $\mathbb{P}^1$ is a $(G,\hbar)$-oper   that is equivalent to the constant element $Z \in H \subset H(u)$  under the $\hbar$-gauge action of $G(u)$, i.e. if $A(u)$ is the
  meromorphic oper $\hbar$-connection (with respect to a particular
  trivialization of the underlying bundle), there exists $g(u) \in
  G(u)$ such that
\begin{eqnarray}    \label{Ag}
A(u)=g(\hbar u)Z g(u)^{-1}.
\end{eqnarray}
Following that,  {\em $Z$-twisted Miura $(G,\hbar)$-oper} is a Miura $(G,\hbar)$-oper on
$\mathbb{P}^1$ that is equivalent to the constant element $Z \in H
\subset H(u)$ under the $q$-gauge action of $B_+(u)$, i.e.
\begin{eqnarray}    \label{gaugeA}
A(u)=v(\hbar u)Z v(u)^{-1}, \qquad v(u) \in B_+(u).
\end{eqnarray}

We discuss in \cite{Frenkel:2020} the set of Miura  $(G,\hbar)$-opers corresponding to the given $Z$-twisted $(G,\hbar)$-oper. For example, when $Z$ is regular semisimple, then there is $|W|$ such Miura  $(G,\hbar)$-opers; on the other hand, when $Z=1$ there is an entire flag variety $G/B_-$ of such Miura  $(G,\hbar)$-opers. 

Notice that $Z$ twisted condition implies the restriction on diagonal part of Miura $(G,\hbar)$-oper. Namely, we know that for Miura $(G,\hbar)$-oper the underlying $(G,\hbar)$-connection can be written in the
form \eqref{form of A}. Since it preserves the $B_+$-bundle
$\cF_{B_+}$ underlying this Miura $(G,\hbar)$-oper, it may be viewed as a meromorphic $(B_+,\hbar)$-connection
on $\P^1$. Taking the quotient of $\cF_{B_+}$ by $N_+ = [B_+,B_+]$ and
using the fact that $B/N_+ \simeq H$, we obtain an $H$-bundle
$\cF_{B_+}/N_+$ and the corresponding $(H,\hbar)$-connection, which we
denote by $A^H(u)$. According to formula \eqref{form of A}, it is
given by the formula
\begin{equation}    \label{AH}
A^H(u)=\prod_ig_i(u)^{\check{\alpha}_i}.
\end{equation}
We call $A^H(u)$ the {\it associated Cartan $\hbar$--connection} of the
Miura $(G, \hbar)$-oper $A(u)$.
The $Z$-twisted condition immediately implies:
\begin{eqnarray}
g_i(u)=\zeta_i\frac{P_i(\hbar u)}{P_i(u)}, \quad {\rm where} \quad P_i(u)\in \mathbb{C}^{\times}.
\end{eqnarray}
Let us relax $Z$-twisted condition by considering a partial condition with the use of associated bundles.

\subsection{Z-twisted Miura-Pl\"ucker $(G,\hbar)$-opers and $QQ$-systems}
Now, let us relax the $Z$-twisted condition by considering a partial condition with the use of associated bundles.

Let $V_i$ be the fundamental representation of $G$ with highest 
weight $\omega_i$ with respect to $B_+$. It comes with a line $L_i
\subset V_i$ (of highest weight vectors) stable under the action of
$B_+$. Let
$\nu_{\omega_i}$ be a generator of the line $L_i \subset V_i$. It is a
vector of weight $\omega_i$ with respect to our Cartan subgroup $H
\subset B_+$. Then the subspace of $V_i$ of weight $\omega_i-\alpha_i$
is one-dimensional and is spanned by $f_i \cdot
\nu_{\omega_i}$. Therefore, the two-dimensional subspace $W_i$ of
$V_i$ spanned by the weight vectors $\nu_{\omega_i}$ and $f_i \cdot
\nu_{\omega_i}$ is a $B_+$-invariant subspace of $V_i$.

Let $(\cF_G,A,\cF_{B_-},\cF_{B_+})$ be a Miura $(G,\hbar)$-oper.  
Since $\cF_{B_+}$ is a $B_+$-reduction of a $G$-bundle $\cF_G$ on $\P^1$
preserved by the $(G,\hbar)$-connection $A$, for each
$i=1,\ldots,r$, the vector bundle
$$
\cV_i = \cF_{B_+} \underset{B_+}\times V_i = \cF_G \underset{G}\times
V_i
$$
associated to $V_i$ contains a rank two
subbundle
$$
\cW_i = \cF_{B_+} \underset{B_+}\times W_i
$$
associated to $W_i \subset V_i$, and $\cW_i$ in
turn contains a line subbundle
$$
\cL_i = \cF_{B_+} \underset{B_+}\times L_i
$$
associated to $L_i \subset W_i$.

If $\phi_i(A)$ denotes the $(GL(V_i),\hbar)$-connection on the vector bundle $\cV_i$ corresponding to the 
above Miura $(G,\hbar)$-oper connection $A$, then since $A$ preserves $\cF_{B_+}$, we obtain that $\phi_i(A)$ preserves the
subbundles $\cL_i$ and $\cW_i$ of $\cV_i$. We further denote by $A_i$ the  
corresponding $\hbar$-connection on the rank 2 bundle $\cW_i$. 
Let us trivialize $\cF_{B_+}$ on a Zariski open subset of $\P^1$ so
that $A(u)$ has the form \eqref{form of A} with respect to this
trivialization. This trivializes the
bundles $\cV_i$, $\cW_i$, and $\cL_i$, so that the $\hbar$-connection
$A_i(u)$ becomes a $2 \times 2$ matrix whose entries are in $\C(u)$:
\begin{equation}    \label{2flagformula}
A_i(u)=\begin{pmatrix}
  g_i(u) &  &\phi_i(z) \prod_{j>i} g_j(u)^{-a_{ji}}\\
&&\\  
  0 & &g^{-1}_i(u) \prod_{j\neq i} g_j(u)^{-a_{ji}}
 \end{pmatrix},
\end{equation}
where we use the ordering of the simple roots determined by the
Coxeter element $c$.

Using the trivialization of $\cW_i$ in which $A_i(z)$ has the form
\eqref{2flagformula}, we represent $\cW_i$ as the direct sum of two
line subbundles: the first is $\cL_i$, generated by the basis vector
$\begin{pmatrix} 1 \\ 0 \end{pmatrix}$; the second, which we denote by
$\wt\cL_i$, is generated by the basis vector $\begin{pmatrix} 0 \\
  1 \end{pmatrix}$. The subbundle $\cL_i$ is $A_i$-invariant, whereas
the subbundle $\wt\cL_i$ and $A_i$ satisfy the 
Miura $(GL(2),\hbar)$-oper condition.
Note, that if $(G,\hbar)$-oper is $Z$-twisted, then $A_i(u)$ is $Z_i$-twisted, where  $Z_i = Z|_{W_i}$. 

Turning this around,  we can formulate the folowing definition, which relaxes the original $Z$-twisted condition.

  A $Z$-{\em twisted Miura-Pl\"ucker $(G,\hbar)$-oper} is a meromorphic
  Miura $(G,\hbar)$-oper on $\P^1$ with the underlying $\hbar$-connection
  $A(z)$, such that there exists $v(z) \in B_+(z)$ such that for all
  $i=1,\ldots,r$, the Miura $(GL(2),\hbar)$-opers $A_i(u)$ associated to
  $A(u)$ by formula \eqref{2flagformula} are $Z_i$-twisted, i.e. can be written in the form:
\begin{equation}    \label{gaugeA3}
A_i(u) = v(\hbar u) Z v(u)^{-1}|_{W_i} = v_i(\hbar u)Z_iv_i(u)^{-1}
\end{equation}
where $v_i(z) = v(z)|_{W_i}$ and $Z_i = Z|_{W_i}$.

Assuming this condition, the element $v(z)$ can be expressed in the form 
\begin{equation}    \label{vdots}
v(u) = \prod_{i=1}^r P_i(u)^{\check\alpha_i} \prod_{i=1}^r
e^{-\frac{Q^j_{-}(u)}{Q^j_{+}(u)} e_i} \dots ,
\end{equation}
where the dots stand for the exponentials of higher commuatator terms
in the Lie algebra ${\mathfrak n}_+$ of $N_+$ (these terms will not
matter in the computations below) and $\{Q^i_+(u),Q^i_-(u)\}$ are
relatively prime polynomials and $Q^i_+(u)$ is a monic polynomial for
each $i=1,\ldots,r$.
Acting on the two-dimensional subspace $W_i$, $v(z)$ has the following form:
\begin{eqnarray}    \label{vzz}
v(u)|_{W^i}=
\begin{pmatrix}
  P_i(u) & 0\\
  0& P_i^{-1}(u)\prod_{j\neq i} P_j^{-a_{ji}}(u)
 \end{pmatrix}
 \begin{pmatrix}
1 & - \frac{Q^i_{-}(u)}{Q^i_{+}(u)}\\
 0& 1
 \end{pmatrix}
\end{eqnarray}
while $Z$ has the form
\begin{eqnarray}
Z|_{W_i}=\begin{pmatrix}
\zeta_i & 0\\
  0& \zeta_i^{-1}\prod_{j\neq i} \zeta_j^{-a_{ji}}
   \end{pmatrix}.
\end{eqnarray}

From formulas \eqref{2flagformula}  we obtain the
following equations on $P_i(z)$ and $Q^i_{\pm}(z)$. First, comparing
the diagonal entries on both sides of $A_i(u) = v_i(\hbar u)Z_iv_i(u)^{-1}$, we obtain:
\begin{eqnarray}    \label{giz}
g_i(u)=\zeta_i\frac{P_i(\hbar u)}{P_i(u)},
\end{eqnarray}
which we have seen already for $Z$-twisted $(H, \hbar)$-connections. 
Second, by comparing the upper triangular entires on both sides of
$A_i(u) = v_i(\hbar u)Z_iv_i(u)^{-1}$, we obtain:
\begin{multline}    \label{Lambdai}
\phi_i(u)\prod_{j>i}\Big[\zeta_j\frac{P_j(\hbar u)}{P_j(u)}\Big]^{-a_{ji}} \; = \\ P_i(u)P_i(\hbar u)\prod_{j\neq
  i}P_j(u)^{a_{ji}}\left[
\zeta_i\frac{Q^i_{-}(u)}{Q^i_{+}(u)}-\zeta_i^{-1}\prod_{j\neq i} \zeta_j^{-a_{ji}}
\frac{Q^i_{-}(\hbar u)}{Q^i_{+}(\hbar u)}\right].
\end{multline}

We call this system of equations a $PQ$-system. 
Let us bring in the following {\it nondegeneracy conditions} on $P_i(u)$:\\

\begin{itemize}
\item $P_i(u)$ are polynomial,\\

\item $P_i(u)$ and $P_j(u)$ have $\hbar$-distinct roots if $i,j$ correspond to neighboring nodes on a Dynkin diagram. \\
\end{itemize}

We call such $Z$-twisted Miura-Pl\"ucker oper {\it nondegenerate}. 
If the $(G, \hbar)$-oper has regular singularities, i.e. $\phi_i(u)=\Lambda_i(u)\in \mathbb{C}^{\times}$, $i=1, \dots, r$, we obtain that 
\begin{eqnarray}
P_i(u)=Q_+^i(u)
\end{eqnarray}
and the $PQ$-system transforms into $QQ$-system:
\begin{multline}\label{qq}
\wt{\xi}_iQ^i_{-}(u)Q^i_{+}(\hbar u)-\xi_iQ^i_{-}(\hbar u)Q^i_{+}(u) = \\
\Lambda_i(u)\prod_{j> i}\Big[Q^j_{+}(\hbar u)\Big]^{-a_{ji}}
\prod_{j< i}\Big[Q^j_{+}(u)\Big]^{-a_{ji}}, \qquad
i=1,\ldots,r,
\end{multline}
where
\begin{equation}    \label{xi}
\wt{\xi}_i=\zeta_i \prod_{j>i} \zeta_j^{a_{ji}}, \qquad
{\xi}_i=\zeta^{-1}_i\prod_{j< i} \zeta_j^{-a_{ji}}
\end{equation}
and we use the ordering of simple roots from the definition of the
$(G,\hbar)$-opers.
The resulting connection, corresponding to nondegenerate Miura-Pl\"ucker $(G, \hbar)$-oper has the following form:
\begin{equation}    \label{key}
A(u) =\prod_j\left[ \zeta_j\frac{Q^j_+(\hbar u)}{Q^j_+(u)}
\right]^{\check{\alpha}_j} e^{\frac{\Lambda_j(u) Q^j_+(u)}{\zeta_j
    Q^j_+(\hbar u)}e_j} = \prod_j
\Big[\zeta_jQ_{+}^j(\hbar u )\Big]^{\check{\alpha}_j}
e^{\frac{\Lambda_j(u)}{\zeta_j Q_+^j(\hbar u)Q_+^j(u)}e_j}
\Big[{Q_{+}^j(u)}\Big]^{-\check{\alpha}_j}\,.
\end{equation}

Thus we have a one-to-one correspondence between $Z$-twisted Miura-Pl\"ucker $(G,\hbar)$-opers and nondegenerate polynomial solutions of the $QQ$-system \eqref{qq}.
There is a stronger version of this statement, extending it from nondegenerate $Z$-twisted Miura-Pl\"ucker $(G,\hbar)$-opers to $Z$-twisted Miura $(G,\hbar)$-opers. 

To do that one has to address the structure of $QQ$-system as eigenvalues of Baxter $Q$-operators and the limit $z_i=\prod_{i=1} \zeta_j^{a_{ij}}\rightarrow 0$. An important feature of the proof are the so-called quantum B\"acklund transformations, which were studied in \cite{Mukhin_2005} in the case of additive parameter shift $u\to u+\hbar$. For $Z$-twisted Miura $(G,\hbar)$-opers these transformations allow to travel between various Miura opers corresponding to 
a given $(G, \hbar)$-opers. Let us have a look at how it works explicitly.

\begin{Prop}    \label{fiter}
  Consider the $\hbar$-gauge transformation of the $\hbar $-connection $A$
  given by formula \eqref{key}:
\begin{equation}
A(u) \mapsto A^{(i)}(u)=e^{\mu_i(\hbar u)f_i}A(u)e^{-\mu_i(z)f_i},
\quad \operatorname{where} \quad \mu_i(z)=\frac{\prod\limits_{j\neq
    i}\Big[Q_+^j(u)\Big]^{-a_{ji}}}{Q^i_{+}(u)Q^i_{-}(u)}\,.
\label{eq:PropDef}
\end{equation}
Then $A^{(i)}(u)$ can be obtained from $A(u)$ by
substituting 
\begin{align}
Q^j_+(u) &\mapsto Q^j_+(u), \qquad j \neq i, \nonumber\\
Q^i_+(u) &\mapsto Q^i_-(u), \qquad Z\mapsto s_i(Z)\,.\nonumber
\end{align}
in the formula \eqref{key}.
\end{Prop}

Producing all such transformations gives rise to all $Z$-twisted Miura $(G,\hbar)$-opers for given $(G,\hbar)$-oper. In \cite{Frenkel:2020} using B\"acklund transformation, we produce an element $v(u)\in B_+(u)$, diagonalizing $A(u)$, i.e.  $A(u)=v(\hbar u)Zv(u)^{-1}$, thus producing the following result:
\begin{Thm}    
  There is a one-to-one correspondence between the set of
  nondegenerate $Z$-twisted Miura $(G,\hbar)$-opers and the set
  of nondegenerate polynomial solutions of the $QQ$-system \eqref{qq}.
\end{Thm}

\begin{Rem}
We note that the $QQ$-system we obtained in (\ref{qq}) corresponds to Bethe ansatz equations of $\mathcal{U}(^L\widehat{\mathfrak{g}})$ and the relations within the Grothendieck ring of category $\mathcal{O}$ of Borel subalgebra only in ADE case \cite{Frenkel:2016}. At the same time, recently the integrable models corresponding to such $QQ$-systems were proposed in \cite{FHRnew}.
\end{Rem}

\begin{Rem}
We call the collection of $QQ$-systems which correspond to all Miura $(G,\hbar)$-opers for a given $Z$-twisted $(G,\hbar)$-oper the {\it full} $QQ$-system \cite{KoroteevZeitlinCrelle}. Recently it was conjectured that this has a representation-theoretic meaning within category $\mathcal{O}$ of  
the Borel subalgebra of $\mathcal{U}(\hat{\mathfrak{g}})$, i.e. there are elements in this category, labled by element s of Weyl  group satisfying the corresponding relations.
\end{Rem}

\subsection{Quantum/calssical duality and quantum K-theory of flag varieties}

The generalization of example of $(SL(2), \hbar)$-oper considered in the beginning of this section is as follows.

  A meromorphic $(GL(r+1),\hbar)$-{\em oper}  on
  $\mathbb{P}^1$ is a triple $(\mathcal{A},E, \mathcal{L}_{\bullet})$, where $E$ is  a vector bundle of rank $r+1$ and $\mathcal{L}_{\bullet}$ is the corresponding complete flag of the vector bundles, 
  $$\mathcal{L}_{r+1}\subset ...\subset \mathcal{L}_{i+1}\subset\mathcal{L}_i\subset\mathcal{L}_{i-1}\subset...\subset \mathcal{L}_1=E,$$ 
  where $\mathcal{L}_{r+1}$ is a line bundle, so that $\hbar$-connection 
  $\mathcal{A}\in \Hom_{\cO_{U}}(E,E^{\hbar})$ 
  satisfies the following conditions:\\
  
  \begin{itemize} 
\item $\mathcal{A}\cdot \mathcal{L}_i\subset \mathcal{L}_{i-1} $,\\
\item There exists a Zariski open dense subset $U \subset \P^1$, such that the restriction of $\mathcal{A}\in Hom(\mathcal{L}_{\bullet}, \mathcal{L}^\hbar_{\bullet})$ to $U \cap M_{\hbar}^{-1}(U)$ is invertible and satisfies the condition that the induced maps  $$\bar{\mathcal{A}}_i:\mathcal{L}_{i}/\mathcal{L}_{i+1}\to \mathcal{L}_{i-1}/\mathcal{L}_{i}$$ 
are isomorphisms on $U \cap M_\hbar^{-1}(U)$. \\
\end{itemize}

  An $(SL(r+1),\hbar)$-$oper$ is a $(GL(r+1),\hbar)$-oper with the condition that $det(\mathcal{A})=1$ on $U \cap M_\hbar^{-1}(U)$.\\

The {\it Miura $(SL(r+1),\hbar)$-oper} is the following qudruple $(\mathcal{A},E, \mathcal{L}_{\bullet}, \hat{\mathcal{L}}_{\bullet})$:
\begin{itemize} 
\item  $(\mathcal{A},E, \mathcal{L}_{\bullet})$ is $(SL(r+1),\hbar)$-{ oper},

\item $\hat{\mathcal{L}}_{\bullet}$ is a full flag of subbundles preserved by $\mathcal{A}$.
\end{itemize}

The regular singularities condition in this case means that $\bar{\mathcal{A}}_i(u)$ is not invertible at zeroes of polynomials $\Lambda_i(u)$.

The equivalence of the general definition can be proven along the same lines as the equivalence of the analogous definitions in the case of classical opers. 
One can derive the second definition from the first by considering the associated bundle $E=(\cF_{SL(r+1)}\times V_{\omega_1})/SL(r+1)$, where $V_{\omega_1}$ in the defining representation of $G$. That immediately provides a flag of subbundles in $E$, preserved by $B_-$. From the chosen order in the Coxeter element we obtain that the induced   $\hbar$-connection on $E$ locally has the form of the matrix with coefficients in $\mathbb{C}(u)$ so that it has zeroes above the superdiagonal. That immediately leads to two  conditions of the definition above. 

Explicitly, considering the determinants  
\begin{equation}\label{altqW} 
 \left(\Big(\prod_{j=0}^{i-2}(A(\hbar^{i-2-j}u)\Big)s(u)\wedge\dots\wedge A(\hbar^{i-2}u) s(\hbar^{i-2}u)\wedge s(q^{i-1}u)
  \right)\bigg|_{\Lambda^i\cL_{r-i+2}^{\hbar^{i-1}}}
\end{equation}
for $i=1,\dots, r+1$,  
where $s$ is a local section of $\cL_{r+1}$, we claim that 
$(E,A,\cL_\bullet)$ is an $(SL(r+1),\hbar)$-oper if and only if at every point of $U \cap M_\hbar^{-1}(U)$, there 
exists local section for which each  such determinant  is nonzero (see \cite{KSZ}).   
When we encounter the case of regular singularities, each  
$\bar{A}_i$ is an isomorphism except at zeroes of $\Lambda_i(u)$ and thus we require the determinants to vanish at zeroes of the following polynomial $W_k(s)$: 
\begin{eqnarray}
W_k(s)=P_1(u) \cdot P_2(\hbar^2u)\cdots P_{k}(\hbar^{k-1}u),  \qquad 
P_i(z)=\Lambda_{r}(u)\Lambda_{r-1}(u)\cdots\Lambda_{r-i+1}(u)\,.\nonumber
\end{eqnarray}  

Now we discuss the $Z$-twisted Miura condition. Miura  condition implies that there exist a flag $\hcL_\bullet$ which is preserved by the $\hbar$-connection $A$. The $Z$-twisted condition implies that in the gauge when $A$ is given by fixed semisimple diagonal element $Z\in H$ such flag is formed by the standard basis $e_1, \dots, e_{r+1}$. 

The relative position between two flags is generic on $U \cap M_\hbar^{-1}(U)$. The regular singularity condition implies that {\it quantum Wronskians}, namely  determinants
\begin{equation}\label{qD}
\mathcal{D}_k(s)(u)=e_1\wedge\dots\wedge{e_{r+1-k}}\wedge
Z^{k-1}s(u)\wedge Z^{k-2} s(\hbar u)\wedge\dots\wedge Z s(\hbar ^{k-2}u)\wedge s(\hbar^{k-1}u)\,
\end{equation}
have a subset of zeroes, which coincide with those of
$\cW_k(s)$.  To be more explicit, for
$k=1,\dots,r+1$, we have nonzero constants $\alpha_k$ and polynomials
\begin{equation} \mathcal{V}_k(u) = \prod_{a=1}^{r_k}(u-v_{k,a})\,,
\label{eq:BaxterRho}
\end{equation}
for which 
\begin{equation}
\det\begin{pmatrix} \,     1 & \dots & 0 & \xi_1^{k-1}s_{1}(u) & \cdots & \xi_{1} s_{1}(\hbar^{k-2}u)  &  s_{1}(\hbar^{k-1}u) \\ 
 \vdots & \ddots & \vdots& \vdots & \vdots & \ddots & \vdots \\  
0 & \dots & 1&\xi_{k}^{k-1}s_{r+1-k}(z) &\dots & \xi_{k} s_{r+1-k}(\hbar^{k-2}u) &   s_{k}(\hbar^{k-1}u)  \\  
0 & \dots & 0&\xi^{k-1}_{k+1}s_{r+1-k+1}(u) & \dots & \xi_{r+1-k+1} s_{k+1}(\hbar^{k-2}u)  &  s_{k+1}(\hbar^{k-1}u)  \\
\vdots & \ddots & \vdots&\vdots & \vdots & \ddots & \vdots \\
0 & \dots & 0&\xi_{r+1}^{k-1}s_{r+1}(u) & \dots &\xi_{r+1} s_{r+1}(\hbar^{k-2}u)  & s_{r+1}(\hbar^{k-1}u)  \, \end{pmatrix} =\alpha_{k} W_{k}
\cV_{k} \,,
\label{eq:MiuraQOperCond}
\end{equation}

\vspace*{2mm}
\noindent where $Z=\diag(\xi_1, \dots, \xi_r)$ so that  $\xi_i=\zeta_i/\zeta_{i-1}$ and $\alpha_k$ are constants. 
Since $\cD_{r+1}(s)=W_{r+1}(s)$, we have
$\cV_{r+1}=1$.  In addition, $\cV_0=1$ and 
$\cD_0=e_1\wedge\dots\wedge e_{r+1}$.

One can identify $\mathcal{V}_k(u)$ with $Q_+^k(u)$ for the $QQ$-system corresponding to $\mathfrak{g}=\mathfrak{sl}(r+1)$. We note here that the $QQ$-relation can be obtained as the Lewis Carroll identity for the minors (\ref{eq:MiuraQOperCond}) of $\hbar$-deformed Wronskian matrix.

Now let us consider the case when only $\Lambda_{1}(u)=\prod^{r+1}_{i=1}(u-a_i)$ and $s_i(u)=c_i(u-p_i)$ are polynomials of degree 1.

The corresponding Bethe equations correspond to the weight zero subset in $\mathfrak{sl}_{r+1}$ XXZ spin chain, which involves only defining fundamental representations with evaluation parameters $\{a_i\}_{i=1, \dots, r}$:
\begin{eqnarray}
\mathscr{H}=\mathbb{C}^{r+1}(a_1)\otimes \dots \otimes \mathbb{C}^{r+1}(a_{r+1})
\end{eqnarray}
From a point of view of Nakajima varieties this subspace in $\mathscr{H}$ corresponds to the localized K-theory of cotangent bundle of full flags: $K^{loc}_{T}(T^*{\rm Fl}_{r+1})$ with $T=\mathbb{C}^{\times}_{a_1}\times\dots \times \mathbb{C}^{\times}_{a_{r+1}}\times \mathbb{C}^{\times}_{\hbar}$. The corresponding $QQ$-system defines relations in the quantum K-theory of  $T^*Fl_{r+1}$.
We know that the Wronskian relation 
$$\mathcal{D}_{r+1}(u)=\Lambda_1(u)$$ 
is equivalent to the $QQ$-system. 
On the other hand, that leads to  the following condition:
\begin{equation} \label{Lag}
 \{H_k=f_k(\{a_i\})\},\quad k=1, \dots, r+1,
\end{equation}
where 
\begin{eqnarray}\label{trs}
H_k=\sum_{\substack{J\subset \{1,\dots, r+1\}\\|J|=k}}
\prod_{\substack{i\in J\\j \notin J}}\frac{{\xi_i}-{\hbar}{\xi_j}}{{\xi_i}-{\xi_j}}\prod_{m\in J}p_m
\end{eqnarray}
are known as trigonometric Ruijsenaaars-Schneider (tRS) Hamiltonians 
and $f_k$ are elementary symmetric functions of ${a_i}$.
The conditions (\ref{Lag}) correspond to the the Lagrangian subvariety in the 
tRS phase space with symplectic form:
\begin{eqnarray}
\Omega=\sum^{r+1}_{i=1}\frac{d\xi_i}{\xi_i}\wedge \frac{dp_i}{p_i}
\end{eqnarray}
Thus we have the following Theorem.
\begin{Thm}
The quantum equivariant K-theory algebra $K_{T}(T^*{\rm Fl}_{r+1})$ can be expressed as an algebra of functions on the space of $(SL(r+1), \hbar)$-opers with regular singularities described in two equivalent ways:
\begin{itemize}
\item  Symmetric functions of Bethe roots/solutions of the related $QQ$-system,\\

\item  Algebra of functions on Lagrangian subvariety in tRS phase space given by tRS Hamiltonians (\ref{trs}).
 \end{itemize}
\end{Thm}

\subsection{Further remarks}
\subsubsection{Generalized minors and $(G, \hbar)$-opers}
As we have seen, the $Z$-twisted $SL(n)$-opers can be described explicitly using $\hbar$-difference Wronskians. That approach is generalizable for any simple, simply connected $G$ using the so-called generalized minors \cite{FZ}. These are regular functions $\{\Delta_{u\omega_i, v\omega_i}\}^{i=1, \dots, r}_{u,v\in W}$on simple simply connected group $G$ which are labeled by fundamental weights $\omega_i$ and the action of Weyl group elements of $u, v$.  One of the relations they satisfy is the generalization of Lewis Carroll identity, namely: 
 \begin{eqnarray}
\Delta_{u\cdot\omega_i, v\cdot\omega_i}(g)\Delta_{uw_i\cdot \omega_i, vw_i\cdot\omega_i}(g)-
\Delta_{uw_i\cdot\omega_i, v\cdot\omega_i}(g)\Delta_{u\cdot\omega_i, vw_i\cdot\omega_i}(g)=\prod_{j\neq i}\Big[\Delta_{u\cdot\omega_j, v\cdot\omega_j}(g)\Big]^{-a_{ji}}, \nonumber
\end{eqnarray}
where the lengths $\ell$ of Weyl group elements satisfy $\ell(uw_i)=\ell(u)+1$, $\ell(vw_i)=\ell(v)+1$ and $g\in G$. 

The generalized $\hbar$-Wronskian matrix is an element $\mathscr{G}(u)\in G(u)$, which is the meromorphic section of $G$-bundle on $\mathbb{P}^1$ satisfies the following difference equations:
 \begin{eqnarray}
Z^{-1}\mathscr{G}(\hbar u)v_{\omega_i}=
\mathscr{G}(u)\prod^r_{i=1}
\Lambda_i^{\check{\alpha}_i}(u)s_i
v_{\omega_i}, \quad i=1,\dots, r.
\end{eqnarray}
The relations for the generalized minors of $\mathscr{G}(u)$ give the equations of the $QQ$-system. We refer the reader to \cite{KoroteevZeitlinCrelle} (see also \cite{zeitlin2022wronskians}) for further details and to recent preprint \cite{geiss} for the explicit relation to cluster algebras.

\subsubsection{The Zoo of opers and dualities} 
The quantum/classical duality between quantum XXZ  and classical tRS models we considered in this Section can be generalized to other types of spin chains: XXX, trigonometric and rational Gaudin,  and multiparticle systems: trigonometric and rational Calogero-Moser (tCM, rCM), rational Ruijsennars-Schneider (rRS). That can be reformulated again using a change of coordinates for various versions of (deformed) oper connections \cite{koroteev2022zoo}. 

Besides quantum/classical duality, one more duality is intrinsically embedded in the explicit construction due to the article \cite{Oblomkov2004} of multiparticle models (all being limits of the tRS model). This has the relation to 3D mirror symmetry we discussed briefly in the remarks to Section 3: namely, one can exchange the variables $\{a_i\}$ and $\{\xi_i\}$ in tRS description of quantum equivariant K-theory ring of the flag variety described as algebra of functions on a Lagrangian subvariety in tRS phase space. While the tRS and rCM models are self-dual, under this duality, rRS and tCM are dual to each other. We will refer to the paper \cite{Koroteev:2023aa}, \cite{koroteev_zeitlin_2023} for further applications of this duality to the 3D mirror symmetry of various Nakajima varieties. 

\subsubsection{$(G,\hbar)$-opers and quantum $q$-Langlands duality}

We already mentioned in the remarks to Section 2 that there are two types: $a$- and $z$- solutions to qKZ equations (\ref{qkzgeom}), which are related to each other via elliptic stable envelope. It turns out that $a$-solutions have a natural meaning as the matrix elements of the products of intertwining operators for $\mathcal{U}_{\hbar}(\mathfrak{g})$ with nontrivial central charge: the conformal blocks. The $z$-solutions, i.e., vertex functions, have a meaning of the conformal blocks for deformed W-algebras $W_{q,t}(^L\mathfrak{g})$. The correspondence between these sets of conformal blocks gives what is called quantum $q$-Langlands correspondence \cite{Aganagic:2017smx} (see also recent physics generalization \cite{haouzi2023new}). The particular limit $q\rightarrow 1$ corresponds to the critical case, which in the differential limit $\hbar\rightarrow 1$ corresponds to relations between oper connections: the differential operators governing $W(^L\mathfrak{g})$-algebra conformal blocks and spectrum of Gaudin model. In the deformed case, the $\hbar$-difference equations governing conformal blocks of $W_{q,t}(^L\mathfrak{g})$ should be related to $(G,\hbar)$-opers as we know that the $qKZ$ equation gives eigenvalue problem for $XXZ$-model, which is a deformation of Gaudin model. This is part of the bigger story, the deformation of Langlands duality, which is not yet entirely shaped.

\bibliography{biblio_lect}

\begin{bibdiv}
\begin{biblist}

\bib{arnold78}{book}{
      author={Arnold, V.~I.},
       title={{Mathematical Methods of Classical Mechanics}},
   publisher={Springer},
        date={1978},
}

\bib{Aganagic:2017smx}{article}{
      author={Aganagic, M.},
      author={Frenkel, E.},
      author={Okounkov, A.},
       title={{Quantum $q$-Langlands Correspondence}},
        date={2018},
     journal={Trans. Moscow Math. Soc.},
      volume={79},
       pages={1\ndash 83},
      eprint={1701.03146},
}

\bib{Aharony:1997bx}{article}{
      author={Aharony, Ofer},
      author={Hanany, Amihay},
      author={Intriligator, Kenneth~A.},
      author={Seiberg, N.},
      author={Strassler, M.J.},
       title={{Aspects of N=2 supersymmetric gauge theories in
  three-dimensions}},
        date={1997},
     journal={Nucl.Phys.},
      volume={B499},
       pages={67\ndash 99},
      eprint={hep-th/9703110},
}

\bib{Aganagic:2017be}{article}{
      author={Aganagic, Mina},
      author={Okounkov, Andrei},
       title={{Quasimap counts and Bethe eigenfunctions}},
        date={2017},
     journal={Moscow Math. J.},
      volume={17},
      number={4},
       pages={565\ndash 600},
      eprint={1704.08746},
}

\bib{Aganagic:2016aa}{article}{
      author={Aganagic, Mina},
      author={Okounkov, Andrei},
       title={Elliptic stable envelopes},
        date={2021},
     journal={J. Amer. Math. Soc.},
      volume={34},
       pages={79\ndash 133},
      eprint={1604.00423},
         url={https://arxiv.org/pdf/1604.00423.pdf},
}

\bib{Baxter:1982zz}{book}{
      author={Baxter, R.~J.},
       title={{Exactly solved models in statistical mechanics}},
   publisher={Academic Press},
        date={1982},
        ISBN={0486462714, 9780486462714},
         url={http://www.amazon.com/dp/0486462714},
}

\bib{babelon}{book}{
      author={Babelon, O.},
      author={Bernard, D.},
      author={Talon, M.},
       title={{Introduction to Classical Integrable Systems}},
   publisher={CUP},
        date={2003},
}

\bib{BFLMS}{article}{
      author={Bazhanov, V.},
      author={Frassek, R.},
      author={Lukowski, T.},
      author={Meneghelli, C.},
      author={Staudacher, M.},
       title={{Baxter Q-Operators and Representations of Yangians}},
        date={2011},
     journal={Nucl. Phys.},
      volume={B850},
       pages={148\ndash 174},
      eprint={1010.3699},
}

\bib{blz}{article}{
      author={Bazhanov, V.},
      author={Lukyanov, S.},
      author={Zamolodchikov, A.},
       title={{Higher-level eigenvalues of Q-operators and Schrodinger
  equation}},
        date={2003},
        ISSN={1095-0753},
     journal={Advances in Theoretical and Mathematical Physics},
      volume={7},
      number={4},
       pages={711\ndash 725},
      eprint={hep-th/0307108},
         url={http://dx.doi.org/10.4310/ATMP.2003.v7.n4.a4},
}

\bib{Bazhanov:1996dr}{article}{
      author={Bazhanov, Vladimir~V.},
      author={Lukyanov, Sergei~L.},
      author={Zamolodchikov, Alexander~B.},
       title={{Integrable Structure of Conformal Field Theory. 2. Q operator
  and DDV Equation}},
        date={1997},
     journal={Commun. Math. Phys.},
      volume={190},
       pages={247\ndash 278},
      eprint={hep-th/9604044},
}

\bib{Bazhanov:1998dq}{article}{
      author={Bazhanov, Vladimir~V.},
      author={Lukyanov, Sergei~L.},
      author={Zamolodchikov, Alexander~B.},
       title={{Integrable Structure of Conformal Field Theory. 3. The
  Yang-Baxter Relation}},
        date={1999},
     journal={Commun. Math. Phys.},
      volume={200},
       pages={297\ndash 324},
      eprint={hep-th/9805008},
}

\bib{Zabrodin:}{article}{
      author={Beketov, M.},
      author={Liashuk, A.},
      author={Zabrodin, A.},
      author={Zotov, A.},
       title={{Trigonometric version of quantum-classical duality}},
        date={2016},
     journal={Nucl. Phys. B},
      volume={903},
      eprint={1201.3990},
         url={https://arxiv.org/abs/1501.07509},
}

\bib{Braverman:2010ei}{article}{
      author={Braverman, Alexander},
      author={Maulik, Davesh},
      author={Okounkov, Andrei},
       title={{Quantum Cohomology of the Springer Resolution}},
        date={2010},
      eprint={1001.0056},
}

\bib{Brinson:2021ww}{article}{
      author={Brinson, Ty~J.},
      author={Sage, Daniel~S.},
      author={Zeitlin, Anton~M.},
       title={{Opers on the Projective line, Wronskian Relations, and the Bethe
  Ansatz}},
        date={202112},
      eprint={2112.02711},
         url={https://arxiv.org/pdf/2112.02711.pdf},
}

\bib{Cherednik:1991mg}{article}{
      author={Cherednik, Ivan},
       title={{Integration of quantum many body problems by Affine
  Knizhnik-Zamolodchikov equations}},
        date={1994},
     journal={Adv. Math.},
      volume={106},
       pages={65\ndash 95},
}

\bib{Ciocan-Fontanine:2011tg}{article}{
      author={Ciocan-Fontanine, I.},
      author={Kim, B.},
      author={Maulik, D.},
       title={Stable quasimaps to git quotients},
        date={2014},
     journal={J. Geom. Phys.},
      volume={75},
       pages={17\ndash 47},
}

\bib{CG}{book}{
      author={Chriss, N.},
      author={Ginzburg, V.},
       title={{Representation Theory and Complex Geometry}},
   publisher={Birkh\"{a}user Boston, Inc., Boston, MA},
        date={1997},
        ISBN={0-8176-3792-3},
}

\bib{Chari:1994pz}{book}{
      author={Chari, V.},
      author={Pressley, A.},
       title={{A guide to quantum groups}},
   publisher={CUP},
        date={1994},
}

\bib{Dinkins:2020aa}{article}{
      author={Dinkins, Hunter},
       title={{3d Mirror Symmetry of the Cotangent Bundle of the Full Flag
  Variety}},
        date={2022},
     journal={Lett. Math. Phys.},
      volume={122},
       pages={1\ndash 31},
      eprint={2011.08603},
         url={https://arxiv.org/pdf/2011.08603.pdf},
}

\bib{Dinkins:2020ab}{article}{
      author={Dinkins, Hunter},
       title={{Symplectic Duality of $T^*Gr(k,n)$}},
        date={2022},
     journal={Math. Res. Lett.},
      volume={29},
      number={3},
       pages={663\ndash 690},
      eprint={2008.05516},
         url={https://arxiv.org/pdf/2008.05516.pdf},
}

\bib{Drinfeld:1986in}{article}{
      author={Drinfeld, V.~G.},
       title={Quantum groups},
        date={1988},
     journal={J. Math. Sci.},
      volume={41},
       pages={898},
}

\bib{Dorey:2007zx}{article}{
      author={Dorey, P.},
      author={Dunning, C.},
      author={Tateo, R.},
       title={{The ODE/IM Correspondence}},
        date={2007},
     journal={J.Phys.},
      volume={A40},
       pages={R205},
      eprint={hep-th/0703066},
}

\bib{DF}{article}{
      author={Ding, J.T.},
      author={Frenkel, I.B},
       title={{Isomorphism of two realizations of quantum affine algebra
  $U_q(\widehat{\mathfrak{gl}}(n))$}},
        date={1993},
     journal={Commun. Math. Phys.},
      volume={156},
       pages={277\ndash 300},
}

\bib{dedushenko2023interfaces}{article}{
      author={Dedushenko, M.},
      author={Nekrasov, N.},
       title={{Interfaces and Quantum Algebras, I: Stable Envelopes}},
      eprint={2109.10941},
}

\bib{Drinfeld:1985}{article}{
      author={Drinfeld, V.},
      author={Sokolov, V.},
       title={{Lie Algebras and Equations of Korteweg-de Vries Type}},
        date={1985},
     journal={J. Sov. Math.},
      volume={30:2},
       pages={1975\ndash 2036},
}

\bib{etingof1998lectures}{book}{
      author={Etingof, Pavel~I.},
      author={Frenkel, Igor~B.},
      author={Kirillov, Alexander~A.},
       title={{Lectures on representation theory and Knizhnik-Zamolodchikov
  equations}},
   publisher={Math. Surv. Mon. AMS},
        date={1998},
      number={58},
}

\bib{ekhammar2021extended}{article}{
      author={Ekhammar, S.},
      author={Shu, H.},
      author={Volin, D.},
       title={{Extended systems of Baxter Q-functions and fused flags I:
  simply-laced case}},
      eprint={2008.10597},
}

\bib{Faddeev:aa}{article}{
      author={Faddeev, L.~D.},
       title={{How Algebraic Bethe Ansatz works for integrable model}},
      eprint={hep-th/9605187},
         url={https://arxiv.org/pdf/hep-th/9605187.pdf},
}

\bib{Frenkel:2003qx}{article}{
      author={Frenkel, E.},
       title={{Opers on the Projective Line, flag Manifolds and Bethe Ansatz}},
        date={2003},
     journal={Mosc. Math. J.},
      volume={4},
       pages={655\ndash 705},
      eprint={math/0308269},
}

\bib{Frenkel:2004qy}{article}{
      author={Frenkel, E.},
       title={{Gaudin Model and Opers}},
        date={2004},
     journal={{Infinite Dimensional Algebras and Quantum Integrable Systems,
  eds. P. Kulish, e.a., Progress in Math.}},
      volume={237},
       pages={1\ndash 60},
      eprint={math/0407524},
}

\bib{frenkel2007langlands}{book}{
      author={Frenkel, Edward},
       title={{Langlands Correspondence for Loop Groups}},
   publisher={Cambridge University Press},
        date={2007},
      volume={103},
}

\bib{Frenkel_LanglandsLoop}{book}{
      author={Frenkel, Edward},
       title={Langlands correspondence for loop groups},
   publisher={Cambridge University Press},
        date={2007},
}

\bib{Feigin:1994in}{article}{
      author={Feigin, B.},
      author={Frenkel, E.},
      author={Reshetikhin, N.},
       title={{Gaudin Model, Bethe Ansatz and Critical Level}},
        date={1994},
     journal={Commun. Math. Phys.},
      volume={166},
       pages={27\ndash 62},
      eprint={hep-th/9402022},
}

\bib{Rybnikov:2010}{article}{
      author={Feigin, B.},
      author={Frenkel, E.},
      author={Rybnikov, L.},
       title={{Opers with Irregular Singularity and Spectra of the Shift of
  Argument Subalgebra}},
        date={2010},
     journal={Duke Math. J.},
      volume={155},
       pages={337\ndash 363},
}

\bib{FFTL:2010}{article}{
      author={Feigin, B.},
      author={Frenkel, E.},
      author={Toledano~Laredo, V.},
       title={{Gaudin models with irregular singularities}},
        date={2010},
     journal={Adv. Math.},
      volume={223},
       pages={873\ndash 948},
}

\bib{frenkel2023extended}{article}{
      author={Frenkel, E.},
      author={Hernandez, D.},
       title={{Extended Baxter relations and $QQ$-systems for quantum affine
  algebras}},
      eprint={2312.13256},
}

\bib{Frenkel:2013uda}{article}{
      author={Frenkel, E.},
      author={Hernandez, D.},
       title={{Baxter's relations and spectra of quantum integrable models}},
        date={2015},
     journal={Duke Math. J.},
      volume={164},
      number={12},
       pages={2407\ndash 2460},
      eprint={1308.3444},
}

\bib{Frenkel:2016}{article}{
      author={Frenkel, E.},
      author={Hernandez, D.},
       title={{Spectra of Quantum KdV Hamiltonians, Langlands Duality, and
  Affine Opers}},
        date={2018},
     journal={Commun. Math. Phys.},
      volume={362},
       pages={361\ndash 414},
      eprint={1606.05301},
         url={https://arxiv.org/abs/1606.05301},
}

\bib{FHRnew}{article}{
      author={Frenkel, E.},
      author={Hernandez, D.},
      author={Reshetikhin, N.},
       title={{Folded quantum integrable models and deformed W-algebras}},
        date={2022},
     journal={Lett. Math. Phys.},
      volume={112},
       pages={80},
      eprint={2110.14600},
}

\bib{Frenkel:2020}{article}{
      author={Frenkel, Edward},
      author={Koroteev, Peter},
      author={Sage, Daniel~S.},
      author={Zeitlin, Anton~M.},
       title={{q-Opers, QQ-Systems, and Bethe Ansatz}},
        date={2020},
     journal={J. Europ. Math. Soc.},
       pages={to appear},
      eprint={2002.07344},
         url={https://arxiv.org/pdf/2002.07344.pdf},
}

\bib{FR1998}{article}{
      author={Frenkel, Igor},
      author={Reshetikhin, Nicolai},
       title={Quantum affine algebras and holonomic difference equations},
        date={1992},
     journal={Comm. Math. Phys.},
      volume={146},
       pages={1\ndash 60},
}

\bib{Frenkel1998}{article}{
      author={Frenkel, E.},
      author={Reshetikhin, N.},
      author={Semenov-Tian-Shansky, M.},
       title={{Drinfeld--{S}okolov Reduction for Difference Operators and
  Deformations of $\mathcal{W}$-Algebras I: The Case of {V}irasoro Algebra}},
        date={1998},
        ISSN={1432-0916},
     journal={Commun. Math. Phys.},
      volume={192},
       pages={605\ndash 629},
      eprint={q-alg/9704011},
         url={https://doi.org/10.1007/s002200050311},
}

\bib{frt}{article}{
      author={Faddeev, L.D.},
      author={Reshetikhin, N.Yu.},
      author={Takhtadzhyan, L.A.},
       title={{Quantization of Lie groups and Lie algebras}},
        date={1990},
     journal={Leningrad Math. J.},
      volume={1:1},
       pages={193\ndash 225},
}

\bib{faddeev1987hamiltonian}{book}{
      author={Faddeev, Ludwig~D.},
      author={Takhtajan, Leon~A.},
       title={{Hamiltonian methods in the theory of solitons}},
   publisher={Springer},
        date={1987},
}

\bib{FZ}{article}{
      author={Fomin, S.},
      author={Zelevinsky, A.},
       title={{Double Bruhat Cells and Total Positivity}},
        date={1999},
        ISSN={0894-0347},
     journal={J. Amer. Math. Soc.},
      volume={12},
      number={2},
       pages={335\ndash 380},
      eprint={math/9802056},
  url={https://doi-org.libezp.lib.lsu.edu/10.1090/S0894-0347-99-00295-7},
}

\bib{Ginzburg:}{article}{
      author={Ginzburg, V.},
       title={{Lectures on Nakajima's Quiver Varieties}},
      eprint={0905.0686},
         url={https://arxiv.org/abs/0905.0686},
}

\bib{GGKM}{article}{
      author={Gardner, Clifford~S.},
      author={Greene, John~M.},
      author={Kruskal, Martin~D.},
      author={Miura, Robert~M.},
       title={{Method for Solving the Korteweg-deVries Equation}},
        date={1967},
     journal={Phys. Rev. Lett.},
      volume={19},
       pages={1095\ndash 1097},
}

\bib{Guest}{book}{
      author={Guest, Martin~A.},
       title={From quantum cohomology to integrable systems},
      series={Oxf. Grad. Texts Math.},
   publisher={OUP},
        date={2008},
      volume={15},
}

\bib{geiss}{article}{
      author={Geiss, C.},
      author={Hernandez, D.},
      author={Leclerc, B.},
       title={{Representations of shifted quantum affine algebras and cluster
  algebras I. The simply-laced case}},
      eprint={2401.04616},
}

\bib{givental1995}{article}{
      author={Givental, Alexander},
      author={Kim, Bumsig},
       title={{Quantum Cohomology of Flag Manifolds and Toda Lattices}},
        date={1995},
     journal={Communications in Mathematical Physics},
      volume={168},
      number={3},
       pages={609\ndash 641},
         url={http://projecteuclid.org/euclid.cmp/1104272492},
}

\bib{2001math8105G}{article}{
      author={Givental, A.},
      author={Lee, Y.~P.},
       title={{Quantum K-Theory on Flag Manifolds, Finite-Difference Toda
  Lattices and Quantum Groups}},
        date={2001},
      eprint={math/0108105},
}

\bib{Gorsky:2013xba}{article}{
      author={Gorsky, A.},
      author={Zabrodin, A.},
      author={Zotov, A.},
       title={{Spectrum of Quantum Transfer Matrices via Classical Many-Body
  Systems}},
        date={2014},
     journal={JHEP},
      volume={1401},
       pages={070},
      eprint={1310.6958},
}

\bib{haouzi2023new}{article}{
      author={Haouzi, Nathan},
       title={{A new realization of quantum algebras in gauge theory and
  Ramification in the Langlands program}},
      eprint={2311.04367},
}

\bib{Heisenberg:1928aa}{article}{
      author={Heisenberg, W.},
       title={Zur theorie des ferromagnetismus},
        date={1928},
     journal={Z. Phys.},
      volume={49},
       pages={619\ndash 636},
}

\bib{hernandez2022stable}{article}{
      author={Hernandez, David},
       title={{Stable maps, Q-operators and category $\mathcal{O}$}},
        date={2022},
     journal={Rep. Th. of AMS},
      volume={26},
      number={7},
       pages={179\ndash 210},
      eprint={1902.02843},
}

\bib{HJ}{article}{
      author={Hernandez, D.},
      author={Jimbo, M.},
       title={{Asymptotic Representations and Drinfeld Rational Fractions}},
        date={2012},
     journal={Comp. Math.},
      volume={148},
      number={5},
       pages={1593\ndash 1623},
      eprint={1104.1891},
}

\bib{Hanany:1996ie}{article}{
      author={Hanany, Amihay},
      author={Witten, Edward},
       title={{Type IIB superstrings, BPS monopoles, and three-dimensional
  gauge dynamics}},
        date={1997},
     journal={Nucl.Phys.},
      volume={B492},
       pages={152\ndash 190},
      eprint={hep-th/9611230},
}

\bib{Intriligator:1996ex}{article}{
      author={Intriligator, Kenneth~A.},
      author={Seiberg, N.},
       title={{Mirror symmetry in three-dimensional gauge theories}},
        date={1996},
     journal={Phys.Lett.},
      volume={B387},
       pages={513\ndash 519},
      eprint={hep-th/9607207},
}

\bib{Jimbo1985q}{article}{
      author={Jimbo, Michio},
       title={{A $q$-difference analogue of $U(g)$ and the Yang-Baxter
  equation}},
        date={1985},
     journal={Lett. Math. Phys.},
      volume={10},
       pages={63\ndash 69},
}

\bib{Korepin_1993}{article}{
      author={Korepin, V.},
      author={Bogoliubov, N.},
      author={Izergin, A.},
       title={{Quantum Inverse Scattering Method and Correlation Functions}},
        date={1993},
     journal={Cambridge University Press},
         url={http://dx.doi.org/10.1017/CBO9780511628832},
}

\bib{Kaledinreview}{incollection}{
      author={Kaledin, D.},
       title={Geometry and topology of symplectic resolutions},
        date={2009},
   booktitle={Algebraic geometry---{S}eattle 2005. {P}art 2},
      series={Proc. Sympos. Pure Math.},
      volume={80, Part 2},
   publisher={AMS},
       pages={595\ndash 628},
}

\bib{Kamnitzer:2022aa}{article}{
      author={Kamnitzer, Joel},
       title={Symplectic resolutions, symplectic duality, and coulomb
  branches},
        date={2022},
      eprint={2202.03913},
         url={https://arxiv.org/pdf/2202.03913.pdf},
}

\bib{1996alg.geom7001K}{article}{
      author={{Kim}, B.},
       title={{Quantum cohomology of flag manifolds $G/B$ and quantum Toda
  lattices}},
        date={1999},
     journal={Ann. Math.},
      volume={149},
      number={1},
       pages={129\ndash 148},
      eprint={alg-geom/9607001},
}

\bib{Kirquiv}{book}{
      author={Kirillov, Alexander, Jr.},
       title={Quiver representations and quiver varieties},
      series={Grad. Stud. Math.},
   publisher={AMS},
        date={2016},
      volume={174},
}

\bib{kontsevich1992intersection}{article}{
      author={Kontsevich, Maxim},
       title={Intersection theory on the moduli space of curves and the matrix
  airy function},
        date={1992},
     journal={Communications in Mathematical Physics},
      volume={147},
       pages={1\ndash 23},
}

\bib{Koroteev:2017aa}{article}{
      author={Koroteev, P.},
      author={Pushkar, P.},
      author={Smirnov, A.},
      author={Zeitlin, A.},
       title={{Quantum {K}-theory of Quiver Varieties and Many-Body Systems}},
        date={2021},
     journal={Selecta Math.},
      volume={27},
      number={87},
      eprint={1705.10419},
         url={https://arxiv.org/pdf/1705.10419},
}

\bib{KSZ}{article}{
      author={Koroteev, P.},
      author={Sage, D.S.},
      author={Zeitlin, A.M.},
       title={{(SL(N),q)-Opers, the q-Langlands Correspondence, and
  Quantum/Classical Duality}},
        date={2021},
     journal={Commun. Math. Phys.},
      volume={381},
      number={2},
       pages={641\ndash 672},
      eprint={1811.09937},
         url={https://arxiv.org/abs/1508.02690},
}

\bib{KhoroshkinR}{article}{
      author={Khoroshkin, S.~M.},
      author={Tolstoy, V.~N.},
       title={Universal {$R$}-matrix for quantized (super)algebras},
        date={1991},
     journal={Comm. Math. Phys.},
      volume={141},
      number={3},
       pages={599\ndash 617},
}

\bib{Khoroshkindrinf}{incollection}{
      author={Khoroshkin, S.~M.},
      author={Tolstoy, V.~N.},
       title={On {D}rinfel\cprime d's realization of quantum affine algebras},
        date={1993},
      volume={11},
       pages={445\ndash 452},
        note={Infinite-dimensional geometry in physics (Karpacz, 1992)},
}

\bib{MR3199543}{article}{
      author={Khoroshkin, Sergey},
      author={Tsuboi, Zengo},
       title={The universal {$R$}-matrix and factorization of the {L}-operators
  related to the {B}axter {Q}-operators},
        date={2014},
        ISSN={1751-8113},
     journal={J. Phys. A},
      volume={47},
      number={19},
       pages={192003, 11},
         url={http://dx.doi.org/10.1088/1751-8113/47/19/192003},
      review={\MR{3199543}},
}

\bib{koroteev2022zoo}{article}{
      author={Koroteev, Peter},
      author={Zeitlin, Anton~M},
       title={{The Zoo of Opers and Dualities}},
     journal={to appear in Int. Math. Res. Not.},
      eprint={arXiv:2208.08031},
         url={https://arxiv.org/abs/2208.08031},
}

\bib{Koroteev:2018a}{article}{
      author={Koroteev, Peter},
      author={Zeitlin, Anton~M.},
       title={{qKZ/tRS Duality via Quantum K-Theoretic Counts}},
        date={2021},
     journal={Math. Res. Lett.},
      volume={28,},
      number={2},
       pages={435\ndash 470},
      eprint={1802.04463},
         url={https://arxiv.org/pdf/1802.04463.pdf},
}

\bib{Koroteev:2023aa}{article}{
      author={Koroteev, Peter},
      author={Zeitlin, Anton~M.},
       title={{3D Mirror Symmetry for Instanton Moduli Spaces}},
        date={2023},
     journal={Commun. Math. Phys.},
       pages={1005\ndash 1068},
      eprint={2105.00588},
         url={https://doi.org/10.1007/s00220-023-04831-5},
}

\bib{KoroteevZeitlinCrelle}{article}{
      author={Koroteev, Peter},
      author={Zeitlin, Anton~M.},
       title={{q-Opers, QQ-Systems, and Bethe Ansatz II: Generalized Minors}},
        date={2023},
     journal={J. Reine Angew. Math. (Crelles Journal)},
      volume={2023},
      number={795},
       pages={271\ndash 296},
      eprint={2108.04184},
         url={https://doi.org/10.1515/crelle-2022-0084},
}

\bib{koroteev_zeitlin_2023}{article}{
      author={Koroteev, Peter},
      author={Zeitlin, Anton~M.},
       title={{Toroidal q-Opers}},
        date={2023},
     journal={J. Inst. Math. Jussieu},
      volume={22},
      number={2},
       pages={581\ndash 642},
      eprint={2007.11786},
}

\bib{2012arXiv1210.2315L}{article}{
      author={{Li}, J.~R.},
      author={{Tarasov}, V.},
       title={{XXZ-type Bethe Ansatz Equations and Quasi-Polynomials}},
        date={2013},
     journal={Journal of Physics: Conference Series},
      volume={411},
      number={1},
       pages={012020},
      eprint={1210.2315},
}

\bib{Maninfrob}{book}{
      author={Manin, Yuri~I.},
       title={Frobenius manifolds, quantum cohomology, and moduli spaces},
      series={Amer. Math. Soc. Coll. Publ.},
   publisher={AMS},
        date={1999},
      volume={47},
}

\bib{Matsuo1992}{article}{
      author={Matsuo, Atushi},
       title={Integrable connections related to zonal spherical functions.},
        date={1992},
     journal={Inv. Math.},
      volume={110},
      number={1},
       pages={95\ndash 122},
         url={http://eudml.org/doc/144044},
}

\bib{McGerty:2016kir}{article}{
      author={McGerty, Kevin},
      author={Nevins, Thomas},
       title={Kirwan surjectivity for quiver varieties},
        date={2018},
     journal={Inv. Math.},
      volume={212},
       pages={161\ndash 187},
      eprint={1610.08121},
}

\bib{2012arXiv1211.1287M}{article}{
      author={{Maulik}, D.},
      author={{Okounkov}, A.},
       title={{Quantum Groups and Quantum Cohomology}},
        date={2019},
     journal={Ast\'erisque},
      volume={408},
      eprint={1211.1287},
}

\bib{MRV1}{article}{
      author={Masoero, D.},
      author={Raimondo, A.},
      author={Valeri, D.},
       title={{Bethe Ansatz and the Spectral Theory of Affine Lie
  Algebra-Valued Connections I. The Simply-Laced Case}},
        date={2016},
        ISSN={1432-0916},
     journal={Commun. Math. Phys.},
      volume={344},
      number={3},
       pages={719\ndash 750},
      eprint={1501.07421},
         url={http://dx.doi.org/10.1007/s00220-016-2643-6},
}

\bib{MRV2}{article}{
      author={Masoero, D.},
      author={Raimondo, A.},
      author={Valeri, D.},
       title={{Bethe Ansatz and the Spectral Theory of Affine Lie
  algebra--Valued Connections II: The Non Simply--Laced Case}},
        date={2017},
     journal={Commun. Math. Phys.},
      volume={349},
      number={3},
       pages={1063\ndash 1105},
      eprint={1511.00895},
         url={https://doi.org/10.1007/s00220-016-2744-2},
}

\bib{MVflag}{article}{
      author={Mukhin, E.},
      author={Varchenko, A.},
       title={{Critical points of master functions and flag varieties}},
        date={2004},
     journal={Comm. Contemp. Math.},
      volume={6},
      number={01},
       pages={111\ndash 163},
      eprint={math/0209017},
}

\bib{Mukhin_2005}{article}{
      author={Mukhin, E.},
      author={Varchenko, A.},
       title={{Discrete Miura Opers and Solutions of the Bethe Ansatz
  Equations}},
        date={2005},
        ISSN={1432-0916},
     journal={Commun. Math. Phys.},
      volume={256},
       pages={565\ndash 588},
      eprint={math/0401137},
         url={http://dx.doi.org/10.1007/s00220-005-1288-7},
}

\bib{mukhvarmiura}{article}{
      author={Mukhin, E.},
      author={Varchenko, A.},
       title={{Miura Opers and Critical Points of Master Functions}},
        date={2005},
     journal={Cent. Eur. J. Math.},
      volume={3},
       pages={155\ndash 182},
      eprint={math/0312406},
         url={https://arxiv.org/abs/math/0312406},
}

\bib{MR2189873}{article}{
      author={Mukhin, E.},
      author={Varchenko, A.},
       title={{Populations of solutions of the {$XXX$} {B}ethe equations
  associated to {K}ac-{M}oody algebras}},
        date={2005},
     journal={Contemp. Math.},
      volume={392},
      eprint={math/0212092},
}

\bib{2006math......4048M}{article}{
      author={{Mukhin}, E},
      author={{Varchenko}, A},
       title={{Quasi-polynomials and the Bethe Ansatz}},
        date={2008},
     journal={Geom. Top. Mon.},
      volume={13},
      eprint={arXiv:math/0604048},
}

\bib{Novikov:1984}{article}{
      author={Novikov, S.},
      author={Manakov, S.V.},
      author={Pitaievskii, L.P.},
      author={Zakharov, V.E.},
       title={{Theory of Solitons: The Inverse Scattering Method}},
        date={1984},
     journal={Monographs in Contemporary Mathematics},
}

\bib{Nakajima:1999hilb}{book}{
      author={Nakajima, H.},
       title={{Lectures on Hilbert schemes of points on surfaces}},
   publisher={AMS},
        date={1999},
}

\bib{Nakajima:2001qg}{article}{
      author={Nakajima, H.},
       title={{Quiver Varieties and Finite Dimensional Representations of
  Quantum Affine Algebras}},
        date={2001},
     journal={J. Amer. Math. Soc.},
      volume={14},
       pages={145\ndash 238},
      eprint={math/9912158},
}

\bib{nakajima1998}{article}{
      author={Nakajima, Hiraku},
       title={Quiver varieties and kac-moody algebras},
        date={199802},
     journal={Duke Math. J.},
      volume={91},
      number={3},
       pages={515\ndash 560},
         url={http://dx.doi.org/10.1215/S0012-7094-98-09120-7},
}

\bib{Nekrasov:2009ui}{article}{
      author={Nekrasov, N.},
      author={Shatashvili, S.},
       title={{Quantum Integrability and Supersymmetric Vacua}},
        date={2009},
     journal={Prog. Theor. Phys. Suppl.},
      volume={177},
       pages={105\ndash 119},
      eprint={0901.4748},
}

\bib{Nekrasov:2009uh}{article}{
      author={Nekrasov, N.},
      author={Shatashvili, S.},
       title={{Supersymmetric Vacua and Bethe Ansatz}},
        date={2009},
     journal={Nucl. Phys. Proc. Suppl.},
      volume={192-193},
       pages={91\ndash 112},
      eprint={0901.4744},
}

\bib{Oblomkov2004}{article}{
      author={Oblomkov, Alexei},
       title={{Double Affine Hecke Algebras and Calogero-Moser Spaces}},
        date={2004},
     journal={Rep. Th. of AMS},
      volume={8},
      number={10},
       pages={243\ndash 266},
         url={https://arxiv.org/abs/math/0303190},
}

\bib{Okounkov:2015aa}{article}{
      author={Okounkov, A.},
       title={{Lectures on {K}-Theoretic Computations in Enumerative
  Geometry}},
        date={2015},
      eprint={1512.07363},
         url={https://arxiv.org/abs/1512.07363},
}

\bib{Okounkov:2016sya}{article}{
      author={Okounkov, Andrei},
      author={Smirnov, Andrey},
       title={{Quantum Difference Equation for Nakajima Varieties}},
        date={2022},
     journal={Invent. Math.},
      number={229},
       pages={1203\ndash 1299},
      eprint={1602.09007},
}

\bib{Pushkar:2016qvw}{article}{
      author={Pushkar, Petr},
      author={Smirnov, Andrey},
      author={Zeitlin, Anton},
       title={{Baxter Q-operator from quantum K-theory}},
        date={2020},
     journal={Adv. Math.},
      volume={360},
       pages={106919},
      eprint={1612.08723},
}

\bib{Reshetikhin:2010si}{article}{
      author={Reshetikhin, N.},
       title={{Lectures on the Integrability of the 6-Vertex Model}},
        date={2010},
      eprint={1010.5031},
         url={https://arxiv.org/abs/1010.5031},
}

\bib{Rimanyi:2019ab}{article}{
      author={Rim{\'a}nyi, Rich{\'a}rd},
      author={Smirnov, Andrey},
      author={Varchenko, Alexander},
      author={Zhou, Zijun},
       title={{Three-Dimensional Mirror Self-Symmetry of the Cotangent Bundle
  of the Full Flag Variety}},
        date={2019},
     journal={SIGMA},
      volume={15},
      number={093},
      eprint={1906.00134},
         url={https://arxiv.org/pdf/1906.00134.pdf},
}

\bib{Rimanyi:2014ef}{article}{
      author={Rimanyi, R.},
      author={Tarasov, V.},
      author={Varchenko, A.},
       title={Trigonometric weight functions as {K}-theoretic stable envelope
  maps for the cotangent bundle of a flag variety},
        date={201411},
      eprint={1411.0478},
         url={https://arxiv.org/abs/1411.0478},
}

\bib{olivier1998quantum}{article}{
      author={Schiffmann, Olivier},
       title={Quantum affine algebras at roots of unity and equivariant
  k-theory},
        date={1998},
      eprint={math/9808038},
}

\bib{1998CMaPh.192..631S}{article}{
      author={{Semenov-Tian-Shansky}, M.},
      author={{Sevostyanov}, A.},
       title={{Drinfeld-Sokolov Reduction for Difference Operators and
  Deformations of $\mathcal{W}$-Algebras II. The General Semisimple Case}},
        date={1998},
     journal={Commun. Math. Phys.},
      volume={192},
       pages={631\ndash 647},
      eprint={q-alg/9702016},
}

\bib{Smirnov:2020aa}{article}{
      author={Smirnov, Andrey},
      author={Zhou, Zijun},
       title={{3d Mirror Symmetry and Quantum $K$-theory of Hypertoric
  Varieties}},
        date={2022},
     journal={Adv. Math.},
      volume={395},
       pages={108081},
      eprint={2006.00118},
         url={https://arxiv.org/pdf/2006.00118.pdf},
}

\bib{Takhtajan:1979iv}{article}{
      author={Takhtajan, L.~A.},
      author={Faddeev, L.~D.},
       title={{The Quantum method of the inverse problem and the Heisenberg XYZ
  model}},
        date={1979},
     journal={Russ. Math. Surveys},
      volume={34},
      number={5},
       pages={11\ndash 68},
        note={[Usp. Mat. Nauk34,no.5,13(1979)]},
}

\bib{Tarasov_2002}{article}{
      author={Tarasov, V.},
      author={Varchenko, A.},
       title={{Duality for Knizhnik--Zamolodchikov and Dynamical Equations}},
        date={2002},
        ISSN={0167-8019},
     journal={Acta Appl. Math.},
      volume={73},
      number={1/2},
       pages={141\ndash 154},
         url={http://dx.doi.org/10.1023/A:1019787006990},
}

\bib{Tarasov_2005}{article}{
      author={Tarasov, V},
      author={Varchenko, A.},
       title={{Dynamical Differential Equations Compatible with Rational qKZ
  Equations}},
        date={2005},
        ISSN={1573-0530},
     journal={Lett. Math. Phys.},
      volume={71},
      number={2},
       pages={101\ndash 108},
         url={http://dx.doi.org/10.1007/s11005-004-6363-z},
}

\bib{varagnolo2000quiver}{article}{
      author={Varagnolo, Michela},
       title={{Quiver varieties and Yangians}},
        date={2000},
     journal={Lett. Math. Phys.},
      volume={53},
       pages={273\ndash 283},
      eprint={math/0005277},
}

\bib{Vasserot:wo}{article}{
      author={Vasserot, Eric},
       title={{Affine Quantum Groups and Equivariant K-theory}},
        date={1998},
     journal={Transform. Groups},
      volume={3},
      number={3},
       pages={269\ndash 299},
      eprint={math/9803024},
}

\bib{witten1990two}{article}{
      author={Witten, Edward},
       title={Two-dimensional gravity and intersection theory on moduli space},
        date={1990},
     journal={Surveys in differential geometry},
      volume={1},
      number={1},
       pages={243\ndash 310},
}

\bib{ZahFAd}{article}{
      author={Zaharov, V.~E.},
      author={Faddeev, L.~D.},
       title={{The {K}orteweg-de {V}ries equation is a fully integrable
  {H}amiltonian system}},
        date={1971},
     journal={Funkcional. Anal. i Prilo\v{z}en.},
      volume={5},
      number={4},
       pages={18\ndash 27},
}

\bib{zeitlin2022wronskians}{article}{
      author={Zeitlin, Anton~M.},
       title={{On Wronskians and $qq$-systems}},
        date={2024},
     journal={Contemp. Math.; ``Algebraic and Topological Aspects of
  Representation Theory"},
      volume={791},
      eprint={2208.08018},
}

\bib{Zabrodin:2017td}{article}{
      author={Zabrodin, A.},
      author={Zotov, A.},
       title={{{KZ}-Calogero Correspondence Revisited}},
        date={2017},
     journal={J. Phys. A: Math. Theor.},
      volume={50},
       pages={205202},
      eprint={1701.06074},
         url={https://arxiv.org/pdf/1701.06074.pdf},
}

\bib{Zabrodin:2017vt}{article}{
      author={Zabrodin, A.},
      author={Zotov, A.},
       title={{{QKZ}-Ruijsenaars Correspondence Revisited}},
        date={2017},
     journal={Nuclear Physics B,},
      volume={922},
       pages={113\ndash 125},
      eprint={1704.04527},
         url={https://arxiv.org/pdf/1704.04527.pdf},
}

\end{biblist}
\end{bibdiv}

\end{document}